\documentclass[reqno]{amsart}
\usepackage{amssymb,latexsym,amsmath,amsthm,enumerate}
\usepackage[mathscr]{eucal}
\usepackage{mathrsfs}
\usepackage{graphicx}
\usepackage{cite}
\makeatletter
\@namedef{subjclassname@2020}{\textup{2020} Mathematics Subject Classification}
\makeatother
\theoremstyle{plain}
\newtheorem{thm}{Theorem}[section]
\newtheorem{lem}[thm]{Lemma}
\newtheorem{cor}[thm]{Corollary}
\newtheorem{pro}[thm]{Proposition}

\theoremstyle{definition}

\newtheorem{remark}[thm]{Remark}
\newcommand{\Var}{\mathsf{V}}
\newcommand{\Mod}{\operatorname{Mod}}

\newcommand{\init}{\operatorname{in}}

\newcommand{\bu}{\mathbf{u}}
\newcommand{\bv}{\mathbf{v}}
\newcommand{\bw}{\mathbf{w}}
\newcommand{\bp}{\mathbf{p}}
\newcommand{\bq}{\mathbf{q}}
\newcommand{\br}{\mathbf{r}}
\newcommand{\bs}{\mathbf{s}}
\newcommand{\bt}{\mathbf{t}}
\renewcommand{\le}{\leqslant}
\renewcommand{\ge}{\geqslant}
\renewcommand{\leq}{\leqslant}
\renewcommand{\geq}{\geqslant}

\begin{document}
\title[A smallest generator of a limit variety]
{An additively idempotent semiring of least order generating a limit variety}

\author{Qirui Ma}
\address{School of Mathematics, Northwest University, Xi'an, 710127, Shaanxi, P.R. China}
\email{qiruimath@yeah.net}

\author{Mengya Yue}
\address{School of Mathematics, Northwest University, Xi'an, 710127, Shaanxi, P.R. China}
\email{myayue@yeah.net}


\begin{abstract}
We construct an additively idempotent semiring of order five that generates
 a limit variety and is not strongly nonfinitely based. The generated variety
 has exactly two maximal proper subvarieties, both hereditarily finitely based.
 We also prove that no ai-semiring of order at most four generates a limit
 variety. Consequently, five is the least possible order of a finite
 ai-semiring generating a limit variety.
\end{abstract}
\maketitle
\section{Introduction}

An \emph{additively idempotent semiring} (or \emph{ai-semiring} for short)
is an algebra $(S,+,\cdot)$ such that $(S,+)$ is a semilattice,
$(S,\cdot)$ is a semigroup, and multiplication distributes over addition
on both sides. Thus an ai-semiring satisfies the identities
\begin{gather}
 x+y\approx y+x,\qquad (x+y)+z\approx x+(y+z),\qquad x+x\approx x,
 \label{eq:ai-add}\\
 (xy)z\approx x(yz),\qquad x(y+z)\approx xy+xz,
 \qquad (x+y)z\approx xz+yz.\label{eq:ai-mult}
\end{gather}
Throughout this paper, the signature consists of the two binary operation
symbols $+$ and $\cdot$, without constant symbols. A zero or an identity
element for multiplication is therefore not required.

The finite basis problem asks whether all identities of a given algebra follow
from finitely many of them. For additively idempotent semirings, this problem
depends on the interaction between the two operations. Nonfinite basis results
for finite semirings include Dolinka's example \cite{Dolinka2007}, the Brandt
monoid \cite{Volkov2021}, and examples whose multiplicative reducts are
finitely based \cite{Jackson2022}. Related results concern power semirings
\cite{GusevVolkov2023}, flat nil-semirings \cite{Wu2023}, and the semiring
$S_7^0$ \cite{Wu2024}.

A \emph{limit variety} is a nonfinitely based variety all of whose proper
subvarieties are finitely based. Limit varieties
arising from flat extensions of groups were studied in \cite{RenJackson2023}.
The preprints \cite{Lyu2026,Shao2026} give six-element examples denoted by
$SR_6$ and $TR_6$, respectively. In the present paper we give a five-element
example and prove that its order is least possible. The proof describes identities directly and treats arbitrary
subvarieties, including subvarieties not specified by a finite generating
algebra.

The study of small additively idempotent semirings provides useful finite
factors. The varieties generated by the two-element semirings were studied
in \cite{ShaoRen2015}, and the three-element case in \cite{Zhao2020}.
Several classes of four-element semirings are treated in
\cite{Ren2025I,Yue2025,Ren2026III}. We use a particular finite basis from
\cite{Yue2025}, after giving the correspondence between the two labellings
of its algebra.
For the minimum-order assertion, we use the complete four-element
classification in the preprint of Wang, Wang, Yin and Wu~\cite{Wang2026}.
We also give the finite constructions needed to exclude each of its
nonfinitely based ai-semirings from generating a limit variety.

The multiplicatively idempotent case has a finite subvariety lattice.
Its structure is described through ordered bands and semilattice-ordered
semigroups in \cite{Ghosh2005,Pastijn2005,KurilPolak2005,PastijnZhao2005};
see also \cite{Pastijn2006}. We use the explicit description in
\cite{Polak2006}. Another input is the hereditary finite basis theorem of
Ren and Zeng \cite{RenZeng2019}. Quotients defined by languages, related to
the syntactic semiring construction \cite{Polak2001}, are useful for
separating identities. The quotient construction and every separation
argument needed here are given explicitly.

The natural order on an ai-semiring is defined by
$a\leq_S b$ if $a+b=b$. For terms, $\bp\preceq\bq$ abbreviates the
identity $\bp+\bq\approx\bq$. Both operations preserve the natural order
in each argument.

For a class $\mathcal A$ of algebras, $\Var(\mathcal A)$ denotes the
variety it generates, and $\operatorname{Id}(\mathcal A)$ its identities.
Joins and meets of varieties are denoted by $\vee$ and $\cap$.
We use $\operatorname{Mod}_{\rm ai}(\Sigma)$, or simply
$\operatorname{Mod}(\Sigma)$ when the background is fixed, for the
ai-semirings satisfying $\Sigma$. The standard facts about varieties,
free algebras, and equational deductions used below can be found in
\cite{BurrisSankappanavar1981}.

A variety is \emph{finitely based} if it admits a finite equational basis.
An algebra is finitely based if the variety it generates is finitely based.
Since the ai-semiring axioms form a finite set, finite basability
relative to these axioms is equivalent to finite basability in the signature
$\{+,\cdot\}$. A variety is \emph{hereditarily finitely based} if all its
subvarieties are finitely based.

A finite ai-semiring $A$ is \emph{strongly nonfinitely based} if every
finite ai-semiring $B$ such that $A\in\Var(B)$ is nonfinitely based;
see \cite{Shao2026}. To prove that $A$ is not strongly nonfinitely based,
it therefore suffices to find a finite, finitely based ai-semiring $B$
such that $A\in\Var(B)$. This condition concerns finite generating
algebras, whereas the definition of a limit variety concerns proper
subvarieties.

The main result of this paper is the following.

\begin{thm}\label{thm:main}
The five-element ai-semiring $S$ defined in Section~\ref{sec:nfb}
generates a limit variety. More precisely, $\Var(S)$ is nonfinitely based
and has exactly two maximal proper subvarieties, both of which are
hereditarily finitely based.
Moreover, $S$ is not strongly nonfinitely based: there is a finitely
based ai-semiring $H_5$ of order five such that $S$ is a homomorphic
image of an eight-element subsemiring of $H_5^2$.
Finally, no ai-semiring of order at most four generates a limit variety.
Thus $S$ has the least possible order among finite ai-semirings
generating limit varieties.
\end{thm}

Section~\ref{sec:nfb} gives the operation tables, an exact criterion for
identities, and a proof of nonfinite basability. Section~\ref{sec:foundations}
determines the two maximal proper subvarieties, denoted by $\mathcal C$
and $\mathcal E$, and gives them finite bases. After recording the required
multiplicatively idempotent varieties, we prove hereditary finite basability
of two subvarieties of $\mathcal C$. The final separation arguments extend
these results to $\mathcal C$ and $\mathcal E$ themselves.
In Section~\ref{sec:not-strong}, we give an equational basis for $H_5$
and construct the required homomorphism onto $S$.
Section~\ref{sec:minimality} establishes the minimum-order assertion.
The proof of the main theorem is completed in Section~\ref{sec:conclusion}.

\section{Preliminaries}\label{sec:preliminaries}

Let $X$ be a countably infinite set of variables. We write $X^+$ for
the free semigroup on $X$ and $X^*=X^+\cup\{\varepsilon\}$, where
$\varepsilon$ denotes the empty word. Words and terms are denoted by
bold lowercase letters, while individual variables are denoted by
ordinary lowercase letters. By distributivity and additive idempotence,
every ai-semiring term can be written as a finite nonempty sum
$\bu=\bu_1+\cdots+\bu_m$, where $\bu_i\in X^+$. We identify such a
term with its finite set of summands. Thus $\bw\in\bu$ means that
$\bw$ is a summand of $\bu$, and the same notation may be used in
set operations on summands.

For a word $\bw$, let $c(\bw)$ be its content, that is, the set of
variables occurring in it, and let $\ell(\bw)$ be its length.
We put $c(\varepsilon)=\varnothing$ and $\ell(\varepsilon)=0$.
For a term $\bu$, put $c(\bu)=\bigcup_{\bw\in\bu}c(\bw)$.
The symbols $h(\bw)$ and $t(\bw)$ denote the first and last letters
of a nonempty word, respectively. We write $e(\bw)$ for the set of
letters occurring in positions before its last position. In particular,
$e(x)=\varnothing$ for a variable $x$. For a term $\bu$, put
\[
 E(\bu)=\bigcup_{\bw\in\bu}e(\bw),\qquad
 T(\bu)=\{t(\bw):\bw\in\bu\}.
\]
The relation $\preceq$ is used only as an abbreviation for a formal
identity. Orders on elements and numerical inequalities are denoted
by $\leq$, with the algebra as a subscript when needed, and inclusion
of varieties is denoted by $\subseteq$.

We allow an empty left or right context around a word: this means omitting
the corresponding factor. Variables are always replaced by nonempty terms.
Thus the notation for contexts does not add an identity element or allow
substitution of an empty word.

\begin{lem}\label{lem:local-finite}
Let $A$ be a finite algebra. Every finitely generated free algebra in
$\Var(A)$ is finite. Suppose that $\mathcal U\subseteq\Var(A)$ is finitely
based and there is an integer $r$ such that, for every
$\mathcal W\subseteq\mathcal U$, all identities of $\mathcal W$ follow
from the identities of $\mathcal U$ and the identities of $\mathcal W$
in at most $r$ variables. Then $\mathcal U$ is hereditarily finitely based.
\end{lem}
\begin{proof}
An $r$-variable term induces a function $A^r\to A$, and there are at most
$|A|^{|A|^r}$ such functions. Two terms represent the same element of the
free algebra precisely when they induce the same function on $A$.
Thus this free algebra is finite. The corresponding free algebra of
any subvariety is its quotient and is also finite.

Choose a finite set of terms representing all $r$-variable term functions
modulo $\mathcal U$. Among the finitely many pairs of representatives,
retain those identities that hold in $\mathcal W$. These identities and
a finite basis of $\mathcal U$ imply every $r$-variable identity of
$\mathcal W$, hence, by hypothesis, all its identities. This is a finite
basis for $\mathcal W$.
\end{proof}
\section{A nonfinitely based ai-semiring of order five}
\label{sec:nfb}

Let $S=\{0,1,2,3,4\}$, with the operations
given in Table~\ref{tab:S}. The labels $0$ and $1$ are names of elements;
they are not distinguished constants in the signature.

\begin{table}[ht]
\caption{The Cayley tables of $S$}\label{tab:S}
\begin{tabular}{c|ccccc}
+&0&1&2&3&4\\\hline
0&0&0&0&0&0\\
1&0&1&0&0&1\\
2&0&0&2&2&2\\
3&0&0&2&3&3\\
4&0&1&2&3&4
\end{tabular}
\qquad
\begin{tabular}{c|ccccc}
$\cdot$&0&1&2&3&4\\\hline
0&0&0&0&0&0\\
1&0&0&0&0&0\\
2&0&0&2&2&4\\
3&0&1&2&3&4\\
4&4&4&4&4&4
\end{tabular}
\end{table}

\subsection{Verification of the semiring axioms}

\begin{lem}\label{lem:S-ai}
The algebra $S$ is an additively idempotent semiring.
\end{lem}
\begin{proof}
Consider the partial order determined by
\[
4<3<2<0,\qquad 4<1<0.
\]
The only incomparable pairs are $\{1,2\}$ and $\{1,3\}$, and the supremum
of each is $0$. The addition table is precisely the join operation of this
order. Hence addition is associative, commutative and idempotent. Its least
and greatest elements are $4$ and $0$, respectively.

For $s\in S$, let $L_s(t)=st$. The multiplication table gives
\[
L_0(t)=L_1(t)=0,\qquad L_3(t)=t,\qquad L_4(t)=4\quad(t\in S).
\]
The remaining map $L_2=f$ satisfies
\[
f(0)=f(1)=0,\qquad f(2)=f(3)=2,\qquad f(4)=4.
\]
In particular, $f^2=f$, and $f$ fixes $0$ and $4$. We verify
$L_{st}=L_s\circ L_t$. If $s\in\{0,1,4\}$, both sides are the corresponding
constant map. If $s=3$, the equality follows from $3t=t$. If $s=2$, the
three possibilities $t\in\{0,1\}$, $t\in\{2,3\}$ and $t=4$ give,
respectively, the constant map with value $0$, the map $f$, and the constant
map with value $4$, on both sides. Thus
\[
(st)u=L_{st}(u)=L_s(L_t(u))=s(tu)
\]
for all $s,t,u\in S$.

Left distributivity is equivalent to preservation of binary joins by each
$L_s$. Constant maps and the identity map preserve binary joins. The map
$f$ is order preserving on both displayed chains. For $j\in\{2,3\}$,
\[
f(1+j)=f(0)=0=f(1)+f(j),
\]
which also handles the incomparable pairs. Consequently $f$ preserves joins.

To verify right distributivity, put $R_t(s)=st$. The images of the long
chain and of the remaining element $1$ are
\[
\begin{array}{c|c|c}
t&(R_t(4),R_t(3),R_t(2),R_t(0))&R_t(1)\\\hline
0&(4,0,0,0)&0\\
1&(4,1,0,0)&0\\
2&(4,2,2,0)&0\\
3&(4,3,2,0)&0\\
4&(4,4,4,0)&0
\end{array}.
\]
Each listed sequence is increasing in the specified partial order. The
short chain has image $4\leq_S0\leq_S0$. Finally, for $j\in\{2,3\}$,
\[
R_t(1+j)=R_t(0)=0=R_t(1)+R_t(j).
\]
Thus every right translation preserves joins, and so the right distributive law holds.
This completes the proof.
\end{proof}

\subsection{A criterion for identities of $S$}

A \emph{word} over a set $X$ of variables is a nonempty finite sequence of
variables. For $\bw=x_1\cdots x_m$, write
\[
\ell(\bw)=m,\qquad c(\bw)=\{x_1,\ldots,x_m\},\qquad
h(\bw)=x_1,\qquad t(\bw)=x_m,
\]
and put $e(\bw)=\{x_1,\ldots,x_{m-1}\}$. Thus $e(\bw)=\varnothing$ when
$\ell(\bw)=1$. The set $c(\bw)$ is the \emph{content} of $\bw$. We distinguish a
one-letter word from a word of length at least two; the latter need not
contain a repeated variable.

By distributivity, every term can be represented as a nonempty finite sum
of words. We identify a term with the set of its summands, so that
\[
\bu=\bp_1+\cdots+\bp_k=\{\bp_1,\ldots,\bp_k\},
\qquad c(\bu)=\bigcup_{\bp\in\bu}c(\bp).
\]
Here $k\geq1$, and each $\bp_i$ is a nonempty word.
Throughout the proof, $\bq \preceq \bu$ is an abbreviation for the identity
$\bu+\bq \approx \bu$. It does not add a relation symbol to the signature.
This notation agrees with the natural order $a\leq_S b$ defined by $a+b=b$
in every additively idempotent semiring.

\begin{lem}\label{lem:S-word-value}
Let $s_1,\ldots,s_m\in S$, where $m\geq1$, and put $B_0=\{2,3\}$.
\begin{enumerate}
\item If every $s_i$ belongs to $B_0$, then $s_1\cdots s_m$ is $2$ if
at least one $s_i$ is $2$, and is $3$ otherwise.
\item Suppose that $j$ is the least index for which $s_j\notin B_0$.
If $s_j=0$ or $s_j=4$, the product is $0$ or $4$, respectively.
If $s_j=1$, the product is $1$ exactly when $j=m$ and
$s_1=\cdots=s_{m-1}=3$; in every other case it is $0$.
\end{enumerate}
\end{lem}
\begin{proof}
Multiplication on $\{2,3\}$ is the two-element semilattice operation with
$2$ as its greatest element. Upon reaching the first factor $0$ or $4$,
the product becomes that factor and remains unchanged, because both are
left zero elements. Before the first factor $1$, a nonempty prefix has
value $2$ or $3$. Now $2\cdot1=0$, $3\cdot1=1$, and $1$ multiplied by
any further factor is $0$. The same description holds if $1$ is the first
factor; in particular it also covers $m=1$.
\end{proof}

For a fixed occurrence of a letter $b$ in a word $\bq$, let $U$ denote the
set of letters occurring before that position. If this is the first
occurrence of $b$, then $b\notin U$. The
first letter of $\bp$ outside $U$ is $b$ if and only if
\[
\bp=\bv b\bw,\qquad c(\bv)\subseteq U,
\]
where the prefix $\bv$ and suffix $\bw$ are allowed to be empty.
Empty prefixes and suffixes describe positions in a word; they are not
terms or allowed substitution values in our constant-free signature.

\begin{lem}[Criterion for an identity $\bq \preceq \bu$]\label{lem:S-criterion}
Let $\bu$ be a term and $\bq$ a word. If $\ell(\bq)\geq2$ and $C=c(\bq)$, then
$S\models \bq \preceq \bu$ if and only if the following conditions hold.
\begin{enumerate}
\item[\textnormal{(A)}]
\[
C=\bigcup\{c(\bp):\bp \in\bu,\ c(\bp)\subseteq C\}.
\]
\item[\textnormal{(B)}]
At the first occurrence of each letter $b$ in $\bq$, with preceding content
$U$, there is a summand $\bp \in\bu$ whose first letter outside $U$
is $b$.
\item[\textnormal{(C)}]
If that first occurrence is not the last position of $\bq$, then either
some summand $\bp \in\bu$ satisfies $c(\bp)\subseteq U$, or there is a
summand $\bp=\bv b\bw \in\bu$ with $c(\bv)\subseteq U$ and $\bw$ nonempty.
\end{enumerate}
For a one-letter word $\bq=x$, the identity $x\preceq \bu$ holds in $S$ if and
only if $\bu$ contains a positive power of $x$.
\end{lem}
\begin{proof}
We first prove necessity for $\ell(\bq)\geq2$. If \textnormal{(A)} fails, choose
a letter $b\in C$ that belongs to none of the sets in the union. Assign
$2$ to $b$, $3$ to the remaining letters of $C$, and $4$ to all other
variables. The value of $\bq$ is $2$. Every summand whose content is
contained in $C$ has value $3$, and every other summand has value $4$.
Thus $\bu$ has value $3$ or $4$, contradicting $\bq \preceq \bu$.

If \textnormal{(B)} fails at the first occurrence of $b$, assign $3$ to
the letters of $U$, $0$ to $b$, and $4$ to all other variables. The value
of $\bq$ is $0$. A summand with content contained in $U$ has value $3$.
In every other summand, the first letter outside $U$ is different from $b$,
and hence the summand has value $4$. This contradicts the inequality.

If \textnormal{(C)} fails at a first occurrence of $b$ that is not the
last position of $\bq$, assign $3$ to $U$, $1$ to $b$, and $4$ elsewhere.
Then $\bq$ has value $0$. There is no summand with content contained in $U$.
Any summand whose first letter outside $U$ is $b$ ends at that occurrence
and has value $1$; all other summands have value $4$. Their sum is $1$
or $4$, and is therefore not greater than or equal to $0$.

For sufficiency, fix an arbitrary assignment $\varphi:X\to S$, extended to terms in the usual way. If $\varphi(\bq)=4$,
then $\varphi(\bq)\leq_S\varphi(\bu)$ because $4$ is the least element of the natural
order. If all letters of $\bq$ have values in $B_0$, condition
\textnormal{(A)} gives a summand with content contained in $C$ and value
$3$ when $\varphi(\bq)=3$. When $\varphi(\bq)=2$, choose a letter of $\bq$ with value $2$;
condition \textnormal{(A)} gives a summand containing that letter and no
letters outside $C$, so that summand also has value $2$. Hence the
inequality holds in both cases.

Suppose next that a letter of $\bq$ has value outside $B_0$, and let $b$
occur at the first position with such a value. This position is the first
occurrence of $b$, and all letters of its preceding content $U$ have
values in $B_0$. If $\varphi(b)=0$, the summand supplied by
\textnormal{(B)} has value $0$, whence $\varphi(\bu)=0$. If $\varphi(b)=4$, then
$\varphi(\bq)=4$, already considered. It remains to treat $\varphi(b)=1$.

If this position is not the last one, then $\varphi(\bq)=0$. A summand of the
form $\bv b\bw$ in \textnormal{(C)}, with $\bw$ nonempty, also has value $0$.
If no such summand exists, \textnormal{(C)} provides a summand with
content contained in $U$, whose value is $2$ or $3$. Condition
\textnormal{(B)} also provides a summand with value $0$ or $1$.
These two values have sum $0$, since $1+2=1+3=0$. Thus $\varphi(\bu)=0$.

Now suppose that $b$ occupies the last position of $\bq$. If all preceding
letters have value $3$, then $\varphi(\bq)=1$, and the summand in
\textnormal{(B)} has value $0$ or $1$, which proves the required
inequality. Otherwise a preceding letter $y$ has value $2$, and
$\varphi(\bq)=0$. If either alternative of \textnormal{(C)} holds,
the preceding argument still applies. Suppose neither holds. Every
summand with content contained in $C=U\cup\{b\}$ must contain $b$,
since otherwise its content would be contained in $U$. Moreover, its
first occurrence of $b$ must be its last position. By \textnormal{(A)},
one such summand contains $y$. Its prefix has value $2$, and its last
letter has value $1$, so its value is $2\cdot1=0$. Thus $\varphi(\bu)=0$ in
this final case as well.

For a one-letter word $x$, the multiplication table gives
\[
s\leq_S s^2,\qquad s^k=s^2\quad(s\in S,\ k\geq2).
\]
Consequently $x\preceq x^k$ holds for every $k\geq1$. Conversely, if
$\bu$ contains no positive power of $x$, assign $3$ to $x$ and
$4$ to every other variable. Every summand of $\bu$ contains a variable
different from $x$, so its value is $4$, whereas $x$ has value $3$.
Thus $x\preceq \bu$ fails. This completes the proof.
\end{proof}

\begin{cor}\label{cor:contained-cover}
If $S\models \bq \preceq \bu$, then, for each $x\in c(\bq)$, some summand $\bp$
of $\bu$ satisfies
\[
x\in c(\bp)\subseteq c(\bq).
\]
This holds for every nonempty word $\bq$, including a one-letter word.
\end{cor}
\begin{proof}
For $\ell(\bq)\geq2$, this is condition \textnormal{(A)}. For $\ell(\bq)=1$, use
the positive-power condition in Lemma~\ref{lem:S-criterion}.
\end{proof}

\subsection{An infinite family of identities}

For $n\geq2$, take pairwise distinct variables
\[
x_1,\ldots,x_n,a_1,\ldots,a_n,z
\]
and define
\[
\mathbf{x}_n=x_1\cdots x_n,\qquad
\bu_n=\sum_{i=1}^{n}x_i a_i+\mathbf{x}_n a_1z,\qquad
\bq_n=\mathbf{x}_n a_1a_2\cdots a_n.
\]

\begin{lem}\label{lem:critical-family}
The identity $\bq_n\preceq \bu_n$ holds in $S$ for every $n\geq2$.
\end{lem}
\begin{proof}
The summands $x_i a_i$ have contents whose union is $c(\bq_n)$, so
\textnormal{(A)} holds. At the first occurrence of $x_i$ in $\bq_n$, the
summand $x_i a_i$ starts with $x_i$ and has a nonempty suffix. At the
first occurrence of $a_1$, all $x_i$ have already appeared, so the
summand $\mathbf{x}_n a_1z$ satisfies both \textnormal{(B)} and the second
alternative of \textnormal{(C)}. At the first occurrence of $a_i$ for
$i\geq2$, the summand $x_i a_i$ satisfies \textnormal{(B)}. If $i<n$,
the preceding content contains $c(x_1a_1)$, so the first alternative
of \textnormal{(C)} holds. The first occurrence of $a_n$ is the last
position and imposes no requirement under \textnormal{(C)}.
Lemma~\ref{lem:S-criterion} now applies.
\end{proof}

For these fixed variables, consider the following property of a word $\bw$:
\begin{equation}\label{eq:nfb-property}
z\notin c(\bw)
\quad\text{and}\quad
\bigl|c(\bw)\cap\{a_1,\ldots,a_n\}\bigr|\geq2.
\end{equation}
Condition~\eqref{eq:nfb-property} refers to the letters of the word itself.

\begin{lem}\label{lem:critical-gap}
Suppose that $S\models \bw \preceq \bu_n$ and $\bw$ satisfies
\eqref{eq:nfb-property}. Then $\bw$ contains all the letters
$x_1,\ldots,x_n,a_1,\ldots,a_n$. More precisely, the first letter from
$\{a_1,\ldots,a_n\}$ to occur in $\bw$ is $a_1$, and every $x_i$ occurs
before it.
\end{lem}
\begin{proof}
Corollary~\ref{cor:contained-cover} shows that $c(\bw)$ is the union of
the contents of those summands of $\bu_n$ whose contents are contained
in $c(\bw)$. Since $z\notin c(\bw)$, the summand $\mathbf{x}_n a_1z$ cannot occur
in this union. Hence $c(\bw)$ is a union of some of the pairs
$\{x_i,a_i\}$. In particular, $\bw$ contains no other variables.

Let $a_j$ be the first letter from $\{a_1,\ldots,a_n\}$ occurring in
$\bw$, and let $U$ be its preceding content. The word contains at least
two distinct letters from this set, so this occurrence of $a_j$ is
not its last position. Since $U$ contains none of the $a_i$, no
summand of $\bu_n$ has content contained in $U$. Condition
\textnormal{(C)} therefore requires a summand whose first letter
outside $U$ is $a_j$ and which continues after that occurrence.
No summand $x_i a_i$ has this property, because its only occurrence
of an $a_i$ is its last position. The required summand must be
$\mathbf{x}_n a_1z$. Consequently $j=1$ and
\[
\{x_1,\ldots,x_n\}\subseteq U.
\]
All the $x_i$ therefore belong to $c(\bw)$. Since this content is a union
of whole pairs $\{x_i,a_i\}$, all the $a_i$ belong to it as well.
\end{proof}

\subsection{The nonfinite basis property}

The following elementary consequence of the criterion will be applied
to the variables in a proposed finite basis, before any substitution
is made.

\begin{lem}\label{lem:formal-cover}
Suppose that $\bs \approx \bt$ is an identity of $S$, and expand both sides
as sums of nonempty words. For every summand $\bv$ of $\bt$ and every
$t\in c(\bv)$, there is a summand $\bp_t$ of $\bs$ such that
\[
t\in c(\bp_t)\subseteq c(\bv).
\]
In particular, some summand of $\bs$ has content contained in $c(\bv)$.
\end{lem}
\begin{proof}
Every ai-semiring satisfies the inequality whose lower side is a summand
of a term and whose upper side is that term. Thus $S\models \bv \preceq \bt \approx \bs$, and
Corollary~\ref{cor:contained-cover} gives the assertion. All contents
in this argument are taken in the variables of $\bs \approx \bt$ itself.
\end{proof}

\begin{thm}\label{thm:S-nfb}
The semiring $S$ is nonfinitely based; equivalently,
$\Var(S)$ has no finite equational basis.
\end{thm}
\begin{proof}
Suppose for contradiction that a finite set $\Sigma$ is an equational basis
for $S$. We may include the finitely many additively idempotent
semiring axioms in $\Sigma$. Every identity of $\Sigma$ holds in $S$.
Expand both sides of each identity in $\Sigma$ by distributivity, and
let $L$ be the largest length of any word in these finitely many
expansions. Thus $L$ bounds the lengths of the words occurring in the chosen basis
identities. The lengths of the words in their substitution instances
may be arbitrarily large. Choose $n>L+1$.

By Lemma~\ref{lem:critical-family}, $\bu_n\approx \bu_n+\bq_n$ is an
identity of $S$, so it has a finite equational derivation from $\Sigma$.
Such a derivation may be written as
\[
\bu_n=\mathbf{d}_0\longrightarrow \mathbf{d}_1\longrightarrow\cdots
\longrightarrow \mathbf{d}_k=\bu_n+\bq_n,
\]
where each step replaces one occurrence of a substitution instance of
one side of an identity in $\Sigma$ by the corresponding instance of
the other side, inside a term. Both directions of each identity are
allowed. This is the usual description of the congruence generated
by substitution instances of the basis: closure under a basic
operation is obtained by placing one replacement inside that operation,
and simultaneous replacements can be made successively. The
variables may be replaced by arbitrary terms.

Expand every $\mathbf{d}_i$ as a sum of words. Since $S\models \mathbf{d}_i\approx \bu_n$,
each word $\bw$ in its expansion satisfies $S\models \bw \preceq \bu_n$.
We do not assume that the original summands of $\bu_n$ remain present
in these expansions. No summand of the initial term $\bu_n$
satisfies \eqref{eq:nfb-property}, whereas the summand $\bq_n$ of the
last term does. Choose the first step whose resulting expansion
has a word satisfying \eqref{eq:nfb-property}, and let $\bw$ be one
such word. Lemma~\ref{lem:critical-gap} implies that $\bw$ contains all
$n$ distinct letters $a_i$.

Write the chosen step as replacement of $\sigma(\bs)$ by $\sigma(\bt)$
inside a term with one specified position, where $\bs \approx \bt$ is a
basis identity oriented in the direction used and $\sigma$ is an
arbitrary term substitution. After expansion by distributivity, the
surrounding term consists of summands unaffected by the replacement
and expressions of the form
\[
\bp_0 \,[\ ]\,\br_0,
\]
where $\bp_0$ and $\br_0$ are words and either may be absent. Indeed, as
one moves from the specified position to the root of the term, an
addition introduces unaffected summands, and a multiplication attaches
each expanded word of the other factor on the appropriate side.
This proves the description by induction on the surrounding term.
The absence of a factor here does not permit an empty word as a
substitution value.

The word $\bw$ cannot arise from an unaffected summand,
since the preceding expansion has no word satisfying
\eqref{eq:nfb-property}. There is therefore a summand
$\bv=t_1\cdots t_m$ of the expansion of $\bt$, with $m\leq L$, such that
\[
\bw=\bp_0\mathbf{a}_1\cdots \mathbf{a}_m \br_0,
\]
where $\mathbf{a}_j$ is a word in the expansion of $\sigma(t_j)$. None of the
$\mathbf{a}_j$, $\bp_0$ and $\br_0$ contains $z$. For the same fixed
words $\bp_0,\br_0$, the preceding expansion contains all words in
$\bp_0 \sigma(\bp)\br_0$, for every summand $\bp$ of the expansion of $\bs$.
In forming each of these words, the expansion of a product permits
an independent choice of a summand at each occurrence of a substituted
variable.

We claim that each $\mathbf{a}_j$ contains at most one distinct letter from
$\{a_1,\ldots,a_n\}$. Otherwise put $t=t_j$ and apply
Lemma~\ref{lem:formal-cover} to obtain a summand $\bp_t$ of $\bs$ with
$t\in c(\bp_t)\subseteq c(\bv)$. At one occurrence of $t$ in $\bp_t$,
choose the word $\mathbf{a}_j$ from $\sigma(t)$. For every other occurrence
of every variable in $\bp_t$, choose a substitution summand that was
used at some occurrence of that variable in the fixed word $\bw$.
Such a choice exists because $c(\bp_t)\subseteq c(\bv)$, and none of
the chosen words contains $z$. The resulting word in the expansion
of $\bp_0 \sigma(\bp_t)\br_0$ contains the whole word $\mathbf{a}_j$, so it satisfies
\eqref{eq:nfb-property}. It already occurred before the chosen step,
a contradiction.

This argument does not require repeated occurrences of a variable
to choose the same summand. For example, if a variable is replaced
by $\mathbf{a} +\mathbf{b}$, then two occurrences expand to the four products
$\mathbf{a} \mathbf{a},\mathbf{a} \mathbf{b},\mathbf{b} \mathbf{a},\mathbf{b} \mathbf{b}$. Thus the independent choices used in the preceding
paragraph are available under an ordinary term substitution.

We also claim that $\bp_0$ and $\br_0$ together contain at most one distinct
letter from $\{a_1,\ldots,a_n\}$. If they contained two, choose by
Lemma~\ref{lem:formal-cover} any summand $\bp$ of $\bs$ with
$c(\bp)\subseteq c(\bv)$. For each variable occurrence in $\bp$, choose
a word used for that variable in $\bw$. The resulting word in
$\bp_0 \sigma(\bp)\br_0$ contains no $z$, while the two surrounding words
already supply two distinct $a_i$. This gives a word
satisfying \eqref{eq:nfb-property} before the chosen step, a contradiction.

It follows that the number of distinct letters $a_i$ in $\bw$ is at most
\[
\underbrace{1+\cdots+1}_{m\ \text{substitution factors}}+1
=m+1\leq L+1.
\]
On the other hand, Lemma~\ref{lem:critical-gap} implies that all $n>L+1$ distinct letters
$a_i$ occur in $\bw$. This is a contradiction, and the proof is complete.
\end{proof}

The preceding proof applies to arbitrary substitutions of terms and to
applications of identities in either direction within a term. In the subsequent
sections, we prove that every proper subvariety of $\Var(S)$ is
finitely based.
\section{Two maximal proper subvarieties and their finite bases}
\label{sec:foundations}

All algebras in this section have the constant-free signature $(+,\cdot)$.
We use the following subsemirings of $S$:
\[
 \begin{gathered}
 A_4=S|_{\{0,1,3,4\}},\qquad B_3=S|_{\{2,3,4\}},\qquad M_0=S|_{\{0,1,2,3\}},\\
 M_1=S|_{\{0,1,2,4\}},\qquad C_4=S|_{\{0,2,3,4\}},\qquad N_3=S|_{\{0,1,4\}}.
 \end{gathered}
\]
Closure under both operations follows directly from the tables of $S$.
Put $\mathcal V=\Var(S)$.
A term is a nonempty finite sum of nonempty words, with repeated summands removed.
As before, $\bq\preceq \bu$ abbreviates the identity $\bu+\bq\approx \bu$.
For a term $\bu$ we use
\[
 E(\bu)=\bigcup_{\bw\in \bu}e(\bw),\qquad
 T(\bu)=\{t(\bw):\bw\in \bu\}.
\]
In each derivation below, a word $\bw$ is added to a term $\bu$ only after
the identity $\bu+\bw\approx \bu$ has been established.

We introduce the following notation for the conditions used in the proof.
For a term $\bu$, a set $U$ of letters, and $b\notin U$, put
\[
\begin{aligned}
 \bu[U]&=\{\mathbf{v}_1\in \bu:c(\mathbf{v}_1)\subseteq U\},\\
 \bu(U,b)&=\{\mathbf{v}_1\in \bu:\mathbf{v}_1=\bp{}b\br\text{ with }c(\bp)\subseteq U\},\\
 \bu^{+}(U,b)&=\{\mathbf{v}_1\in \bu:\mathbf{v}_1=\bp{}b\br\text{ with }c(\bp)\subseteq U,\ \br\ne\varepsilon\}.
\end{aligned}
\]
In these decompositions only, $\bp$ or $\br$ may be empty, denoted by $\varepsilon$.
This notation does not introduce a constant into the semiring signature and does not permit substitution of an empty word for a variable.
Thus $\mathbf{v}_1\in \bu(U,b)$ means that the first letter of $\mathbf{v}_1$ outside $U$ is $b$;
$\mathbf{v}_1\in \bu^{+}(U,b)$ additionally means that this occurrence is not the final position.
The conditions in Lemma~\ref{lem:S-criterion} can consequently be written using
$\bu(U,b)\ne\varnothing$ and
\[
 \bu[U]\ne\varnothing\quad\text{or}\quad \bu^{+}(U,b)\ne\varnothing.
\]
We use these abbreviations throughout the proof.

We shall also use the following elementary observation.
For an element $a$ of a join semilattice, its \emph{principal down-set} is
$\mathord\downarrow a=\{b:b\le a\}$.
It is closed under addition and satisfies
\[
 b_1+\cdots+b_m\in\mathord\downarrow a
 \quad\Longleftrightarrow\quad
 b_i\in\mathord\downarrow a\text{ for every }i.
\]
For elements $b,c$ of a finite join semilattice, $b\le c$ if and only if every proper principal down-set containing $c$ also contains $b$.
Indeed, necessity is immediate; if $b\not\le c$, the principal down-set of $c$ is proper and separates them.
We shall apply this observation to evaluations of terms.
No closure under multiplication is required of a principal down-set.

\subsection{Normal forms and a finite set of auxiliary identities}

\begin{lem}[Word normal forms]\label{lem:wordnormal}
The semiring $S$ satisfies
\[
 \mathrm{W1}:x^2y\approx xy,\qquad
 \mathrm{W2}:xyx\approx xy^2.
\]
The identities W1 and W2 imply $xyxz\approx xyz$.
Let $\bt$ be the word consisting of the distinct letters of a nonempty word $\bw$, listed in order of first occurrence.
If the final letter of $\bw$ does not occur earlier, then $\bw\approx \bt$ follows from W1 and W2.
If the final letter occurs earlier, then $\bw\approx \bt^2$ follows from them.
\end{lem}
\begin{proof}
W1 and W2 follow from the evaluation formula in Lemma~\ref{lem:S-word-value}.
Right multiplication of W2 by $z$, followed by W1, gives
$xyxz\approx xy^2z\approx xyz$.
Consequently every repeated occurrence in the proper prefix of a word can be deleted; adjacent repetitions are deleted by W1.
After these deletions, the word is either $\bt$ or $\bt{}a$, where $a$ occurs in $\bt$.
Write $\bt=\mathbf{a}a\mathbf{b}$.
If $\mathbf{b}$ is nonempty, W2 gives $\bt{}a=\mathbf{a}a\mathbf{b}a\approx \mathbf{a}a\mathbf{b}^2$;
deleting repetitions from the proper prefix of $\bt^2=\mathbf{a}a\mathbf{b}\mathbf{a}a\mathbf{b}$ gives the same result.
If $\mathbf{b}$ is empty, deleting repetitions from the proper prefix of $\bt^2$ gives $\bt{}a$ directly.
An empty $\mathbf{a}$ means that a left factor is absent, not that a variable has been replaced by an empty word.
\end{proof}

In particular, $\bw\approx \bw^2$ holds whenever $t(\bw)\in e(\bw)$.
We distinguish words according to whether their final letter has occurred
earlier. All references to length mean the number of letter occurrences.

Let $\Gamma$ consist of the ai-semiring axioms, W1, W2, and the following fifteen identities:
\begin{align*}
 \mathrm L&:\ x\preceq x^2,&
 \mathrm{ST}&:\ xyz\preceq xz+yz,\\
 \mathrm{J2}&:\ x^2y^2\preceq x^2+y,&
 \mathrm H&:\ yxyz\preceq xyz+yt,\\
 \mathrm I&:\ xyzt\preceq xzyt+yu,&
 \mathrm{Blk}&:\ xayzt\preceq xazyt+y+a,\\
 \mathrm{TI}_{p,s}&:\ psar\preceq par+sa,&& (p,s\text{ independently present or absent}),\\
 \mathrm{Mark}&:\ paya\preceq pay+a,&&\\
 \mathrm{MN}_{p,v}&:\ pvt^2\preceq pvt+vtz,&& (p,v\text{ independently present or absent}).
\end{align*}
Here ``present or absent'' specifies four ordinary identities in each of the two indicated families.
It never allows empty-word substitutions.
For example, the four MN identities are
\[
 pvt^2\preceq pvt+vtz,\quad vt^2\preceq vt+vtz,\quad
 pt^2\preceq pt+tz,\quad t^2\preceq t+tz.
\]
The letters naming these identities are labels only.

\begin{lem}\label{lem:helpers}
Every identity in $\Gamma$ holds in $S$.
Moreover, $\Gamma$ implies
\[
 (x+y)^2\approx x^2+y^2,\qquad
 (xy)^2\approx x^2y^2,\qquad x^4\approx x^2.
\]
Hence $a\mapsto a^2$ is an idempotent endomorphism of every semiring satisfying $\Gamma$.
\end{lem}
\begin{proof}
We verify the inequalities by Lemma~\ref{lem:S-criterion}.
For L the required summand is a pure power.
For ST the summands $xz,yz$ cover all letters of the lower word; they belong respectively to the required sets $\bu^{+}(\varnothing,x)$ and $\bu^{+}(\{x\},y)$, and supply the condition for the final $z$ as well.
For J2 the summand $x^2$ supplies the nonfinal occurrence of $x$ and belongs to $\bu[\{x\}]$ before the first $y$.
For H the summand $xyz$ covers all letters of the lower word, while $yt$ supplies a nonfinal first $y$; the other first occurrences use $xyz$.
For I the summand $xzyt$ covers all letters of the lower word; $yu$ supplies a nonfinal occurrence of $y$ after $x$, and $xzyt$ supplies the other required decompositions.
For Blk the summand of length at least two contains all letters of the word on the left.
The set of letters preceding the first $y$ contains $a$, which is itself a summand;
thus the summand $y$ suffices even though its occurrence is final.
The other required decompositions use the summand $xazyt$.
For TI the summands $par,sa$ together contain all letters of the lower word and supply the required nonfinal occurrences.
If $s$ is absent, the lower word is itself the first summand; absence of $p$ does not alter the verification.
For Mark, $pay$ contains all letters of the word on the left; before the first $y$, the preceding content contains $a$.
For MN, $pvt$ contains all letters of the lower word, and $vtz$ supplies a nonfinal occurrence of $t$.
The same verification applies when $p$ or $v$ is absent, including the instance $t^2\preceq t+tz$.
These are identities with ordinary variables, so their validity includes arbitrary term substitutions.

By L, the normal form lemma, and J2,
\[
 xy\preceq (xy)^2\approx x^2y^2\preceq x^2+y\preceq x^2+y^2.
\]
Interchanging $x,y$ gives $yx\preceq x^2+y^2$.
Expansion of $(x+y)^2$ now gives the first asserted equality.
The second follows from the normal forms, since $(xy)^2$ and $x^2y^2$ have the same first-occurrence list and a repeated final letter.
The third follows from W1.
\end{proof}

\begin{lem}[Constructing auxiliary words]\label{lem:carrier-close}
Fix a term $\bu$ and an equational theory containing $\Gamma$.
All factors in the following statements are nonempty words, except for
the factors explicitly allowed to be absent. The following deductions
are valid in that theory.
\begin{enumerate}
\item If $\bw\preceq \bu$, $\mathbf{a}\preceq \bu$, and $c(\mathbf{a})\subseteq e(\bw)$, then $\bw^2\preceq \bu$.
\item If $\mathbf{u}_0t\preceq \bu$ and $\mathbf{v}_0t\mathbf{z}\preceq \bu$, where $t\notin c(\mathbf{u}_0)$, $\mathbf{z}$ is nonempty, and $c(\mathbf{v}_0)\subseteq c(\mathbf{u}_0)$, then $(\mathbf{u}_0t)^2\preceq \bu$.
The factors $\mathbf{u}_0,\mathbf{v}_0$ may be absent.
\item If $\bt^2\preceq \bu$ and $\mathbf{v}_1\preceq \bu$, then $\bt^2\mathbf{v}_1^2\preceq \bu$.
This word has content $c(\bt)\cup c(\mathbf{v}_1)$ and a repeated final letter.
\item If $\bs{}a\preceq \bu$ and $\bp{}a\br\preceq \bu$, where $\br$ is nonempty and $\bs,\bp$ may be absent, then one can deduce $\bw\preceq \bu$ for a word $\bw$ with
$c(\bw)=c(\bs{}a)\cup c(\bp{}a\br)$ and $t(\bw)\in e(\bw)$.
\end{enumerate}
\end{lem}
\begin{proof}
For (1), the assertion follows immediately if the final letter is repeated.
Otherwise write $\bw=\mathbf{u}_0t$.
The hypothesis forces $\mathbf{u}_0$ to be nonempty, and the normal forms give $\mathbf{u}_0t\approx \mathbf{u}_0\mathbf{a}t$.
Substitute $p=\mathbf{u}_0,a=\mathbf{a},y=t$ in Mark:
\[
 \mathbf{u}_0\mathbf{a}t\mathbf{a}\preceq \mathbf{u}_0\mathbf{a}t+\mathbf{a}\preceq \bu.
\]
The word on the left has the same first-occurrence list as $\bw$ and a repeated final letter, so it is equivalent to $\bw^2$.
For (2), the appropriate MN identity gives
$\mathbf{u}_0\mathbf{v}_0t^2\preceq \mathbf{u}_0\mathbf{v}_0t+\mathbf{v}_0t\mathbf{z}\preceq \bu$; the normal forms identify its left side with $(\mathbf{u}_0t)^2$.
Statement (3) is J2 with word substitutions.
For (4), TI gives $\bp{}\bs{}a\br\preceq \bp{}a\br+\bs{}a\preceq \bu$.
The proper prefix of $\bp{}\bs{}a\br$ contains every letter of $\bs{}a$, so (1) applies.
\end{proof}

\begin{lem}[Reordering a proper prefix]\label{lem:carrier-reorder}
Suppose $\bw\preceq \bu$ has been deduced, $e(\bw)=D$, and $t(\bw)=t$.
Let $b_1,\ldots,b_k$ be distinct letters of $D$.
The case $k=0$ is understood to require no change.
For $U_i=\{b_1,\ldots,b_{i-1}\}$ assume that
\[
 \bu(U_i,b_i)\ne\varnothing,\qquad
 \bu[U_i]\ne\varnothing\ \text{or}\ \bu^{+}(U_i,b_i)\ne\varnothing
 \quad(1\le i\le k).
\]
Using H, I, Blk, W1, and W2, one can deduce $\bw'\preceq \bu$, where $\bw'$ has the same final letter as $\bw$ and the same content in its proper prefix, and its proper prefix begins with $b_1\cdots b_k$.
This holds for $k=|D|$ and for $k<|D|$.
If $E(\bu)\cap T(\bu)=\varnothing$ and all $b_i\in E(\bu)$, the second displayed condition follows from the first.
\end{lem}
\begin{proof}
First remove repetitions from the proper prefix.
For $h=b_1$, the empty set cannot contain the letters of a nonempty summand.
The hypotheses therefore give an original summand $h\mathbf{z}$ with $\mathbf{z}$ nonempty.
If the current word is $\mathbf{a}h\mathbf{r}_0$ with $\mathbf{a}$ nonempty, then $\mathbf{r}_0$ is also nonempty, since it contains the final letter.
Apply H with $y=h,x=\mathbf{a},z=\mathbf{r}_0,t=\mathbf{z}$ to obtain
$h\mathbf{a}h\mathbf{r}_0\preceq \mathbf{a}h\mathbf{r}_0+h\mathbf{z}\preceq \bu$.
Deleting the repeated occurrence of $h$ gives a word beginning with $h$.

Suppose a nonempty prefix $\mathbf{u}_0$ is already in the required order and the next letter to be moved is $b$.
If $b$ does not immediately follow $\mathbf{u}_0$, the current word has the form $\mathbf{u}_0\mathbf{a}b\mathbf{r}_0$, with $\mathbf{a},\mathbf{r}_0$ nonempty.
If $\mathbf{v}_0b\bw_0$ is an original summand with $c(\mathbf{v}_0)\subseteq c(\mathbf{u}_0)$ and $\bw_0$ nonempty, use I with
\[
 x=\mathbf{u}_0,\quad y=\mathbf{v}_0b,\quad z=\mathbf{a},\quad t=\mathbf{r}_0,\quad u=\bw_0.
\]
Its first summand $\mathbf{u}_0\mathbf{a}\mathbf{v}_0b\mathbf{r}_0$ is equivalent to the current word, and its conclusion $\mathbf{u}_0\mathbf{v}_0b\mathbf{a}\mathbf{r}_0$ is equivalent to $\mathbf{u}_0b\mathbf{a}\mathbf{r}_0$.
If no such summand exists, there are original summands $\mathbf{a}_0,\mathbf{v}_0b$ with
$c(\mathbf{a}_0),c(\mathbf{v}_0)\subseteq c(\mathbf{u}_0)$.
Apply Blk with $x=\mathbf{u}_0,a=\mathbf{a}_0,y=\mathbf{v}_0b,z=\mathbf{a},t=\mathbf{r}_0$.
Its summand $\mathbf{u}_0\mathbf{a}_0\mathbf{a}\mathbf{v}_0b\mathbf{r}_0$ is equivalent to the current word; its two other summands are $\mathbf{v}_0b,\mathbf{a}_0$; its conclusion again is equivalent to $\mathbf{u}_0b\mathbf{a}\mathbf{r}_0$.
All deleted letters have already occurred in the proper prefix $\mathbf{u}_0$.
Repeating this argument proves the assertion. Every substitution used above replaces variables by nonempty words.

For the last assertion, a letter of $E(\bu)$ cannot be the last letter of any original summand, since $E(\bu)\cap T(\bu)=\varnothing$.
Thus every summand belonging to $\bu(U_i,b_i)$ also belongs to $\bu^{+}(U_i,b_i)$.
\end{proof}

\subsection{An infinite equational basis}

For each $n\ge2$, choose which of the formal letters $u_1,\ldots,u_n$ are present.
Put $\bw_i=u_i a_i$, interpreting $\bw_i=a_i$ when $u_i$ is absent.
Let $\mathbf{x}$ be the product, in index order, of the present $u_i$.
Consider the family of identities
\begin{equation}\label{eq:theta}
 \Theta_n:\qquad (\mathbf{x}a_1\cdots a_n)^2
 \preceq \bw_1+\cdots+\bw_n+\mathbf{x}a_1z.
\end{equation}
If all $u_i$ are absent, the last summand is $a_1z$ and the lower word is $(a_1\cdots a_n)^2$.
There are finitely many specified identities for each $n$.

\begin{lem}\label{lem:infinite-complete}
The set $\Gamma$, together with all identities in \eqref{eq:theta}, is an equational basis for $\mathcal V$.
\end{lem}
\begin{proof}
First, each identity in \eqref{eq:theta} holds in $S$ by Lemma~\ref{lem:S-criterion}.
The summands $\bw_i$ cover the letters of the lower word.
Each present $u_i$ has its required nonfinal occurrence in $\bw_i$;
$a_1$ has its required nonfinal occurrence in $\mathbf{x}a_1z$.
For every later $a_i$, the summand $\bw_i$ supplies the required first-letter decomposition, and the letters of $\bw_1$ have already occurred.

For completeness, suppose $S\models \bq\preceq \bu$, and first replace all words by their normal forms.
If $\bq$ is a single letter, the identity follows from a positive-power summand and L.
Assume $\ell(\bq)\ge2$ and put
\[
 F=\{\bw\in \bu:c(\bw)\subseteq c(\bq)\},\qquad
 E_F=\bigcup_{\bw\in F}e(\bw),\qquad T_F=\{t(\bw):\bw\in F\}.
\]
The criterion gives $\bigcup_{\bw\in F}c(\bw)=c(\bq)$.

Suppose first that $E_F\cap T_F\ne\varnothing$.
Lemma~\ref{lem:carrier-close}(4) gives a word with a repeated final letter.
By repeated applications of J2 to the remaining words of $F$, we obtain a square with content exactly $c(\bq)$.
We check that $\bq^2\preceq \bu$ satisfies all the conditions needed for reordering.
Only when the final letter $t$ of $\bq$ is unrepeated can a new condition arise: its first occurrence is no longer final in $\bq^2$.
Choose $a\in E_F\cap T_F$, with $\mathbf{s}_0\in F$ ending in $a$ and $\mathbf{v}_1\in F$ containing $a$ nonfinally.
If $a=t$, then $\mathbf{v}_1\in \bu^{+}(c(\bq)\setminus\{t\},t)$.
If $a\ne t$, either $\mathbf{s}_0$ contains $t$ nonfinally or $c(\mathbf{s}_0)\subseteq c(\bq)\setminus\{t\}$.
This proves the additional condition in both cases.
Lemma~\ref{lem:carrier-reorder} therefore gives $\bq^2\preceq \bu$, and L gives $\bq\preceq \bu$.

Next suppose $E_F\cap T_F=\varnothing$ and $T_F=\{t\}$.
Repeated applications of ST combine the nonempty proper prefixes of the words of $F$ into a word $\mathbf{u}_0t\preceq \bu$ with content $c(\bq)$.
If every word is the single letter $t$, simply take $\mathbf{u}_0t=t$.
If the final letter of $\bq$ is $t$ and is unrepeated, Lemma~\ref{lem:carrier-reorder} gives $\bq\preceq \bu$ directly.
Otherwise the first occurrence of $t$ in $\bq$ is nonfinal.
The letters preceding it cannot contain an entire original summand: such a summand would belong to $F$, yet would lack the common final letter $t$ of all members of $F$.
The criterion consequently gives an original summand $\mathbf{v}_0t\mathbf{z}$, with $c(\mathbf{v}_0)\subseteq E_F$ and $\mathbf{z}$ nonempty.
Lemma~\ref{lem:carrier-close}(2) gives $(\mathbf{u}_0t)^2\preceq \bu$.
If $\bq$ has an unrepeated final letter $b\ne t$, condition~\textnormal{(A)} ensures that some member of $F$ contains $b$ nonfinally.
This supplies the only possible new condition for $\bq^2$.
Lemma~\ref{lem:carrier-reorder} and L again yield $\bq\preceq \bu$.

It remains to consider $E_F\cap T_F=\varnothing$ and $|T_F|=n\ge2$.
Group the words of $F$ by their final letters and use ST within each group to obtain
$\mathbf{u}_i a_i\preceq \bu$; some $\mathbf{u}_i$ may be absent.
Renumber so that $a_1$ is the first of these final letters to occur in $\bq$.
Its first occurrence is nonfinal, and no original summand is contained in the letters preceding it, since those letters contain none of the $a_i$.
Thus there is an original summand $\mathbf{g}=\mathbf{v}_0a_1\mathbf{z}$ with $c(\mathbf{v}_0)\subseteq E_F$ and $\mathbf{z}$ nonempty.
TI combines $\mathbf{u}_1a_1$ and $\mathbf{g}$ into $\mathbf{v}_0\mathbf{u}_1a_1\mathbf{z}\preceq \bu$.
Since its proper prefix contains $c(\mathbf{u}_1a_1)$, Lemma~\ref{lem:carrier-close}(1) gives its square.
J2 incorporates the other groups and produces a square with content $c(\bq)\cup c(\mathbf{z})$.

Let $\mathbf{x}=\mathbf{u}_1\cdots \mathbf{u}_n$, omitting absent factors, and consider
\[
 (\mathbf{x}a_1\mathbf{z}a_2\cdots a_n)^2.
\]
Every first occurrence in $\mathbf{x}$ has a suitable nonfinal occurrence in the corresponding $\mathbf{u}_i a_i$.
The required occurrence of $a_1$ is supplied by $\mathbf{g}$.
For each letter first occurring in $\mathbf{z}$, the letters preceding its corresponding occurrence in $\mathbf{g}$ have already occurred; if this is the last position of $\mathbf{g}$, the preceding content also contains $c(\mathbf{u}_1a_1)$.
For each remaining $a_i$, the summand $\mathbf{u}_i a_i$ satisfies the required condition, since the preceding content again contains $c(\mathbf{u}_1a_1)$.
It follows from Lemma~\ref{lem:carrier-reorder} that the square just constructed can be replaced by this specified square.
Write the resulting word as $\mathbf{x}a_1\mathbf{y}$, with $\mathbf{y}$ nonempty.
Apply \eqref{eq:theta}, substituting the nonempty words $\mathbf{u}_i$ for the present $u_i$ and $\mathbf{y}$ for $z$.
All summands in that instance have already been proved to lie below $\bu$.
It follows that $(\mathbf{x}a_1\cdots a_n)^2\preceq \bu$.
Its content is exactly $c(\bq)$, so the letters contributed only by $\mathbf{z}$ have now been removed.

The conditions for reordering to $\bq^2$ still hold.
If $\bq$ has an unrepeated final letter $t$, choose a different group final letter $a_i\ne t$.
Its group word either does not contain $t$, or contains it nonfinally; these are precisely the two alternatives needed before the first $t$ in $\bq^2$.
A final application of Lemma~\ref{lem:carrier-reorder} and L proves $\bq\preceq \bu$.

The preceding cases exhaust all possibilities for arbitrary numbers of variables and summands, and for arbitrary word lengths.
Every term identity is equivalent to requiring each summand on either side to lie below the other side.
This proves completeness.
\end{proof}

\subsection{A finite basis for $\mathcal C$}

Put
\[
 \chi_A:\ yx^2\preceq x+z+yx+zy,\qquad
 \mathcal C=\Var(M_0,M_1,C_4).
\]
For a term $\bu$, write
\[
 E(\bu)=\bigcup_{\bw\in \bu}e(\bw),\qquad T(\bu)=\{t(\bw):\bw\in \bu\}.
\]
We shall consider separately the cases
$E(\bu)\cap T(\bu)=\varnothing$ and $E(\bu)\cap T(\bu)\ne\varnothing$.

\begin{lem}[Identity criteria for the three generators]\label{lem:factor-criteria}
Let $F=\{\bw\in \bu:c(\bw)\subseteq c(\bq)\}$.
\begin{enumerate}
\item $M_0\models \bq\preceq \bu$ if and only if $c(\bq)\subseteq c(\bu)$ and
\[
 E(\bu)\cap T(\bu)\ne\varnothing\quad\text{or}\quad
 \bigl(e(\bq)\subseteq E(\bu)\text{ and }t(\bq)\in T(\bu)\bigr).
\]
\item If $\ell(\bq)\ge2$, then $M_1\models \bq\preceq \bu$ if and only if $F\ne\varnothing$ and, for every first occurrence of a letter $b$ of $\bq$, with $U$ the set of preceding letters,
\[
 \bu(U,b)\ne\varnothing,\qquad
 \bu[U]\ne\varnothing\text{ or }\bu(U,b)\text{ contains a word of length at least two}.
\]
This includes an unrepeated final letter of $\bq$.
The last condition concerns the length of the whole summand, not the position of $b$ in it.
\item $C_4\models \bq\preceq \bu$ if and only if $\bigcup_{\bw\in F}c(\bw)=c(\bq)$ and $\bu(U,b)\ne\varnothing$ at every first occurrence in $\bq$.
\end{enumerate}
For a lower word of length one $x$, the criterion in both $M_1$ and $C_4$ is the existence in $\bu$ of a positive power of $x$.
Thus, for $\ell(\bq)\ge2$, the criterion in $\mathcal C$ consists of the following four conditions:
\begin{description}
\item[(K1)] $\bigcup_{\bw\in F}c(\bw)=c(\bq)$;
\item[(K2)] $\bu(U,b)\ne\varnothing$ at each first occurrence;
\item[(K3)] $\bu[U]\ne\varnothing$ or some word of $\bu(U,b)$ has length at least two, at each first occurrence;
\item[(K4)] $E(\bu)\cap T(\bu)\ne\varnothing$, or $e(\bq)\subseteq E(\bu)$ and $t(\bq)\in T(\bu)$.
\end{description}
\end{lem}
\begin{proof}
In the additive order of $M_0$, the proper principal down-sets are
$\{1\}$, $\{2,3\}$ and $\{3\}$.
All summands of $\bu$ evaluate to $1$ exactly when every letter of $E(\bu)$ takes value $3$ and every letter of $T(\bu)$ takes value $1$.
Such assignments exist precisely when $E(\bu)\cap T(\bu)=\varnothing$.
All summands lie in $\{2,3\}$, respectively $\{3\}$, exactly when all letters take values in that same set.
The observation concerning principal down-sets gives (1), including the case in which $\bq$ is a single letter.

In $M_1$, a word of length at least two can be evaluated by omitting its initial factors equal to $2$.
The first $0$ or $1$ makes the value permanently $0$, and the first $4$ makes it permanently $4$; a word consisting entirely of $2$'s has value $2$.
The proper principal down-sets are $\{1,4\},\{2,4\},\{4\}$.
If $F$ is empty, assign $2$ to the letters of the lower word and $4$ to the others.
If $\bu(U,b)$ is empty at a first occurrence, assign $2$ to $U$, $0$ to $b$, and $4$ to all other letters.
If $\bu[U]$ is empty and all summands in $\bu(U,b)$ are single-letter words, use the same assignment with $b=1$.
The term then takes value $1$ or $4$, whereas the lower word takes value $0$.
These assignments prove necessity in (2).
Conversely, a value $2$ of the lower word is handled by a summand in $F$.
A first $0$ is handled by a summand in $\bu(U,b)$.
A first $1$ either occurs in a summand of length at least two, giving value $0$, or combines with the value $2$ of a summand in $\bu[U]$ to give $1+2=0$.
If $\bq$ has value $4$, the inequality holds because $4$ is the least element in the additive order.
This proves sufficiency.

In $C_4$, the first $0$ or $4$ determines a word's value; if all letters are $2$ or $3$, the value is $2$ exactly when some letter is $2$.
If $\bigcup_{\bw\in F}c(\bw)\ne c(\bq)$, assign $2$ to a letter of $\bq$ outside this union, $3$ to the other letters of $\bq$, and $4$ elsewhere.
If a condition $\bu(U,b)\ne\varnothing$ fails, assign $3$ to $U$, $0$ to $b$, and $4$ elsewhere.
Conversely, the conditions on contents and first occurrences give the inequality for each possible value $0,2,3,4$ of $\bq$.
This gives (3).
For a lower word of length one, necessity follows by assigning that letter $2$ and all others $4$; sufficiency follows from $x\preceq x^2\approx x^3$.
Finally, an identity holds in the variety generated by several algebras exactly when it holds in every generator.
Combining (1)--(3) proves the final assertion.
\end{proof}

\begin{lem}\label{lem:C-closing}
In any equational theory containing $\Gamma\cup\{\chi_A\}$, if
$\mathbf{z}\preceq \bu$, $t(\mathbf{z})\in e(\mathbf{z})$, and $\bw\preceq \bu$ with $\ell(\bw)\ge2$, then $\bw^2\preceq \bu$.
No inclusion between $c(\mathbf{z})$ and $c(\bw)$ is required.
\end{lem}
\begin{proof}
Using J2, first replace $\mathbf{z}$ by a word still lying below $\bu$ whose final letter is repeated and whose content contains $c(\bw)$.
Write this enlarged word again as $\mathbf{z}$.
Let $h=h(\bw)$ and substitute $x=\bw,y=h,z=\mathbf{z}$ in $\chi_A$.
The word normal forms give $h\bw\approx \bw$ and $h\bw^2\approx \bw^2$.
They also give $\mathbf{z}h\approx \mathbf{z}$, since $\mathbf{z}$ has a repeated final letter and contains $h$.
The resulting instance is therefore $\bw^2\preceq \bw+\mathbf{z}\preceq \bu$.
\end{proof}

\begin{thm}\label{thm:C-fb}
We have
\[
 \mathcal V\cap\operatorname{Mod}(\chi_A)=\mathcal C,
\]
and $\Gamma\cup\{\chi_A\}$ is a finite equational basis for $\mathcal C$.
\end{thm}
\begin{proof}
The identity $\chi_A$ holds in $M_0,M_1,C_4$ by Lemma~\ref{lem:factor-criteria}.
For its term $\bu=x+z+yx+zy$, the intersection $E(\bu)\cap T(\bu)$ is nonempty.
The summands $x,yx$ cover the letters of the lower word; $yx$ supplies a summand of length two at the first $y$ and at the first $x$.
Thus $\mathcal C\subseteq\mathcal V\cap\operatorname{Mod}(\chi_A)$.

To prove completeness, suppose $\bq\preceq \bu$ satisfies (K1)--(K4).
The case in which $\bq$ is a single letter has already been proved.
If $E(\bu)\cap T(\bu)=\varnothing$, (K4) forces the final letter $t$ of $\bq$ to be unrepeated.
Every member of $F$ ends in $t$: its final letter lies in $T(\bu)$, while all other letters of $\bq$ lie in the disjoint set $E(\bu)$.
By (K1), ST combines these summands into a word with the same content and final letter as $\bq$.
Every occurrence needed for reordering is nonfinal in its original summand, since its letter belongs to $E(\bu)$.
Lemma~\ref{lem:carrier-reorder} gives $\bq\preceq \bu$.

Suppose $E(\bu)\cap T(\bu)\ne\varnothing$.
Lemma~\ref{lem:carrier-close}(4) and J2 give a word $\mathbf{z}\preceq \bu$ with a repeated final letter and content containing $c(\bu)$.
Lemma~\ref{lem:C-closing} then gives $\bw^2\preceq \bu$ for every original summand $\bw$ of length at least two.
Consequently, every summand of length at least two used in (K3) may be replaced by its square, in which the relevant first occurrence is nonfinal.
We also need a square lying below $\bu$ with content exactly $c(\bq)$.
If $F$ contains a word of length at least two, start with its square.
Otherwise all members of $F$ are single letters, including $h=h(\bq)$ by (K1).
Condition (K3) at the empty prefix supplies an original summand $h\mathbf{v}_0$ with $\mathbf{v}_0$ nonempty.
The instance $h^2\preceq h+h\mathbf{v}_0$ of MN supplies the required initial square.
J2 incorporates the remaining members of $F$ without adding letters outside $c(\bq)$.
All first occurrences, including a formerly unrepeated final letter, now satisfy the hypotheses of Lemma~\ref{lem:carrier-reorder}.
It gives $\bq^2\preceq \bu$, and L gives $\bq\preceq \bu$.

Thus every identity of $\mathcal C$ follows from $\Gamma\cup\{\chi_A\}$.
This finite set holds in $\mathcal V\cap\operatorname{Mod}(\chi_A)$, and its completeness identifies its model class with $\mathcal C$.
Both assertions follow.
\end{proof}

\subsection{A finite basis for $\mathcal E$}

Put
\[
 \varepsilon:\ xy\preceq x+xyz,\qquad
 \mathcal E=\mathcal V\cap\operatorname{Mod}(\varepsilon).
\]
\begin{thm}\label{thm:E-fb}
The set $\Gamma\cup\{\varepsilon\}$ is a finite equational basis for $\mathcal E$.
\end{thm}
\begin{proof}
By Lemma~\ref{lem:infinite-complete}, it suffices to derive every instance of \eqref{eq:theta}.
Write its summands as $\bw_i=u_i a_i$ and $\mathbf{g}=\mathbf{x}a_1r$, where each $u_i$ is a formal letter or is absent.
If all $u_i$ are absent, MN gives $a_1^2\preceq a_1+a_1r$.
J2 incorporates the other $a_i$, and the first-occurrence conditions used to verify \eqref{eq:theta} allow reordering to the square on the lower side.

Otherwise choose a present $u_i$.
Apply Lemma~\ref{lem:carrier-close}(4) to $\bw_1,\mathbf{g}$ and then J2 to the other groups.
This gives a word $\mathbf{z}$ below the term, with a repeated final letter and containing all letters of \eqref{eq:theta}, including $r$.
Consider
\[
 \bw_0=u_i\mathbf{x}a_1\mathbf{z}a_i.
\]
It has the same content as $\mathbf{z}$ and a repeated final letter.
Its initial $u_i$, and each new letter in $\mathbf{x}$, have suitable nonfinal occurrences in their corresponding $\bw_j$.
The occurrence of $a_1$ is supplied by $\mathbf{g}$.
For every letter first occurring in the displayed factor $\mathbf{z}$, the required summand is $\mathbf{z}$ itself: the preceding letters in $\mathbf{z}$ have already occurred, and its first occurrence in $\mathbf{z}$ cannot be final because the final letter of $\mathbf{z}$ is repeated.
The last $a_i$ is repeated and imposes no new first-occurrence condition.
Lemma~\ref{lem:carrier-reorder} therefore gives $\bw_0$ below the term.

Substitute $x=u_i,y=\mathbf{x}a_1,z=\mathbf{z}$ in $\varepsilon$ and right-multiply by $a_i$:
\[
 u_i\mathbf{x}a_1a_i\preceq u_i a_i+u_i\mathbf{x}a_1\mathbf{z}a_i=\bw_i+\bw_0.
\]
All substituted factors are nonempty.
If $i=1$, the final letter on the left is already repeated.
If $i\ne1$, its proper prefix contains all letters of $\bw_1=u_1a_1$, even when $u_1$ is absent; Lemma~\ref{lem:carrier-close}(1) gives its square.
This square contains only letters of the lower word and includes all present $u_j$.
Applying J2 with the other $\bw_j$ produces a square whose content is exactly the content of the lower word; $r$ has been eliminated.
The first-occurrence conditions for \eqref{eq:theta} and Lemma~\ref{lem:carrier-reorder} give $(\mathbf{x}a_1\cdots a_n)^2$ below the term.
The argument applies to every $n$ and to every choice of the present letters $u_i$.
Hence the finite set derives the full basis of $\mathcal V$ together with $\varepsilon$, as required.
\end{proof}

\subsection{Subvarieties not containing $A_4$ or $B_3$}

\begin{lem}\label{lem:A4-characteristic}
For every $\mathcal W\subseteq\mathcal V$,
\[
 A_4\notin\mathcal W\quad\Longleftrightarrow\quad
 \mathcal W\models\chi_A.
\]
\end{lem}
\begin{proof}
Put $\bv=x+z+yx+zy$ and evaluate $(x,y,z)=(1,3,4)$ in $A_4$.
Use the principal down-set $I=\{1,4\}$.
The words $x^2,x,y,z$ represent $0,1,3,4$, respectively.
Membership in $I$ separates every pair having one element inside $I$ and the other outside.
The two remaining pairs $0,3$ and $1,4$ are separated after right multiplication by $1$: their images are $0,1$ and $0,4$.

Under this evaluation, a word belongs to $I$ exactly when its first letter different from $y$ is $z$, or when it is $y^kx$ for some $k\ge0$.
It lies outside $I$ exactly when it is a positive power of $y$, or is $y^kx\mathbf{v}_0$ with $\mathbf{v}_0$ nonempty.
This description applies to arbitrary word lengths.
By W1 and W2, every word over three letters has one of the thirty forms $\bt,\bt^2$, where $\bt$ is a nonempty list of distinct letters.
The forms evaluating into $I$ are exactly
\[
 x,\ yx,\ yz,\ (yz)^2,\ yzx,\ (yzx)^2,
 \qquad \bt,\bt^2\ (\bt=z,zx,zy,zxy,zyx).
\]
Each listed word $\mathbf{v}_1$ satisfies $\mathbf{v}_1\preceq \bv$ in $S$.
For $x,yx$ this is immediate.
For a word starting with $z$, the summands $z,zy,x$ satisfy condition~\textnormal{(A)}, $zy$ supplies a nonfinal initial $z$, and the summand $z$ is contained in the letters preceding every later first occurrence.
For a word starting with $yz$, use $yx$ for the initial $y$ and $zy$ for the first $z$ and for covering $y$.
If an $x$ follows, the summand $x$ supplies the required decomposition, while the preceding content already contains $z$.
These verify all conditions of Lemma~\ref{lem:S-criterion}.

The other fourteen forms are
\[
\begin{gathered}
 x^2,y,y^2,xy,(xy)^2,xz,(xz)^2,(yx)^2,\\
 xyz,(xyz)^2,xzy,(xzy)^2,yxz,(yxz)^2.
\end{gathered}
\]
For each such form $\br$, the identity $yx^2\preceq \bv+\br$ holds in $S$.
Indeed, $x,yx$ already cover the letters of the lower word, and $yx$ supplies the required nonfinal initial $y$.
The only missing condition concerns the first $x$, whose preceding content is $\{y\}$.
If $\br=y$ or $y^2$, its entire content is contained in $\{y\}$.
Every other listed $\br$ supplies a nonfinal $x$ after a possibly empty prefix of $y$'s.
Thus the assertion follows from Lemma~\ref{lem:S-criterion} for words of arbitrary length.

Suppose an identity $\rho$ fails in $A_4$.
After choosing the separating right multiplication just described, represent its assigned elements by $x^2,x,y,z$ and interchange the two sides if necessary.
This gives a consequence $\mathbf{r}_0\approx \bt$ of $\rho$ whose left value lies in $I$ and whose right value does not.
Since $I$ is a principal down-set, every summand of $\mathbf{r}_0$ lies in $I$ and at least one summand $\br$ of $\bt$ lies outside it.
The preceding argument gives $\bv+\mathbf{r}_0\approx \bv$ and $yx^2\preceq \bv+\br$ in $S$.
Using $\mathbf{r}_0\approx \bt$, one obtains $\bv+\bt\approx \bv$, hence $\br\preceq \bv$ and then $yx^2\preceq \bv$.
If $A_4\notin\mathcal W$, some identity of $\mathcal W$ fails in $A_4$, so $\mathcal W\models\chi_A$.
Conversely, at the specified evaluation $\bv=1$ but $yx^2=0$; thus $A_4$ fails $\chi_A$.
\end{proof}

\begin{lem}\label{lem:B3-characteristic}
For every $\mathcal W\subseteq\mathcal V$,
\[
 B_3\notin\mathcal W\quad\Longleftrightarrow\quad
 \mathcal W\models\varepsilon.
\]
\end{lem}
\begin{proof}
Consider the following elements of the three-generated free algebra:
\[
 b_2=y+yx+yxz,\qquad b_3=y+yxz,\qquad b_4=yxzx.
\]
Under every evaluation in $S$, the triple $(b_2,b_3,b_4)$ is one of
\[
 (0,0,0),(2,2,2),(3,3,3),(4,4,4),
 (2,2,4),(3,3,4),(2,3,4).
\]
To verify this assertion, apply Lemma~\ref{lem:S-word-value} as follows.
If $y=0,1,4$, the triple is respectively constant $0,0,4$.
If $y=2$, then $x=0,1$ gives the constant $0$ triple and $x=4$ gives $(2,2,4)$;
for $x=2,3$, the three cases $z\in\{0,1\}$, $z\in\{2,3\}$, and $z=4$ give respectively constant $0$, constant $2$, and $(2,2,4)$.
If $y=3$, then $x=0,1$ gives constant $0$ and $x=4$ gives $(3,3,4)$.
For $x=2$, the same three cases for $z$ give constant $0$, constant $2$, and $(2,3,4)$.
For $x=3$, the values $z=0,1,2,3,4$ give respectively constant $0$, constant $0$, constant $2$, constant $3$, and $(3,3,4)$.

Each listed triple defines, by $i\mapsto b_i$, a homomorphism from $B_3$ into $S$.
The constant triples take multiplicatively idempotent values;
$(2,2,4)$ and $(3,3,4)$ are obtained by identifying $2,3$;
$(2,3,4)$ is the identity map.
Consequently the identities expressing preservation of the two operations hold in $S$, and hence in every $\mathcal W\subseteq\mathcal V$.
The three free elements therefore define a homomorphism
$B_3\to F_{\mathcal W}(3)$.

Every nontrivial congruence of $B_3$ identifies $2,3$.
If $2\sim4$, addition of $3$ gives $2\sim3$.
If $3\sim4$, left multiplication by $2$ first gives $2\sim4$, and then the preceding argument applies.
The partition $\{2,3\}\mid\{4\}$ is itself a congruence.
Thus, if $b_2\ne b_3$ in the free algebra, the displayed homomorphism is injective and $B_3\in\mathcal W$.
It follows that $B_3\notin\mathcal W$ forces $b_2=b_3$, namely
$yx\preceq y+yxz$.
Interchanging $x,y$ gives $\varepsilon$.
Conversely, in $B_3$ the assignment $(x,y,z)=(3,2,4)$ gives $xy=2$ and $x+xyz=3$;
hence $B_3$ fails $\varepsilon$.
\end{proof}

\begin{thm}\label{thm:two-coatoms}
The variety $\mathcal V$ has exactly two maximal proper subvarieties, namely $\mathcal C$ and $\mathcal E$.
Every proper subvariety is contained in at least one of them.
\end{thm}
\begin{proof}
First,
\[
 \mathcal V=\Var(A_4)\vee\Var(B_3).
\]
An explicit proof is supplied by the subsemiring of $A_4\times B_3$ with elements
\[
 (0,2),(0,3),(0,4),(1,2),(1,3),(3,2),(3,3),(4,4).
\]
In this order, map them to $0,0,0,0,1,2,3,4$.
This is a surjective homomorphism to $S$.
Both operations are computed coordinatewise: the three pairs with first coordinate $0$, together with $(1,2)$, form the fiber of $0$; all other fibers are singletons.
The operation tables verify closure and show that the operations on these five fibers are exactly those of $S$.
The reverse inclusion holds because $A_4,B_3$ are subsemirings of $S$.

We also give a proof using identities, since the resulting criterion for $A_4$ will be used later.
In $B_3$, the identity $\bq\preceq \bu$ holds if and only if $c(\bq)$ is the union of the contents of those summands of $\bu$ whose contents are contained in $c(\bq)$.
Indeed, in $B_3$ a word has value $4$ if some letter has value $4$; otherwise it has value $2$ if some letter has value $2$, and value $3$ in the remaining case.
A letter of the lower word not covered by these summands can therefore be assigned $2$, the other letters of the lower word $3$, and all external letters $4$, disproving the inequality.
Conversely, a covering summand containing a letter assigned $2$ gives value $2$; if the lower word has value $3$, any contained summand has value $3$. A value $4$ of the lower word is least in the additive order and requires no further condition.
For $A_4$ and a lower word of length at least two, the criterion is: at least one such summand exists; $\bu(U,b)\ne\varnothing$ at every first occurrence; and
$\bu[U]\ne\varnothing$ or $\bu^{+}(U,b)\ne\varnothing$ whenever that occurrence is nonfinal in $\bq$.
To see this, evaluate at the first letter whose value differs from $3$.
A first $4$ gives the least additive value.
A first $0$ requires a summand in $\bu(U,b)$.
A first $1$ at the last position again requires such a summand; at a nonfinal position it requires either a summand in $\bu^{+}(U,b)$, or the sum of a value $3$ from $\bu[U]$ and a value $1$.
If every letter of the lower word takes value $3$, the required contained summand gives value $3$.
Necessity follows by assigning $3$ to the preceding letters and $0$, $1$, or $4$ as in these cases.
For a lower word of length one in $A_4$, the criterion is the existence of a pure-power summand.
For necessity assign $3$ to that letter and $4$ elsewhere; sufficiency follows from L and the word normal forms.
Intersecting the $A_4$ and $B_3$ criteria gives exactly Lemma~\ref{lem:S-criterion}, proving the equality of varieties again.

A proper subvariety cannot contain both $A_4$ and $B_3$.
The preceding two lemmas therefore place every proper subvariety inside $\mathcal C$ or $\mathcal E$.
Moreover, $B_3\in\mathcal C\setminus\mathcal E$, since $B_3$ is a subsemiring of $C_4$ and fails $\varepsilon$.
Also $A_4\in\mathcal E\setminus\mathcal C$: the displayed $A_4$ criterion verifies $xy\preceq x+xyz$, whereas $A_4$ fails $\chi_A$.
Thus both subvarieties are proper and neither contains the other.
Any proper subvariety strictly above one of them would still have to be contained in one of these two, which is impossible.
They are therefore precisely the two maximal proper subvarieties.
\end{proof}

\subsection{Three further identities defining subvarieties of $\mathcal C$}

\begin{lem}\label{lem:C-omit-M0}
Put $\mathcal D=\Var(M_1,C_4)$ and $\chi_0:xy\approx(xy)^2$.
For every $\mathcal W\subseteq\mathcal C$,
\[
 M_0\notin\mathcal W\quad\Longleftrightarrow\quad\mathcal W\models\chi_0,
 \qquad
 \mathcal C\cap\operatorname{Mod}(\chi_0)=\mathcal D.
\]
In particular, $\Gamma\cup\{\chi_A,\chi_0\}$ is a finite basis for $\mathcal D$.
\end{lem}
\begin{proof}
The subsemiring $S|_{\{0,1,3\}}$ has the same identity criterion as $M_0$, so it generates the same variety.
Its distinct elements can be separated into value $1$ and a value different from $1$ by either the identity map or right multiplication by $1$.
If an identity fails there, choose such a right factor and represent values $3,1,0$ by $x,y,y^2$, respectively.
One obtains a two-variable consequence $\mathbf{r}_0\approx \bt$ with $\mathbf{r}_0(3,1)=1$ and $\bt(3,1)\ne1$.
Every summand of $\mathbf{r}_0$ is $x^ky$ for some $k\ge0$, so $x\mathbf{r}_0\approx xy$.
Some summand of $\bt$ is not of this form.
After prefixing it by $x$, it is either a power of $x$ of length at least two or has both $x,y$ in its proper prefix.
In particular, the resulting term $x\bt$ satisfies
$E(x\bt)\cap T(x\bt)\ne\varnothing$.
The $M_0$ criterion gives $x\bt\approx(x\bt)^2$ in $M_0$.
In $M_1$, all its summands have length at least two and take values in the subsemiring $\{0,2,4\}$, in which multiplication is idempotent.
In $C_4$, every element is already multiplicatively idempotent.
Hence $x\bt\approx(x\bt)^2$ holds in $\mathcal C$.
Using $x\mathbf{r}_0\approx x\bt$ gives $xy\approx(xy)^2$.
Conversely, this identity fails in $M_0$ at $(x,y)=(3,1)$.

Both $M_1$ and $C_4$ satisfy $\chi_0$.
Consider an identity $\bq\preceq \bu$ of $\mathcal D$ with $\ell(\bq)\ge2$.
The $M_1$ criterion at its first position gives an original summand $\bw$ of length at least two.
Using $\chi_0$, replace $\bw$ by $\bw^2$.
This preserves the values in $M_1,C_4$, and it makes $E(\bu)\cap T(\bu)$ nonempty for the modified term.
The remaining $M_0$ condition, $c(\bq)\subseteq c(\bu)$, follows from the corresponding condition for $C_4$.
Thus the modified identity holds in $\mathcal C$ and follows from its finite basis.
For a lower word of length one, a pure-power summand already suffices.
This proves the equality of varieties and the asserted finite basis.
\end{proof}

\begin{lem}\label{lem:C-omit-M1}
Put
\[
 \chi_1:\ yxz\preceq x+z^2+zy^2+yzx^2,\qquad
 \mathcal J^*=\Var(M_0,C_4,N_3).
\]
For every $\mathcal W\subseteq\mathcal C$,
\[
 M_1\notin\mathcal W\quad\Longleftrightarrow\quad\mathcal W\models\chi_1,
 \qquad
 \mathcal C\cap\operatorname{Mod}(\chi_1)=\mathcal J^*.
\]
Consequently $\Gamma\cup\{\chi_A,\chi_1\}$ is a finite basis for $\mathcal J^*$.
\end{lem}
\begin{proof}
Put $\bv=x+z^2+zy^2+yzx^2$, evaluate $(x,y,z)=(1,2,4)$ in $M_1$, and use the principal down-set $I=\{1,4\}$.
The words $x^2,x,y,z$ represent its four elements.
Membership in $I$, together with membership after right multiplication by $2$ or $4$, separates all distinct elements.
In particular, right multiplication by $2$ separates $1,4$, and right multiplication by $4$ separates $0,2$.
A word evaluates into $I$ exactly when it is the single letter $x$ or its first letter different from $y$ is $z$.
All other words are positive powers of $y$ or have length at least two with $x$ as the first letter different from $y$.
Of the thirty normal forms on these variables, the fifteen evaluating into $I$ are
\[
 x,\ yz,(yz)^2,yzx,(yzx)^2,
 \qquad \bt,\bt^2\ (\bt=z,zx,zy,zxy,zyx).
\]
Each of these words $\mathbf{v}_1$ satisfies $\mathbf{v}_1\preceq \bv$ in $\mathcal C$.
An initial $z$ uses $z^2$; an initial $y$ uses $yzx^2$.
After the first $z$, the preceding content contains the content of the summand $z^2$.
Condition (K1) is satisfied by $x,z^2,zy^2$, and (K4) is automatic because $E(\bv)\cap T(\bv)\ne\varnothing$.
The other fifteen forms are
\[
\begin{gathered}
 x^2,y,y^2,xy,(xy)^2,xz,(xz)^2,yx,(yx)^2,\\
 xyz,(xyz)^2,xzy,(xzy)^2,yxz,(yxz)^2.
\end{gathered}
\]
For every such $\br$, one has $yxz\preceq \bv+\br$ in $\mathcal C$.
Conditions (K1),(K2) already hold for $\bv$.
At the initial $y$, there is a suitable summand of length at least two.
At the first $x$, condition (K3) follows from $\br=y,y^2$, whose content is contained in $\{y\}$, or from a summand $\br$ of length at least two whose first letter outside $\{y\}$ is $x$.
Before the first $z$, the preceding content contains $x$.
Condition (K4) remains automatic.
The separating evaluations and term expansions used in Lemma~\ref{lem:A4-characteristic} now show that any identity failing in $M_1$ implies $\chi_1$ over $\mathcal C$.
At the specified assignment, $\bv=1$ and $yxz=0$, proving the reverse implication.

In $N_3$, a word of length at least two depends only on its first letter: a first value $1$ gives $0$, and first values $0,4$ are unchanged.
Thus, for $\ell(\bq)\ge2$, the identity $\bq\preceq \bu$ holds in $N_3$ exactly when $h(\bq)$ is the first letter of an original summand of length at least two.
Necessity follows by assigning $1$ to this letter and $4$ to all other letters; sufficiency follows from equal values of the two words of length at least two with the same first letter.
For a lower word of length one $x$, the condition is that $x$ itself is a summand or is the first letter of a summand of length at least two.
It follows that the criterion in $\mathcal J^*$ when $\ell(\bq)\ge2$ consists of (K1),(K2),(K4), together with the stated first-letter condition.
Each generator satisfies $\chi_1$, so $\mathcal J^*$ is contained in the indicated subvariety.

For the converse, suppose an identity satisfies this criterion.
If $E(\bu)\cap T(\bu)=\varnothing$, (K4) and (K1) imply that every summand whose content is contained in $c(\bq)$ has the same final letter as $\bq$.
ST and Lemma~\ref{lem:carrier-reorder} give the desired identity.
If $E(\bu)\cap T(\bu)\ne\varnothing$, construct $\mathbf{z}\preceq \bu$ with a repeated final letter and containing all original letters.
Using $\chi_A$, obtain the square of each original word of length at least two.
If there is a summand of length at least two whose content is contained in $c(\bq)$, start with its square and construct a square with content exactly $c(\bq)$.
If all such summands are single letters, take $h=h(\bq)$ and use the original summand $h\mathbf{v}_0$ of length at least two supplied by the $N_3$ condition;
MN gives $h^2\preceq h+h\mathbf{v}_0\preceq \bu$.
J2 now produces a square $\bt\preceq \bu$ with content exactly $c(\bq)$.

Apply H so that the resulting word begins with the required first letter.
For later first occurrences, an original summand of length at least two can be replaced by its square and then used with I;
if the content of an original summand is contained in the preceding content, use Blk as in Lemma~\ref{lem:carrier-reorder}.
The remaining case is that only a single-letter summand $b$ supplies the required occurrence.
Let $\mathbf{u}_0$ be the already arranged nonempty prefix, and write the current word as $\bt$, beginning with $\mathbf{u}_0$.
Substitute $x=b,y=\mathbf{u}_0,z=\bt$ in $\chi_1$.
Since $\bt$ has a repeated final letter, begins with $\mathbf{u}_0$, and contains $b$, its term side becomes
\[
 b+\bt^2+\bt{}\mathbf{u}_0^2+\mathbf{u}_0\bt{}b^2\approx b+\bt.
\]
The word $\mathbf{u}_0b\bt$ has $b$ immediately after $\mathbf{u}_0$, has the required content, and has a repeated final letter.
This completes the reordering.
Finally use L if the lower word has an unrepeated final letter.
For a lower word of length one, the assertion follows from the condition on powers of that letter imposed by $C_4$.
Hence all identities of $\mathcal J^*$ follow from the stated finite set.
\end{proof}

\begin{lem}\label{lem:C-omit-N3}
Put $F_0:x^2\preceq x+y^2$ and $\mathcal J_0=\Var(M_0,C_4)$.
For every $\mathcal W\subseteq\mathcal C$,
\[
 N_3\notin\mathcal W\quad\Longleftrightarrow\quad\mathcal W\models F_0,
 \qquad
 \mathcal C\cap\operatorname{Mod}(F_0)=\mathcal J_0.
\]
Thus $\Gamma\cup\{\chi_A,F_0\}$ is a finite basis for $\mathcal J_0$.
Over $\mathcal C$, the identity $F_0$ is equivalent to
\[
 x^2+y\approx x^2+y^2.
\]
Moreover, for a nonempty term $\bu$,
\[
 \mathcal J_0\models \bu\approx \bu^2
 \quad\Longleftrightarrow\quad E(\bu)\cap T(\bu)\ne\varnothing.
\]
\end{lem}
\begin{proof}
Evaluate $(x,y)=(1,4)$ in $N_3$ and use the principal down-set $I=\{1,4\}$.
The words $x^2,x,y$ represent $0,1,4$.
Membership in $I$, together with membership after right multiplication by $4$, separates distinct elements.
The eight two-variable normal forms are
\[
 x,\ y,\ x^2,\ y^2,\ xy,\ yx,\ xy^2,\ yx^2.
\]
Those evaluating into $I$ are $x,y,y^2,yx,yx^2$.
Each lies below $\bv=x+y^2$ in $\mathcal C$:
for $yx,yx^2$, the summand $y^2$ supplies the initial nonfinal $y$, and its content is contained in the set of letters preceding $x$.
The other conditions follow immediately from (K1)--(K4).
The remaining forms are $x^2,xy,xy^2$.
Each provides a summand of length at least two beginning with $x$, so
$x^2\preceq \bv+\br$ holds in $\mathcal C$ for each of them.
The same principal-down-set separation argument as before shows that
$N_3\notin\mathcal W$ implies $\mathcal W\models F_0$.
The chosen evaluation violates $F_0$, proving the converse.

The identity $F_0$ holds in both $M_0$ and $C_4$.
Consider any identity $\bq\preceq \bu$ holding in $\mathcal J_0$.
If $E(\bu)\cap T(\bu)=\varnothing$, combine the summands having the same final letter by ST and apply Lemma~\ref{lem:carrier-reorder}, as above.
If $E(\bu)\cap T(\bu)\ne\varnothing$, first produce a word $\mathbf{z}\preceq \bu$ with a repeated final letter.
The identity $\chi_A$ gives the square of every original summand of length at least two.
For a single-letter summand $a$, the instance
\[
 a^2\preceq a+\mathbf{z}^2\approx a+\mathbf{z}\preceq \bu
\]
of $F_0$ gives its square as well.
Hence $\bw^2\preceq \bu$ has been deduced for every original summand $\bw$.
The $C_4$ criterion supplies original summands covering the content of the lower word, and J2 combines their squares into a square with exactly that content.
Every original summand supplying a required first occurrence can now be replaced by its square, making that occurrence nonfinal.
Lemma~\ref{lem:carrier-reorder} and L give $\bq\preceq \bu$.
This proves the equality of varieties and the finite basis.

L gives $x^2+y\preceq x^2+y^2$; interchanging the variables in $F_0$ gives the reverse inequality.
Conversely, substituting $x\mapsto y$, $y\mapsto x$ in the displayed equality gives
$x+y^2\approx x^2+y^2$, which implies $F_0$.
Finally, if $E(\bu)\cap T(\bu)\ne\varnothing$, the preceding construction gives the square of every summand below $\bu$; additivity of the square map yields $\bu^2\preceq \bu$, and L gives the converse.
If the intersection is empty, evaluate all letters of $E(\bu)$ as $3$ and all letters of $T(\bu)$ as $1$ in $M_0$.
Then $\bu=1$ and $\bu^2=0$, so the equality fails.
\end{proof}
\section{Multiplicatively idempotent subvarieties}
\label{sec:bands}

In this section, we describe the multiplicatively idempotent subvarieties needed below.
We distinguish the variety of distributive lattices,
$\mathsf D_{\mathrm{lat}}$, from $\mathcal D=\Var(M_1,C_4)$.
Let $\mathsf L$ be generated by the two-element ai-semiring with left-zero
multiplication, and let $\mathsf M$ be generated by the two-element
ai-semiring in which both operations are join. The variety
$\mathsf D_{\mathrm{lat}}$ is generated by the two-element lattice, with
join as addition and meet as multiplication. Put
\[
 B=S|_{\{0,2,4\}},\qquad B_3=S|_{\{2,3,4\}},\qquad
 \mathsf B=\Var(B),\quad \mathsf S_{\mathrm{bi}}=\Var(B_3),\quad
 \mathsf K=\Var(C_4).
\]
For brevity, juxtaposed variety symbols denote joins; thus
$\mathsf{DM}=\mathsf D_{\mathrm{lat}}\vee\mathsf M$ and
$\mathsf{LDM}=\mathsf L\vee\mathsf D_{\mathrm{lat}}\vee\mathsf M$.

\begin{thm}[Pastijn; Pol\'ak]\label{thm:band-literature}
In the signature consisting of the two binary operations and no constants,
the subvarieties of the ai-semiring variety defined by $x^2\approx x$
form a finite distributive lattice. Each of these subvarieties is finitely
based \cite{Pastijn2005}. The order of its join-irreducible elements and their defining
identities are given in \cite[Result~4 and Theorem~7, pp.~831--834]{Polak2006}.
\end{thm}

The definition on p.~826 of \cite{Polak2006} does not require neutral
elements and therefore uses the same signature as the present paper.
The following lemma describes the interval needed below.

\begin{lem}\label{lem:band12}
We have $\mathsf K=\mathsf B\vee\mathsf S_{\mathrm{bi}}$.
There are exactly twelve subvarieties of $\mathsf K$:
\[
\begin{gathered}
\mathsf T,\quad \mathsf L,\quad \mathsf D_{\mathrm{lat}},\quad \mathsf M,
\quad \mathsf L\vee\mathsf D_{\mathrm{lat}},\quad \mathsf{LM},\quad
\mathsf{DM},\quad \mathsf{LDM},\\
\mathsf B,\quad \mathsf S_{\mathrm{bi}},\quad
\mathsf L\vee\mathsf S_{\mathrm{bi}},\quad \mathsf K,
\end{gathered}
\]
where $\mathsf T$ is the trivial variety. The maximal proper subvarieties
of $\mathsf K$ are $\mathsf B$ and $\mathsf L\vee\mathsf S_{\mathrm{bi}}$.
The largest proper subvariety of $\mathsf B$ is $\mathsf{LDM}$, and that
of $\mathsf S_{\mathrm{bi}}$ is $\mathsf{DM}$. Every proper subvariety of
$\mathsf{LDM}$ that contains $\mathsf M$ is contained in $\mathsf{LM}$
or $\mathsf{DM}$.
\end{lem}
\begin{proof}
The additive order of $B$ is $4<2<0$. Its multiplication has identity
$2$ and two left-zero elements $0,4$. This is the algebra denoted by
$B$ in the cited classification, with its elements $a,1,b$ corresponding
to $4,2,0$. The algebra $B_3$ is obtained from $\{2,3\}$, with both
operations equal to join, by adjoining $4$ as an additive least element
and a multiplicative zero. It therefore corresponds to $M^0$, whose
variety is denoted by $\mathsf S$ in that classification.

Both algebras are subalgebras of $C_4$. To prove the reverse inclusion, consider
the following subset of $B\times B_3$:
\[
 T_0=\{(0,2),(0,3),(0,4),(2,2),(2,3),(4,4)\}.
\]
Define
\[
 \pi(0,i)=0\ (i=2,3,4),\qquad
 \pi(2,i)=i\ (i=2,3),\qquad \pi(4,4)=4.
\]
For addition, an element with first coordinate $0$ added to any element
still has first coordinate $0$; adding $(4,4)$ leaves every other element
unchanged; and two elements with first coordinate $2$ are added by taking
the join of their second coordinates. It follows that this subset is closed under addition and that addition is
preserved by $\pi$.

For multiplication, a left factor with first coordinate $0$ gives a
product with first coordinate $0$, and the left factor $(4,4)$ gives
$(4,4)$. A left factor $(2,i)$ multiplied by a right factor of first
coordinate $0$, $4$, or $2$ gives, respectively, an element of the
corresponding one of these three classes. In the last case the second
coordinate is the product in $\{2,3\}$. Comparison with the table of $C_4$
proves closure and $\pi(uv)=\pi(u)\pi(v)$ in each case. Thus $\pi$ is a
surjective homomorphism onto $C_4$, and
$C_4\in\operatorname{HS}(B\times B_3)$. This proves the equality of joins.

In the partially ordered set of join-irreducible varieties in
Theorem~\ref{thm:band-literature}, those below
$\mathsf B\vee\mathsf S_{\mathrm{bi}}$ are precisely
\[
 \mathsf L,\quad\mathsf D_{\mathrm{lat}},\quad\mathsf M,\quad
 \mathsf B,\quad\mathsf S_{\mathrm{bi}}.
\]
The first three are pairwise incomparable. The element $\mathsf B$ is
above all three, whereas $\mathsf S_{\mathrm{bi}}$ is above
$\mathsf D_{\mathrm{lat}}$ and $\mathsf M$ only. The downward closed
subsets of this partially ordered set give exactly the twelve joins
listed above. The remaining assertions follow from the inclusions among these joins.
\end{proof}

\begin{thm}[Ren--Zeng]\label{thm:RenZeng}
Let
\[
 \mathcal I=\operatorname{Mod}_{\rm ai}(x^2\approx x),\qquad
 \mathcal Z=\operatorname{Mod}_{\rm ai}(xy\approx zt).
\]
Then $\mathcal I\vee\mathcal Z$ is hereditarily finitely based.
\end{thm}

Theorem~\ref{thm:RenZeng} is the hereditary finite basis theorem of
Ren and Zeng~\cite{RenZeng2019}; it applies to the displayed join.

\begin{lem}\label{lem:square-transfer}
Let $\mathcal A\subseteq\mathcal C$ be finitely based and satisfy
\[
 (x+y)^2\approx x^2+y^2,\qquad (xy)^2\approx x^2y^2,
 \qquad x^4\approx x^2.
\]
For every $\mathcal W\subseteq\mathcal A$, the variety
$\mathcal W\cap\mathcal I$ is generated by the square images of
the algebras in $\mathcal W$. Write
$A^{[2]}=\{a^2:a\in A\}$ for this image, to distinguish it from
the direct product $A\times A$. If $\mathcal U\subseteq\mathsf K$, then
\[
 R_{\mathcal U}=\{A\in\mathcal A:A^{[2]}\in\mathcal U\}
\]
is finitely based. Moreover, if $\bp \approx \bq$ holds in $\mathcal U$,
then $\bp \bt \approx \bq \bt$ holds in $R_{\mathcal U}$ for every nonempty
right factor $\bt$.
\end{lem}
\begin{proof}
The square map is an idempotent endomorphism, so its image is a
multiplicatively idempotent subalgebra in the original variety.
Conversely, every multiplicatively idempotent algebra is its own square
image. This proves the first assertion.

By Theorem~\ref{thm:band-literature}, choose a finite equational basis $\Delta$
for $\mathcal U$. Replace each variable in each identity of $\Delta$
by its square. These finitely many identities, together with a finite
equational basis for $\mathcal A$, define $R_{\mathcal U}$.

To prove the last assertion, apply an identity of $\mathcal U$ to the
squares and multiply on the right by $\bt$. Every inserted square then
has a nonempty factor to its right. The identity $x^2y\approx xy$,
valid in $\mathcal C$, removes these squares successively and gives
$\bp \bt \approx \bq \bt$. All substitutions used here are by nonempty terms.
\end{proof}
\section{Every subvariety of $\mathcal D$ is finitely based}
\label{sec:D-hfb}

In this section, we prove that $\mathcal D=\Var(M_1,C_4)$ is hereditarily finitely based.
The subsemirings of $S$ considered in this section retain the element labels in its Cayley tables. Put
\[
\begin{gathered}
 M_1=S|_{\{0,1,2,4\}},\qquad C_4=S|_{\{0,2,3,4\}},\qquad
 N=S|_{\{0,1,4\}},\\
 B=S|_{\{0,2,4\}},\qquad Z=S|_{\{0,1\}},\qquad
 S_b=S|_{\{2,3,4\}}.
\end{gathered}
\]
Thus $N$ is the semiring denoted by $N_3$ earlier, and $S_b=B_3$.
The element $2$ is a multiplicative identity of $B$, and $0,4$ are its left zero elements.
Let $L,M,D_2$ be the two-element semirings whose addition is join and whose multiplication is, respectively, left projection, join, and meet.
In particular, $D_2$ is the two-element distributive lattice semiring; it is not the variety $\mathcal D$.

We use the finite equational basis for $\mathcal C$, Lemmas~\ref{lem:C-omit-M0}--\ref{lem:C-omit-N3}, the characterization of subvarieties not containing $B_3$, and the subvariety lattice in Lemma~\ref{lem:band12}.
We also use Theorem~\ref{thm:RenZeng}.
These results do not depend on hereditary finite basability of $\mathcal C$ or $\Var(S)$.
For a variety $\mathcal X$ and an identity $\rho$, the notation $\mathcal X+\rho$ means $\mathcal X\cap\operatorname{Mod}(\rho)$.

\subsection{Identities of the relevant subsemirings}

The following identities hold in $\mathcal D$:
\begin{equation}\label{eq:D-word}
\begin{gathered}
 x^3\approx x^2,\qquad x^2y\approx xy\approx xy^2,\qquad xyx\approx xy,\\
 (xy)^2\approx xy,\qquad x\preceq x^2,\qquad
 (x+y)^2\approx x^2+y^2.
\end{gathered}
\end{equation}
The square map also preserves multiplication. These assertions follow from the word-evaluation formulas for $M_1$ and $C_4$.
Consequently, if $\bw$ has length at least two, then $\bw \approx \bw^2$ in $\mathcal D$.
The same is true of a term all of whose summands have length at least two.

Let $\bu$ be a nonempty term, regarded also as its finite set of summands.
Let $L(\bu)$ be the set of variables occurring as one-letter summands, and let $H(\bu)$ be the set of first letters of summands of length at least two.
For a word $\bq$, define the following conditions:
\begin{align*}
\mathrm{WA}(\bu,\bq):\quad&\text{there is }\bp \in \bu \text{ with }c(\bp)\subseteq c(\bq);\\
\mathrm{SA}(\bu,\bq):\quad&\text{there is }\bp \in \bu \text{ with }\ell(\bp)\ge2\text{ and }c(\bp)\subseteq c(\bq),\\
&\text{or }c(\bq)\cap L(\bu)\cap H(\bu)\ne\varnothing.
\end{align*}
We use these symbols as abbreviations for the displayed conditions.
For each letter $b$ of $\bq$, let $A_b$ be the set of letters occurring before its first occurrence in $\bq$.
We shall also use
\begin{align*}
\mathrm{Ex}(\bu,\bq):\quad&\text{for each }b\in c(\bq),\text{ some }\bp \in \bu
\text{ has }b\text{ as its first letter outside }A_b;\\
\mathrm{Cov}(\bu,\bq):\quad&c(\bq)=\bigcup\{c(\bp):\bp \in \bu,\ c(\bp)\subseteq c(\bq)\}.
\end{align*}
The first condition allows $b$ to be the last letter of the summand.

For $\ell(\bq)\ge2$, the conditions for $\bq \preceq \bu$ in several semirings are as follows:
\[
\begin{array}{c|l}
B&\mathrm{WA}(\bu,\bq)\text{ and }\mathrm{Ex}(\bu,\bq),\\
S_b&\mathrm{Cov}(\bu,\bq),\\
N&h(\bq)\in H(\bu),\\
C_4&\mathrm{Cov}(\bu,\bq)\text{ and }\mathrm{Ex}(\bu,\bq).
\end{array}
\]
Indeed, in $B$ a word takes the value of its first letter assigned $0$ or $4$, and takes the value $2$ if all its letters are assigned $2$.
For $S_b$ the two proper principal down-sets of its additive order give the displayed content condition.
In $N$, every word of length at least two has the value of the square of its first letter.
The $C_4$ condition also follows by restricting the word-evaluation formula for $S$.
Here and below, a \emph{principal down-set} is a set of the form $\{a:a\leq d\}$ in the natural semilattice order.

For $M_1$ and $\ell(\bq)\ge2$, the necessary and sufficient conditions are
\begin{equation}\label{eq:M1-criterion}
\begin{gathered}
\mathrm{WA}(\bu,\bq),\qquad\mathrm{Ex}(\bu,\bq),\\
\text{for each }b\in c(\bq),\quad
\left\{
\begin{array}{l}
\text{some }\bp \in \bu \text{ satisfies }c(\bp)\subseteq A_b,\quad\text{or}\\
\text{some }\bp \in \bu \text{ has }\ell(\bp)\ge2,\text{ and }b\\
\qquad\text{is its first letter outside }A_b.
\end{array}\right.
\end{gathered}
\end{equation}
The last condition is imposed also when the first occurrence of $b$ is the final position of $\bq$.
It requires a summand of length at least two, not that the occurrence of $b$ in that summand be nonfinal.

We now verify \eqref{eq:M1-criterion}.
A word of length at least two over $M_1$ takes the value $2$ exactly when all letters are assigned $2$.
Otherwise consider its first letter not assigned $2$: the value is $4$ if this letter is assigned $4$, and is $0$ if it is assigned $0$ or $1$.
The proper principal down-sets are $\{4\}$, $\{1,4\}$ and $\{2,4\}$.
Assigning the letters of $c(\bq)$ to $2$ and all other letters to $4$ proves the necessity of $\mathrm{WA}$.
For a particular $b$, assign $A_b$ to $2$, assign $b$ first to $0$ and then to $1$, and assign all remaining letters to $4$.
The first assignment proves $\mathrm{Ex}$; the second proves the last condition when there is no summand supported in $A_b$.
Conversely, a value $4$ of $\bq$ is below every element, and the case in which the value is $2$ follows from $\mathrm{WA}$.
If the value is $0$, apply the conditions at the first letter not assigned $2$.
An appropriate summand has value $0$, unless it is the one-letter word assigned $1$.
In the latter case the last condition supplies either a longer such summand, again of value $0$, or a summand supported in $A_b$, of value $2$; the latter together with $1$ has sum $0$.
This proves sufficiency.

An identity holds in a join of varieties if and only if it holds in every constituent variety.
For a word $x$ of length one on the left-hand side, the condition in each of $B,C_4,D_2$ is that $\bu$ have a summand supported on $\{x\}$.
Necessity follows by assigning $x$ to $2$ and all other letters to $4$ in $B$, by the same assignment in $C_4$, and by assigning $x$ to the larger element and the other letters to the smaller element in $D_2$.
Sufficiency follows from $x\preceq x^2$ and \eqref{eq:D-word}.
Thus the same condition applies to a join containing one of these factors.
We verify it directly for $U$ below; for $T$ it follows from $M_1\in\Var(T)$ and $B\subseteq M_1$.

By Lemma~\ref{lem:C-omit-M0}, $\mathcal D$ is finitely based: add $xy\approx(xy)^2$ to a finite basis for $\mathcal C$.

\subsection{Subvarieties not containing $N$ or $M_1$}

\begin{lem}\label{lem:D-omit-N}
Let $f$ be the identity $x^2\preceq x+y^2$.
For every $\mathcal W\subseteq\mathcal D$,
\[
 N\notin\mathcal W\iff\mathcal W\models f,
 \qquad \mathcal D+f=\Var(C_4,Z).
\]
Every subvariety of $\mathcal D$ not containing $N$ is therefore finitely based.
\end{lem}
\begin{proof}
The equivalence follows from Lemma~\ref{lem:C-omit-N3}. We give a direct proof for arbitrary identities.
In $N$, take $(x,y)=(1,4)$ and the principal down-set $I=\{1,4\}$.
The words $x^2,x,y$ represent $0,1,4$.
Membership in $I$, together with membership in $I$ after right multiplication by $1$, distinguishes these three elements.
By \eqref{eq:D-word}, every word on $x,y$ is equivalent to one of
\[
x,\quad y,\quad x^2,\quad xy,\quad yx,\quad y^2.
\]
The representatives taking values in $I$ are $x,y,yx,y^2$; their sum equals $\bu_0=x+y^2$ in $\mathcal D$.
For each of the two remaining representatives $\bp =x^2,xy$, the identity $x^2\preceq \bu_0+\bp$ holds in $\mathcal D$, as follows from the conditions just proved.
Choose an identity of $\mathcal W$ that fails in $N$.
Represent the values in an assignment for which it fails by $x^2,x,y$, apply one of the two separating maps, and interchange its sides if necessary.
This gives an identity $\bs \approx \bt$ of $\mathcal W$ for which every summand of $\bs$ takes its value in $I$, whereas some summand of $\bt$ does not.
Here we use the fact that a principal down-set in a join semilattice contains a finite join exactly when it contains every joinand.
Consequently $\bu_0+\bs \approx \bu_0$ and $x^2\preceq \bu_0+\bt$ hold in $\mathcal D$.
Together with $\bs \approx \bt$ they imply $f$.
The assignment $(x,y)=(1,4)$ shows that $N$ fails $f$.

In $\mathcal D+f$, if $\bu$ has a summand $\bw$ of length at least two, then each one-letter summand $x$ satisfies
$x^2\preceq x+\bw^2\approx x+\bw \preceq \bu$.
All summands of length at least two already equal their squares, so square additivity gives $\bu \approx \bu^2$.
The square image of $M_1$ is $B\subseteq C_4$, and that of $C_4$ is $C_4$ itself.
Hence every identity of $C_4$, with all variables replaced by their squares, holds in $\mathcal D$.
In $Z$, a term with a summand of length at least two is constantly $0$, the largest element of its additive order.
A sum of one-letter words is determined by its variable set, and is never identically $0$.
Thus an identity of $\Var(C_4,Z)$ either has the same sum of one-letter words on both sides, or both sides have a summand of length at least two.
In the second case square both sides and apply the corresponding identity in the square image.
This proves $\mathcal D+f=\Var(C_4,Z)$; the reverse inclusion also follows directly from the tables.
Since $C_4$ has idempotent multiplication and $Z$ has constant multiplication, Theorem~\ref{thm:RenZeng} applies.
\end{proof}

Put
\[
\mathcal J=\Var(N,C_4),\qquad
\chi_1:\ yxz\preceq x+z^2+zy^2+yzx^2.
\]
Then $\mathcal J$ is the largest subvariety of $\mathcal D$ not containing $M_1$.
Indeed, Lemma~\ref{lem:C-omit-M1} gives
$\mathcal C+\chi_1=\Var(M_0,C_4,N)$.
Adding $xy\approx(xy)^2$ gives exactly $\mathcal J$, as we now verify.
Suppose $\bq \preceq \bu$ holds in $N,C_4$, with $\ell(\bq)\ge2$.
There is a summand $\bw$ of $\bu$ of length at least two with $h(\bw)=h(\bq)$.
Replace this summand by $\bw^2$.
The replacement preserves values in $N,C_4$, and the resulting term $\bu'$ satisfies $E(\bu')\cap T(\bu')\ne\varnothing$.
Thus the condition for $M_0$ reduces to $c(\bq)\subseteq c(\bu')$, which follows from $\mathrm{Cov}(\bu,\bq)$.
The identity $\bq \preceq \bu'$ therefore holds in $M_0,C_4,N$.
In the variety with the additional identity $xy\approx(xy)^2$, we have $\bu \approx \bu'$, proving the original identity.
If the word on the left has length one, the required summand is a power of that variable.
This proves the asserted equality, the characterization of $M_1\notin\mathcal W$, and a finite basis for $\mathcal J$.

\subsection{Finite bases for $M_1$ and $\Var(B,N)$}

We require a finite equational basis for $M_1$.
The semiring $M_1$ is isomorphic to $S_{(4,427)}$ in Yue--Ren--Zeng--Shao~\cite{Yue2025}:
our labels $0,4,1,2$ correspond, respectively, to their labels $1,2,3,4$.
Proposition~8.2, equations~(111)--(120), of the authors' version gives the following equational basis in addition to the ai-semiring axioms:
\begin{gather}\label{eq:M1-basis}
 x^3\approx x^2,\quad xyx\approx xy,\quad x^2y\approx xy,\quad xy^2\approx xy,\\
 xy\approx(xy)^2,\quad x\preceq x^2,\quad x+xy\approx x^2+xy,\nonumber\\
 yx\preceq x^2+yz,\quad yx\preceq x+yxz,\quad cbad\preceq a+b+cbd.\nonumber
\end{gather}
The paper is published in \emph{Algebra Universalis} \textbf{86}, Article~33 (2025).
We use the displayed equational basis below.

Let $k$ denote $xy\preceq y+xz$, and put $\mathcal K=\Var(M_1)+k$.
Then
\begin{equation}\label{eq:K-NB}
\mathcal K=\Var(N,B).
\end{equation}
To prove this, assume $H=H(\bu)\ne\varnothing$, and put $L_0=L(\bu)\setminus H$.
For $h\in H$, choose a summand $h\br$ of $\bu$.
For every summand $\bp$ of $\bu$, the identity $k$ gives $h\bp \preceq \bp +h\br \preceq \bu$.
Moreover, $h(h\br)\approx h\br$.
Put $\mathbf{h}=\sum_{h\in H}h$. We obtain
\begin{equation}\label{eq:K-normal}
\bu \approx\sum_{x\in L_0}x+\mathbf{h}\bu^2.
\end{equation}
If the first sum is empty, it is omitted; no zero constant is introduced.
More explicitly, $\bu \approx \sum_{x\in L(\bu)}x+\mathbf{h}\bu$; the identity $xy^2\approx xy$ replaces $\mathbf{h}\bu$ by $\mathbf{h}\bu^2$.
If $h\in L(\bu)\cap H$, then $h^3$ is a summand of $\mathbf{h}\bu^2$. Since $h^3\approx h^2$ and $h\preceq h^2$, the separate summand $h$ may be omitted.
In $N$, a term is determined by the pair $(L_0,H)$: it equals the sum of the letters in $L_0$ and the squares of the letters in $H$.
Assigning one variable to $1$ and all other variables to $4$ distinguishes its absence, its membership in $L_0$, and its membership in $H$, by the respective values $4,1,0$.

If $\bu \approx \bv$ holds in $N,B$, the pairs $(L_0,H)$ coincide.
Also $\bu^2\approx \bv^2$ holds in $M_1$, since its square map is an onto homomorphism to $B$.
If $H$ is empty, $\bu,\bv$ are the same one-letter sum; otherwise \eqref{eq:K-normal} proves $\bu \approx \bv$.
This proves \eqref{eq:K-NB}, and \eqref{eq:M1-basis} together with $k$ is a finite basis for $\mathcal K$.
The identities $xy\approx x^2$, $x^3\approx x^2$, $x\preceq x^2$ form a finite basis for $N$ itself, by the same normal-form argument.

\subsection{An auxiliary five-element semiring}

Let $U$ be the algebra with the following tables; its labels are independent of the labels of $S$:
\[
\begin{array}{c|ccccc}+&0&1&2&3&4\\\hline
0&0&0&0&0&0\\1&0&1&0&0&1\\2&0&0&2&0&2\\3&0&0&0&3&3\\4&0&1&2&3&4
\end{array}\qquad
\begin{array}{c|ccccc}\cdot&0&1&2&3&4\\\hline
0&0&0&0&3&3\\1&0&0&0&3&3\\2&2&2&2&4&4\\3&3&3&3&3&3\\4&4&4&4&4&4
\end{array}.
\]
Its square image is $\{0,2,3,4\}\cong L\times D_2$, under
\[
0\longleftrightarrow(1,1),\quad2\longleftrightarrow(0,1),\quad
3\longleftrightarrow(1,0),\quad4\longleftrightarrow(0,0).
\]
On this image addition is coordinatewise join and multiplication is
$(a,b)(c,d)=(a,b\wedge d)$.
The row and column of $1$ in the multiplication table are the same as those of $0$, and $1^2=0$.
The map $r$ sending $1$ to $0$ and fixing the other elements preserves addition, and $uv=r(u)r(v)$ for all $u,v$.
The value of every product lies in $L\times D_2$.
Associativity and distributivity therefore follow from that product semiring; the displayed addition is a join-semilattice operation.
This proves that $U$ is an ai-semiring.

\begin{lem}\label{lem:D-U-criterion}
A finite equational basis for $\Var(U)$ consists of the ai-semiring axioms and
\begin{gather}\label{eq:U-basis}
xyz\approx xzy,\quad x^2y\approx xy,\quad xy^2\approx xy,\quad xyz\preceq xy,\\
xyz\preceq xt+yz,\quad x^2\preceq x+xy,\quad x\preceq x^2.\nonumber
\end{gather}
For $\ell(\bq)\ge2$, the identity $\bq \preceq \bu$ holds in $U$ if and only if
$h(\bq)\in H(\bu)$ and $\mathrm{SA}(\bu,\bq)$.
For $\bq =x$ it holds if and only if some summand of $\bu$ is a power of $x$.
\end{lem}
\begin{proof}
The first three identities give a normal form for words of length at least two: such a word is determined by its first letter and its content.
One first permutes letters after the initial letter and then deletes repetitions; a word supported on a single letter is reduced to its square, not to a one-letter word.
The remaining identities follow from the coordinate description.
For identities involving $1$, note specifically that $1^2=0$ and $1+1v=0$ for every $v$; the other cases lie in the square image.

If $h=h(\bq)\notin H(\bu)$, assign $h$ to $1$ and the other letters to $4$.
Then $\bu$ takes a value in $\{1,4\}$, while $\bq$ takes the value $0$ or $3$, so the identity fails.
Suppose the first-letter condition holds but $\mathrm{SA}$ fails.
Assign the letters of $c(\bq)\cap L(\bu)$ to $1$, those of $c(\bq)\cap H(\bu)$ to $2$, the remaining letters of $c(\bq)$ to $0$, and all letters outside $c(\bq)$ to $4$.
The first two sets are disjoint.
Every summand of length at least two has a letter outside $c(\bq)$, so its second coordinate after applying $r$ is zero; its initial letter has first coordinate zero.
Its value is therefore $4$.
Every one-letter summand takes the value $1$ or $4$.
The word $\bq$ takes the value $2$, proving failure again.

Conversely, choose a summand $h\mathbf{v}_0$.
If $\bw =\bs \bt$ is a summand of length at least two with $c(\bw)\subseteq c(\bq)$, then $\bs,\bt$ are nonempty and the fifth identity gives
$h\bw \preceq h\mathbf{v}_0 +\bw \preceq \bu$.
In the other case of $\mathrm{SA}$, choose $a\in c(\bq)\cap L(\bu)\cap H(\bu)$ and a summand $a\mathbf{w}_0$.
The sixth identity gives $a^2\preceq a+a\mathbf{w}_0 \preceq \bu$, and the fifth then gives $ha^2\preceq h\mathbf{v}_0 +a^2\preceq \bu$.
Thus in either case there is a word of length at least two, beginning with $h$ and supported in $c(\bq)$, whose inequality with $\bu$ follows from the displayed identities.
The fourth identity permits appending every missing letter of $\bq$.
The first three identities then give $\bq \preceq \bu$.
All substitutions in this construction use nonempty words.
For a word $x$ of length one on the left-hand side, necessity follows by assigning $x$ to $2$ and all other letters to $4$; sufficiency follows from $x\preceq x^2$.
Every term identity is equivalent to finitely many inequalities with a single word on the left, one for each summand on either side.
Hence the displayed identities form an equational basis.
\end{proof}

The subsemiring $U|_{\{0,1,3,4\}}$ has a congruence with classes $\{0,3\}$, $\{1\}$, $\{4\}$; its quotient is $N$.
Also $L,D_2\in\Var(U)$ by the square-image description.
Finally $U\in\mathcal J$.
Indeed, an identity $\bq \preceq \bu$ of $\mathcal J$ with $\ell(\bq)\ge2$ satisfies $h(\bq)\in H(\bu)$ and $\mathrm{Cov}(\bu,\bq)$.
If none of the summands supported in $c(\bq)$ has length at least two, every letter of $\bq$, including $h(\bq)$, is a one-letter summand.
The first-letter condition then implies $\mathrm{SA}(\bu,\bq)$.
Lemma~\ref{lem:D-U-criterion} proves the identity in $U$; the case in which the word on the left has length one has already been treated.

\subsection{Descriptions of $\mathcal J\cap\mathcal E$ and $\mathcal D\cap\mathcal E$}

Put
\[
g:\ xy\preceq y^2+xz,\qquad e:\ xy\preceq x+xyz,\qquad
\mathcal J_E=\mathcal J+e,\qquad\mathcal T=\mathcal D+e.
\]
Over $\mathcal D$, the identities $e,g$ are equivalent.
Left multiplication of $g(y,x,z)$ by $x$ gives
$xy\preceq x^2+xyz\approx x+xyz$.
For the equality here, $C_4$ is multiplicatively idempotent; in $M_1$, the two sides both equal $0$ if $x=1$, and $x^2=x$ otherwise.
Conversely, put $\bu =y^2+xz$.
The identity $\mathrm{J2}$ from Lemma~\ref{lem:helpers} gives $yxz\preceq \bu$; also $y\preceq \bu$.
Applying $e(y,x,z)$ yields $yx\preceq \bu$.
The identity $\mathrm H$ gives $xy\preceq yx+xz\preceq \bu$.

\begin{lem}\label{lem:D-JE}
We have $\mathcal J_E=\Var(B,U)$.
For $\ell(\bq)\ge2$, its identities $\bq \preceq \bu$ are characterized by
$\mathrm{SA}(\bu,\bq)$, $\mathrm{Ex}(\bu,\bq)$, and $h(\bq)\in H(\bu)$.
\end{lem}
\begin{proof}
Both $B,U$ belong to $\mathcal J$ and satisfy $g$.
Conversely, consider an identity of their join.
By condition $\mathrm{SA}$, there is a word $\bs$ of length at least two, supported in $c(\bq)$, with $\bs \preceq \bu$ derivable in $\mathcal J$.
Take an appropriate summand if it has length at least two; otherwise use $a^2\preceq a+a\mathbf{v}_0$ for the indicated one-letter summand $a$ and longer summand $a\mathbf{v}_0$.
In either case $\bs \approx \bs^2$.
For each first occurrence of a letter $b$ in $\bq$, choose a summand $\mathbf{v}_0 b\mathbf{w}_0$ whose first letter outside $A_b$ is this $b$.
If $\mathbf{w}_0$ is empty, the summand is already supported in $c(\bq)$.
If $\mathbf{w}_0$ is nonempty, the substitution $x\mapsto \mathbf{v}_0 b$, $y\mapsto \bs$, $z\mapsto \mathbf{w}_0$ in $g$ yields
\[
\mathbf{v}_0 b\bs \preceq \bs^2+\mathbf{v}_0 b\mathbf{w}_0 \preceq \bu.
\]
The added word contains $b$ and is supported in $c(\bq)$.
After adding these derived summands, condition $\mathrm{Cov}$ holds.
The original summands still ensure $\mathrm{Ex}$ and $h(\bq)\in H(\bu)$.
The modified identity therefore holds in $N,C_4$, and follows from a basis for $\mathcal J$.
Each new summand is already below $\bu$, so the original identity follows in $\mathcal J+g$.
The case in which the word on the left has length one follows from the condition already established. This completes the proof.
In particular, adjoining $e$ to a finite basis for $\mathcal J$ gives a finite basis for this join.
\end{proof}

\begin{lem}\label{lem:D-T-criterion}
For $\ell(\bq)\ge2$, an identity $\bq \preceq \bu$ holds in $\mathcal T$ if and only if it satisfies \eqref{eq:M1-criterion} and $\mathrm{SA}(\bu,\bq)$.
Moreover, $\mathcal T$ is generated by the quotient semiring $T$ associated with the language
\[
\mathscr L=\{y,z\}\ \cup\ x^+y\{x,y,z\}^*z\{x,y,z\}^*
\]
by the construction given below.
\end{lem}
\begin{proof}
Let $X=\{x,y,z\}$, and let $F$ be the semiring of nonempty finite subsets of $X^+$, with union as addition and setwise concatenation as multiplication.
Its elements are precisely the nonempty terms on $X$, modulo the ai-semiring axioms.
For $\bu,\bv \in F$, define
\[
\bu \equiv \bv \quad\Longleftrightarrow\quad
\bigl(\mathbf{s}_0  \bu \mathbf{t}_0 \subseteq\mathscr L\iff\mathbf{s}_0  \bv \mathbf{t}_0 \subseteq\mathscr L\bigr)
\quad\text{for every }\mathbf{s}_0,\mathbf{t}_0 \in X^*.
\]
An empty $\mathbf{s}_0$ or $\mathbf{t}_0$ merely omits that factor and is not a constant in the algebraic signature.
The relation is an equivalence. It is compatible with addition, since all words in a union belong to the language if and only if all words in each of its two terms belong to the language.
It is compatible with multiplication because, for a term $\br =\sum_i \br_i$,
$\mathbf{s}_0  \bu \br \mathbf{t}_0 \subseteq\mathscr L$ holds exactly when $\mathbf{s}_0  \bu (\br_i\mathbf{t}_0)\subseteq\mathscr L$ holds for every $i$; left multiplication is analogous.
Thus $\equiv$ is a congruence. Define $T=F/{\equiv}$.
Thus $T$ is the syntactic semiring associated with the displayed language. We shall use its defining congruence directly.

Put $\bv =y+z+xyz$.
The conditions for $S$ in Lemma~\ref{lem:S-criterion} give
\begin{equation*}
\mathscr L=\{\bw \in X^+:S\models \bw \preceq \bv \}.
\end{equation*}
Indeed, the one-letter possibilities are $y,z$.
A longer word must start with $x$, contain both $y,z$ by the condition on the union of the contents of the relevant summands, and have its first $y$ before its first $z$.
These conditions are also sufficient by the same lemma.
It follows that $T\in\Var(S)$.
To see explicitly that arbitrary term substitutions are covered, take any identity of $S$ and substitute elements of $F$ for its variables.
For any fixed $\mathbf{s}_0,\mathbf{t}_0$, its two substituted sides have the same value as terms over $S$, and hence are either both below $\bv$ or both not below $\bv$.
Being below $\bv$ is equivalent to every summand belonging to $\mathscr L$.
Thus the identity holds in the quotient.

We next prove that $T\models xy\approx(xy)^2$.
Every summand of a substituted product $\mathbf{u}_1 \mathbf{u}_2$ is a word of length at least two.
Let $\bp_i$ be these words, and fix $\mathbf{s}_0,\mathbf{t}_0 \in X^*$.
Suppose every $\mathbf{s}_0  \bp_i \mathbf{t}_0$ belongs to $\mathscr L$.
All these words have length at least two, so they belong to the second part of the displayed language: they start with $x$, contain $y,z$, and have no $z$ before the first $y$.
For any nonempty sequence of indices, the word
$\mathbf{s}_0  \bp_{i_1}\cdots \bp_{i_k}\mathbf{t}_0$ has the same three properties.
For the order condition, if $\mathbf{s}_0$ already contains $y$, its prefix through that first $y$ is unchanged.
If $\mathbf{s}_0$ does not contain $y$, then $\mathbf{s}_0$ contains only $x$, and no $\bp_i$ can contain $z$ before its first $y$: such a $z$ would already make $\mathbf{s}_0  \bp_i \mathbf{t}_0$ fail.
The first-letter condition is unchanged, and at least one selected block together with $\mathbf{s}_0,\mathbf{t}_0$ supplies both $y$ and $z$.
This proves preservation under every concatenation, including the cross terms $\bp_i \bp_j$.
Conversely, $\mathbf{s}_0  \bp \bp  \mathbf{t}_0$ and $\mathbf{s}_0  \bp \mathbf{t}_0$ have the same first letter, content, and relative order of the first $y,z$; thus they have the same membership in $\mathscr L$.
If all cross terms are in the language, the diagonal terms $\bp_i \bp_i$ therefore imply that all $\mathbf{s}_0  \bp_i \mathbf{t}_0$ are in it.
This proves the identity for all term substitutions and all contexts.
It fails in $A_4$, so $A_4\notin\Var(T)$.
The earlier characterization of subvarieties not containing $A_4$ gives $T\in\mathcal C$, and then $T\in\mathcal D$.

The quotient also satisfies $e$.
For arbitrary nonempty terms $\mathbf{u}_1,\mathbf{u}_2,\mathbf{u}_3$, suppose that every summand of
$\mathbf{s}_0(\mathbf{u}_1 +\mathbf{u}_1 \mathbf{u}_2 \mathbf{u}_3)\mathbf{t}_0$ belongs to $\mathscr L$.
Choose arbitrary summands $\mathbf{a} \in \mathbf{u}_1,\mathbf{b} \in \mathbf{u}_2$, and any $\mathbf{c} \in \mathbf{u}_3$.
Both $\mathbf{s}_0  \mathbf{a} \mathbf{t}_0$ and $\mathbf{s}_0  \mathbf{a} \mathbf{b} \mathbf{c} \mathbf{t}_0$ are in $\mathscr L$.
The former cannot be the one-letter word $y$ or $z$: then $\mathbf{s}_0,\mathbf{t}_0$ would be empty and $\mathbf{a}$ that one letter, whereas the latter word would begin with that letter and have length greater than one.
Hence $\mathbf{s}_0  \mathbf{a} \mathbf{t}_0$ starts with $x$ and contains $y,z$ in the required order.
The word $\mathbf{s}_0  \mathbf{a} \mathbf{b} \mathbf{t}_0$ still starts with $x$ and still contains those letters.
If the inserted block $\mathbf{b}$ introduced a $z$ before the first $y$, the same initial segment $\mathbf{s}_0  \mathbf{a} \mathbf{b}$ would make $\mathbf{s}_0  \mathbf{a} \mathbf{b} \mathbf{c} \mathbf{t}_0$ fail the order condition.
Thus $\mathbf{s}_0  \mathbf{a} \mathbf{b} \mathbf{t}_0 \in\mathscr L$.
This holds for every $\mathbf{a},\mathbf{b}$, proving $\mathbf{u}_1 \mathbf{u}_2 \preceq \mathbf{u}_1 +\mathbf{u}_1 \mathbf{u}_2 \mathbf{u}_3$ in the quotient.
In particular $T\in\mathcal T$.

We also have $M_1\in\Var(T)$.
Represent its elements $1,2,4$ by variables $p,q,r$, and its element $0$ by $p^2$.
Define a substitution by
\[
\sigma(p)=z,\qquad\sigma(q)=x,\qquad\sigma(r)=xyz.
\]
For every word $\bw$ on $p,q,r$, direct evaluation gives
\[
\begin{array}{c|c}
\bw (1,2,4)\in\{1,4\}&\sigma(\bw)\in\mathscr L\\
\bw (1,2,4)\in\{2,4\}&\sigma(\bw)yz\in\mathscr L\\
\bw (1,2,4)=4&x\sigma(\bw)\in\mathscr L.
\end{array}
\]
We verify these equivalences for every word.
A word consisting entirely of $q$ has value $2$; its image is a positive power of $x$, and exactly the second membership condition holds.
If the first non-$q$ letter is $r$, the value is $4$; the image begins with a positive string of $x$ followed by $yz$, and all three membership conditions hold.
The one-letter word $p$ has value $1$ and image $z$, so exactly the first membership condition holds.
Every other word whose first non-$q$ letter is $p$ has value $0$; its image has a $z$ before any $y$, or starts with $z$ and has length at least two, and none of the membership conditions holds.
These four cases exhaust all words.
Each of the three sets on the left is a principal down-set, so the equivalences also hold for terms with universal membership of their summands on the right.
The three membership conditions distinguish the four elements of $M_1$.
Given an identity failing in $M_1$, represent the values of its variables by $p,q,r,p^2$ and apply a distinguishing row of the table.
The corresponding substitution and context distinguish its two sides in $T$.
Thus every identity of $T$ holds in $M_1$, proving $M_1\in\Var(T)$.

It remains to establish the stated criterion.
The $M_1$ conditions are necessary by this inclusion.
If $\mathrm{SA}(\bu,\bq)$ fails, assign each variable of $c(\bq)\cap L(\bu)$ the word $y$, each other variable of $c(\bq)$ the word $x$, and each variable outside $c(\bq)$ the word $xyz$.
The image of $\bq$ is a word of length at least two without $z$, so is not in $\mathscr L$.
Every summand of $\bu$ of length at least two contains an outside variable, because the first alternative in $\mathrm{SA}$ fails.
Its initial variable is not in $c(\bq)\cap L(\bu)$, because the second alternative fails.
Its image therefore starts with $x$, contains an $xyz$ block, and has its first $y$ before its first $z$.
One-letter summands have images $y$ or $xyz$.
Hence all images of the summands of $\bu$ belong to $\mathscr L$, and the identity fails in $T$.

Conversely, work in $\mathcal D+g$ and assume the stated conditions.
As in Lemma~\ref{lem:D-JE}, use $\mathrm{SA}$ to obtain $\bs \approx \bs^2\preceq \bu$ with $\ell(\bs)\ge2$ and $c(\bs)\subseteq c(\bq)$.
For each required first occurrence use $g$ to add a derived summand $\mathbf{v}_0 b\bs$, unless the original summand already ends at $b$.
This ensures $\mathrm{Cov}$, preserves every condition for $M_1$, and supplies all conditions for $C_4$.
The resulting inequality follows from the finite basis for $\mathcal D$; deleting the derived summands gives $\bq \preceq \bu$.
For a word of length one on the left, use the condition that a summand is a power of that variable; its necessity follows from $B\subseteq M_1\in\Var(T)$.
Thus every identity of $T$ follows in $\mathcal T$, while $T\in\mathcal T$ was proved above.
Consequently $\mathcal T=\Var(T)$.
\end{proof}

\subsection{Three intervals and hereditary finite basability of $\Var(U)$}

For $R$ equal to a one-element semiring, $M$, or $B$, set
\[
\mathcal U_R=\Var(U,R),\qquad
\mathcal K_R=\Var(N,D_2,R).
\]
The following lemma describes three intervals in the subvariety lattice.

\begin{lem}\label{lem:D-U-cover}
If $\mathcal U_R$ is finitely based, then
\[
\mathcal U_R+k=\mathcal K_R,\qquad
[\mathcal K_R,\mathcal U_R]=\{\mathcal K_R,\mathcal U_R\}.
\]
For $R=B$, we have $D_2\in\Var(B)$ and $\mathcal K_B=\mathcal K$.
\end{lem}
\begin{proof}
For a word of length at least two on the left, the condition in $\mathcal K_R$ is
$h(\bq)\in H(\bu)$ and $\mathrm{WA}(\bu,\bq)$, together with the condition in $R$.
Each factor satisfies $k$.
If the only available summand supported in $c(\bq)$ is a one-letter summand $a$, choose an original summand $h(\bq)\mathbf{v}_0$ of length at least two.
Then $k$ gives $h(\bq)a\preceq a+h(\bq)\mathbf{v}_0 \preceq \bu$.
The modified term now satisfies $\mathrm{SA}$, and the added inequality also holds in $R$.
The criterion for $U$ and that for $R$ prove the original identity.
The case of a word of length one follows from the condition already established.
This proves the equality. A finite equational basis is obtained by adjoining $k$.

Take $\mathcal K_R\subseteq\mathcal W\subsetneq\mathcal U_R$, and an identity $\rho$ of $\mathcal W$ failing in $\mathcal U_R$.
Since $R\in\mathcal W$, the failure must occur in $U$.
Let $I=\{1,4\}$, a principal down-set of the additive order of $U$.
The table records membership in $I$ after applying each of four maps:
\[
\begin{array}{c|cccc}
a&a&a3&2a&a2\\\hline
0&0&0&0&0\\1&1&0&0&0\\2&0&1&0&0\\3&0&0&1&0\\4&1&1&1&1
\end{array}.
\]
Since the rows are distinct, these maps distinguish any two unequal values.
Under $(x,y,z)=(2,1,3)$, the elements $0,1,2,3,4$ are represented by $y^2,y,x,z,xz$.
Substitute these representatives into the failing assignment of $\rho$, apply a separating map using words for its fixed factors, and orient the two sides appropriately.
We obtain an identity $\bs \approx \bt$ of $\mathcal W$ such that all summands of $\bs$ take values in $I$ and some summand of $\bt$ does not.
Under the displayed assignment, a word takes its value in $I$ precisely when it is $y$, or when it starts with $x$ and contains $z$.
Put $\bu_0=y+xz$.
Every such word satisfies $\bw \preceq \bu_0$ in $U$ by Lemma~\ref{lem:D-U-criterion}, and in $M$ because its content is contained in $\{x,y,z\}$.
It also does so in $B$: the summand $xz$ supplies $\mathrm{WA}$; the first occurrences of $x,z$ use the summand $xz$, and any first occurrence of $y$ uses $y$.
Hence $\bu_0+\bs \approx \bu_0$ in $\mathcal U_R$.

Since $\rho$ holds in $\mathcal K_R$, every summand of $\bt$ is below $\bu_0$ there.
The $N,D_2$ conditions imply that such a summand is either the one-letter word $y$, or a word of length at least two starting with $x$ and containing $y$ or both $x,z$.
A summand taking a value outside $I$ must therefore contain $y$ but not $z$.
By \eqref{eq:D-word}, it is equivalent to $xy$.
Thus $\rho$ implies $xy\preceq y+xz$, which is $k$.
The equality already proved gives $\mathcal W=\mathcal K_R$.
Finally $D_2\in\Var(B)$ follows either from Lemma~\ref{lem:band12}
or directly from the subsemiring $B|_{\{2,4\}}$, so the assertion for
$R=B$ follows from \eqref{eq:K-NB}.
\end{proof}

\begin{lem}\label{lem:D-U-HFB}
The variety $\Var(U)$ is hereditarily finitely based.
\end{lem}
\begin{proof}
Let $\mathcal W\subseteq\Var(U)$.
If $N\notin\mathcal W$, use Lemma~\ref{lem:D-omit-N}.
Otherwise the subvariety of multiplicatively idempotent members of $\mathcal W$ contains $\Var(L)$ and is contained in $\Var(L,D_2)$.
Lemma~\ref{lem:band12} gives exactly two possibilities.
If it is $\Var(L,D_2)$, then $N,D_2\in\mathcal W$ and Lemma~\ref{lem:D-U-cover} gives
$\mathcal W=\Var(N,D_2)$ or $\mathcal W=\Var(U)$.
Both are finitely based.
If the idempotent subvariety is $\Var(L)$, its square images satisfy left-zero multiplication.
Since $xy\approx x^2y^2$, this implies $xy\approx x^2$ in $\mathcal W$.
The normal form for $N$ applies to every term; $N$ distinguishes the different forms.
As $N\in\mathcal W$, we obtain $\mathcal W=\Var(N)$, also finitely based.
\end{proof}

\subsection{The variety $\Var(U,M)$}

Let $s$ denote $xyz\approx xzy$. We claim that
\begin{equation}\label{eq:D-regularized}
\mathcal R:=\mathcal J_E+s=\Var(U,M).
\end{equation}
For $\ell(\bq)\ge2$, the right-hand side has the conditions
$h(\bq)\in H(\bu)$, $\mathrm{SA}(\bu,\bq)$, and $c(\bq)\subseteq c(\bu)$.
To derive every such identity on the left, first obtain $\bs \approx \bs^2\preceq \bu$ with $\ell(\bs)\ge2$ and $c(\bs)\subseteq c(\bq)$, as above.
Start with the square of an original summand beginning with $h(\bq)$.
Repeated use of $\mathrm{J2}$ joins to this word the content of every original summand, without changing the first letter.
This yields a word $\bt \approx \bt^2\preceq \bu$ with $h(\bt)=h(\bq)$ and $c(\bt)=c(\bu)$.
The identity $s$ permits ordering the remaining letters so that the first-occurrence list $\mathbf{q}_0$ of $\bq$ comes first, followed by a word $\br$ of the extra letters, if any.
If there are no extra letters, the normal forms in \eqref{eq:D-word} give $\bt \approx \bq$.
Otherwise substitute $x\mapsto \mathbf{q}_0$, $y\mapsto \bs$, $z\mapsto \br$ in $g$ to obtain $\mathbf{q}_0 \bs \preceq \bs^2+\mathbf{q}_0 \br\preceq \bu$.
If $\mathbf{q}_0$ has at least two letters, then $\mathbf{q}_0 \bs \approx \bq$ by repeated-letter deletion.
If $\mathbf{q}_0 =h(\bq)$ has one letter, then $\bs$ is a power of that letter of length at least two, and again $\mathbf{q}_0 \bs \approx \bq$.
For a word of length one on the left, the condition is that some summand is a power of that variable.
Conversely $U,M$ satisfy $s$ and belong to $\mathcal J_E$.
This proves \eqref{eq:D-regularized} and gives a finite equational basis for $\mathcal R$.

\begin{lem}\label{lem:D-R-HFB}
The variety $\mathcal R$ is hereditarily finitely based.
\end{lem}
\begin{proof}
Let $\mathcal W\subseteq\mathcal R$.
The case $N\notin\mathcal W$ follows from Lemma~\ref{lem:D-omit-N}.
Otherwise its multiplicatively idempotent subvariety contains $\Var(L)$ and is contained in $\Var(L,D_2,M)$.
The lattice in Lemma~\ref{lem:band12} gives the following four possibilities.

If it is $\Var(L)$, the preceding argument gives $\mathcal W=\Var(N)$.

If it is $\Var(L,D_2)$, the square images satisfy $xy\preceq x$.
Thus $\mathcal W$ satisfies $xy\preceq x^2$.
Within $\mathcal R$ this identity defines exactly $\Var(U)$.
Indeed, in an identity of $U$ whose left-hand side is a word of length at least two, the only condition for $\mathcal R$ that may fail is $c(\bq)\subseteq c(\bu)$.
Construct $\bt \approx \bt^2\preceq \bu$ with the required first letter and content $c(\bu)$ as above.
The additional identity gives $\bt t\preceq \bt$, and consequently permits adding to $\bt$ each missing letter of $\bq$.
Adding the resulting derived word to $\bu$ makes its total content contain $c(\bq)$, while the other required conditions persist.
The criterion for $\mathcal R$ then proves the identity.
The condition for a word of length one remains unchanged.
Hence $\mathcal W\subseteq\Var(U)$, and Lemma~\ref{lem:D-U-HFB} applies.

If the idempotent subvariety is $\Var(L,M)$, square substitution in $x\preceq xy$ gives $x^2\preceq xy$ in $\mathcal W$.
Within $\mathcal R$ this defines exactly $\Var(N,M)$.
For a word of length at least two on the left, the conditions in $N,M$ are only $h(\bq)\in H(\bu)$ and $c(\bq)\subseteq c(\bu)$.
Construct $\bt$ as above, put the letters of $\bq$ first using $s$, and remove the extra suffix by $x^2\preceq xy$.
If the word on the left contains only one variable, use its square for the left side of this last identity.
For a word $x$ of length one on the left-hand side, the $N$ condition says that $x$ is already a one-letter summand or is the initial letter of a longer summand $x\mathbf{v}_0$.
The latter case follows from $x\preceq x^2\preceq x\mathbf{v}_0$.
This proves the claimed equality and gives a finite equational basis.
Since $N,M\in\mathcal W$, here $\mathcal W=\Var(N,M)$.

Finally, if the idempotent subvariety is $\Var(L,D_2,M)$, then
\[
\Var(N,D_2,M)=\mathcal K_M\subseteq\mathcal W.
\]
Lemma~\ref{lem:D-U-cover} gives $\mathcal W=\mathcal K_M$ or $\mathcal W=\mathcal R$, both finitely based.
These cases exhaust all subvarieties of $\mathcal R$.
\end{proof}

\subsection{Hereditary finite basability of $\mathcal J_E$ and $\mathcal T$}

\begin{lem}\label{lem:D-JE-HFB}
The variety $\mathcal J_E$ is hereditarily finitely based.
\end{lem}
\begin{proof}
Let $\mathcal W\subseteq\mathcal J_E$.
If $N\notin\mathcal W$, apply Lemma~\ref{lem:D-omit-N}.
Otherwise its multiplicatively idempotent subvariety contains $\Var(L)$ and is contained in $\Var(B)$.
If it is all of $\Var(B)$, then $\mathcal K=\Var(N,B)\subseteq\mathcal W$.
The case $R=B$ in Lemma~\ref{lem:D-U-cover} gives $\mathcal W=\mathcal K$ or $\mathcal W=\mathcal J_E$; both have finite bases.
If it is proper, it is contained in $\Var(L,D_2,M)$ by Lemma~\ref{lem:band12}.
Its square images satisfy $s:xyz\approx xzy$.
Replacing variables by squares and using \eqref{eq:D-word} yields $s$ in $\mathcal W$ itself.
Thus $\mathcal W\subseteq\mathcal J_E+s=\mathcal R$, and Lemma~\ref{lem:D-R-HFB} applies.
\end{proof}

\begin{lem}\label{lem:D-T-cover}
Let $j$ denote $xy\preceq y+xyz$. Then
\[
\mathcal T+j=\Var(M_1),\qquad
[\Var(M_1),\mathcal T]
=\{\Var(M_1),\mathcal T\}.
\]
\end{lem}
\begin{proof}
The semiring $M_1$ satisfies $j$.
Consider an identity $\bq \preceq \bu$ of $M_1$, with $\ell(\bq)\ge2$.
If $\mathrm{SA}$ holds, use Lemma~\ref{lem:D-T-criterion}.
Otherwise every summand supported in $c(\bq)$ is a one-letter summand, and none of those letters belongs to $H(\bu)$.
Choose the first such letter $a$ in the first-occurrence order of $\bq$.
There is no summand supported entirely on the preceding letters: such a summand would be another one-letter summand of the same kind occurring earlier.
The last condition of \eqref{eq:M1-criterion} supplies a summand $\mathbf{v}_0 a\br$ of length at least two, with $a$ its first letter outside the preceding set.
The word $\mathbf{v}_0$ is nonempty, since otherwise $a\in H(\bu)$; the word $\br$ is nonempty, since otherwise $\mathbf{v}_0 a$ is a longer summand supported in $c(\bq)$.
Substituting into $j$ gives $\mathbf{v}_0 a\preceq a+\mathbf{v}_0 a\br \preceq \bu$.
The new summand satisfies the first alternative of $\mathrm{SA}$, so Lemma~\ref{lem:D-T-criterion} completes the derivation.
For a word of length one on the left, the condition is that some summand is a power of that variable.
This proves $\mathcal T+j=\Var(M_1)$.

Now let $\rho$ be an identity holding in $M_1$ but failing in $T$.
Every value of a variable in $T$ is represented by a term over $x,y,z$.
The defining congruence of $T$ supplies a context and an orientation for which the substituted identity is $\bs \approx \bt$, all summands of $\bs$ belong to $\mathscr L$, and some summand $\bw$ of $\bt$ does not.
Put $\bv =y+z+xyz$.
Since $\bv +\bs \approx \bv$ holds in $S$, it holds in $\mathcal T$.
Since $\rho$ holds in $M_1$, we have $M_1\models \bw \preceq \bv$.

Apply \eqref{eq:M1-criterion} to this fixed term $\bv$.
The one-letter possibilities are $y,z$.
Every longer word must start with $x$, must contain $y$, and, if it contains $z$, must have its first $y$ before its first $z$.
For the latter requirement, a first $z$ before $y$ has neither a summand supported in its preceding letters nor a summand of length at least two whose first letter outside that preceding set is $z$.
Hence a word satisfying the $M_1$ condition but not lying in $\mathscr L$ contains only $x,y$ and starts with $x$.
It equals $xy$ in $\mathcal D$ by \eqref{eq:D-word}.
Thus $\rho$ implies $xy\preceq y+z+xyz$ over $\mathcal T$.
Replacing $z$ by $xyz$ and using $xyxyz\approx xyz$ gives $j$.
Every proper subvariety of $\mathcal T$ containing $\Var(M_1)$ therefore equals $\Var(M_1)$.
\end{proof}

It follows that $\mathcal T$ is hereditarily finitely based.
Indeed, if $\mathcal W\subseteq\mathcal T$ does not contain $M_1$, the characterization proved earlier gives
$\mathcal W\subseteq\mathcal J\cap\mathcal E=\mathcal J_E$, and Lemma~\ref{lem:D-JE-HFB} applies.
If it contains $M_1$, Lemma~\ref{lem:D-T-cover} says it is one of the two finitely based varieties in that interval.

\subsection{Completion of the proofs for $\mathcal J$ and $\mathcal D$}

Put $\mathcal J_\ell=\mathcal J+s$.
We first prove
\begin{equation}\label{eq:D-Jleft}
\mathcal J_\ell=\Var(N,S_b).
\end{equation}
For a word of length at least two on the left, the conditions on the right are
$h(\bq)\in H(\bu)$ and $\mathrm{Cov}(\bu,\bq)$.
If a summand of length at least two is supported in $c(\bq)$, take its square as an auxiliary word already below $\bu$.
Otherwise condition $\mathrm{Cov}$ makes $h=h(\bq)$ a one-letter summand, and there is also a longer summand $h\mathbf{v}_0$; use $h^2\preceq h+h\mathbf{v}_0$ to obtain such an auxiliary word.
Apply $\mathrm{J2}$ successively with all summands supported in $c(\bq)$.
This gives a word equal to its square, below $\bu$, and having content exactly $c(\bq)$.
If its initial letter is not $h$, use its square and write the part through the first $h$ as $\bs h\br$, with $\bs,\br$ nonempty.
The identity $\mathrm H$, with the original summand $h\mathbf{v}_0$, gives $h\bs h\br \preceq \bs h\br +h\mathbf{v}_0 \preceq \bu$.
Thus its initial letter can be made $h$.
The identity $s$ permutes the remaining letters; \eqref{eq:D-word} then gives $\bq$.
The case of a word of length one again follows from the condition that a summand is a power of that variable.
The reverse inclusion follows because $N,S_b\in\mathcal J$ and satisfy $s$.
This proves \eqref{eq:D-Jleft} and a finite basis for $\mathcal J_\ell$.

The multiplicatively idempotent subvariety of $\mathcal J_\ell$ is
$\Var(L,S_b)$.
Take any $\mathcal W\subseteq\mathcal J_\ell$.
If it does not contain $N$, it is finitely based by Lemma~\ref{lem:D-omit-N}.
Otherwise its idempotent subvariety contains $\Var(L)$.
If that subvariety is all of $\Var(L,S_b)$, then $N,S_b\in\mathcal W$, so $\mathcal W=\mathcal J_\ell$.
If it is proper, it cannot contain $S_b$, since it already contains $L$.
Consequently $S_b\notin\mathcal W$, and the earlier characterization gives $\mathcal W\subseteq\mathcal E$.
Hence
\[
\mathcal W\subseteq\mathcal J_\ell\cap\mathcal E
=\mathcal J_E+s=\mathcal R.
\]
By Lemma~\ref{lem:D-R-HFB}, every such $\mathcal W$ is finitely based.
Thus $\mathcal J_\ell$ is hereditarily finitely based.

Next take $\mathcal W\subseteq\mathcal J$ containing $N$.
If its multiplicatively idempotent subvariety is $\Var(C_4)$, then $N,C_4\in\mathcal W$ and $\mathcal W=\mathcal J$.
Otherwise Lemma~\ref{lem:band12} places it in one of the two maximal proper subvarieties
$\Var(B)$ or $\Var(L,S_b)$.
In the former case $S_b\notin\mathcal W$, so $\mathcal W\subseteq\mathcal J_E$.
In the latter case its square images satisfy $s$, and \eqref{eq:D-word} again implies $s$ in $\mathcal W$, so $\mathcal W\subseteq\mathcal J_\ell$.
Both containing varieties are hereditarily finitely based.
Combining these cases with Lemma~\ref{lem:D-omit-N} proves that $\mathcal J$ is hereditarily finitely based.

Finally let $\mathcal W\subseteq\mathcal D$ be arbitrary.
If $M_1\notin\mathcal W$, then $\mathcal W\subseteq\mathcal J$.
If $M_1\in\mathcal W$, its multiplicatively idempotent subvariety contains $\Var(B)$.
The interval $[\Var(B),\Var(C_4)]$ in Lemma~\ref{lem:band12} has two elements.
If the idempotent subvariety is $\Var(B)$, then $S_b\notin\mathcal W$ and $\mathcal W\subseteq\mathcal T$.
If it is $\Var(C_4)$, then $M_1,C_4\in\mathcal W$ and $\mathcal W=\mathcal D$.
The varieties $\mathcal J,\mathcal T$ are hereditarily finitely based and $\mathcal D$ has a finite basis.
Every possibility therefore gives a finite basis for $\mathcal W$.
We conclude that
\begin{equation}\label{eq:D-HFB-conclusion}
\boxed{\mathcal D=\Var(M_1,C_4)\text{ is hereditarily finitely based}.}
\end{equation}
The argument applies to arbitrary subvarieties.
The preceding results establish hereditary finite basability successively for $\Var(U)$, $\mathcal R$, $\mathcal J_E$, $\mathcal T$, $\mathcal J_\ell$, $\mathcal J$, and $\mathcal D$.
Hereditary finite basability of $\mathcal C$ and the limit-variety assertion will be proved in the following sections.
\section{Every subvariety of $\mathcal J_0$ is finitely based}\label{sec:j0}

In this section, we prove that $\mathcal J_0=\Var(M_0,C_4)$ is
hereditarily finitely based. The proof applies to arbitrary subvarieties,
including those without a finite generating algebra, and to identities on
arbitrary finite sets of variables. We use the description of the finite
lattice of multiplicatively idempotent subvarieties in Lemma~\ref{lem:band12}.

Lemmas~\ref{lem:wordnormal}, \ref{lem:helpers}, and \ref{lem:C-omit-N3},
together with Theorem~\ref{thm:C-fb}, give a finite equational basis for $\mathcal J_0$
and the following identities:
\begin{align}
 x^2y&\approx xy,&xyx&\approx xy^2,&x&\preceq x^2,\label{eq:j0-word}\\
 (x+y)^2&\approx x^2+y^2,&(xy)^2&\approx x^2y^2,&x^4&\approx x^2,\label{eq:j0-square}\\
 x^2+y&\approx x^2+y^2,\label{eq:j0-close}\\
 xyz&\preceq xz+yz,\label{eq:j0-ST}\\
 yxyz&\preceq xyz+yt,&xyzt&\preceq xzyt+yu.\label{eq:j0-HI}
\end{align}
The notation $\bw \preceq \bu $ abbreviates the identity $\bu +\bw \approx \bu $.
We also use the reordering assertion of Lemma~\ref{lem:carrier-reorder}.
In its applications below, the last variable is fixed, and the variables preceding
it are put in a prescribed order by their first occurrences. The required
original summands are explicitly supplied whenever that lemma is applied.

\subsection{Notation and subvarieties not containing $M_0$}

For a nonempty word $\bw $, let $c(\bw )$ be its set of variables, $h(\bw )$ its first
variable, $t(\bw )$ its last variable, and $e(\bw )$ the set of variables occurring
before its last position. Thus $e(\bw )=\varnothing$ for a one-letter word.
For a nonempty term $\bu =\sum_{\bw \in\Omega}\bw $, put
\[
 c(\bu )=\bigcup_{\bw \in\Omega}c(\bw ),\qquad
 E(\bu )=\bigcup_{\bw \in\Omega}e(\bw ),\qquad
 T(\bu )=\{t(\bw ):\bw \in\Omega\}.
\]
For use within this section only, write
\[
 \mathrm{(O)}_{\bu}:\quad E(\bu )\cap T(\bu )=\varnothing.
\]
This notation will be used only in this section. Also put
\[
 E_t(\bu )=\bigcup_{\substack{\bw \in\Omega\\t(\bw )=t}}e(\bw ),
 \qquad H(\bu )=\{h(\bw ):\bw \in\Omega,\ \ell(\bw)\ge2\}.
\]
If $\mathrm{(O)}_{\bu}$ holds, the variables occurring before the last position
of a summand and those occurring in a last position form disjoint sets.
In particular, $H(\bu )\subseteq E(\bu )$.

By \eqref{eq:j0-word}, repeated occurrences in a proper prefix can be deleted
after their first occurrence. Whenever a word with prescribed content is
chosen, we choose an order of its variables and use each variable once.
Whenever the first variable or the order matters, it will be specified. Every word substituted for
a variable is nonempty. An empty left or right context means that the
corresponding multiplication is omitted; it does not mean that an empty word
is substituted for a variable.

\begin{lem}\label{lem:j0-M0}
The identity $\bq \preceq \bu $ holds in $M_0$ if and only if
\[
 c(\bq )\subseteq c(\bu )
 \quad\text{and}\quad
 \left[
 E(\bu )\cap T(\bu )\ne\varnothing
 \ \text{or}\
 \bigl(e(\bq )\subseteq E(\bu ),\ t(\bq )\in T(\bu )\bigr)
 \right].
\]
Moreover, $\mathcal J_0\models \bu \approx \bu ^2$ if and only if
$E(\bu )\cap T(\bu )\ne\varnothing$.
\end{lem}
\begin{proof}
Consider the three-element subsemiring $D_3=S|_{\{0,1,3\}}$.
A word has value $3$ precisely when all its variables have value $3$.
It has value $1$ precisely when every position except the last has value $3$
and its last position has value $1$. All other evaluations have value $0$.
The proper principal down-sets of its additive order are $\{1\}$ and $\{3\}$.
For $\{3\}$, requiring $\bq $ to have value $3$ whenever all summands of
$\bu $ do is equivalent to $c(\bq )\subseteq c(\bu )$.

All summands of $\bu $ can have value $1$ if and only if
$E(\bu )\cap T(\bu )=\varnothing$. In that case the relevant assignment gives value
$3$ to $E(\bu )$ and value $1$ to $T(\bu )$; the word $\bq $ has value $1$ exactly
when $e(\bq )\subseteq E(\bu )$ and $t(\bq )\in T(\bu )$, once the content condition
has excluded variables outside $c(\bu )$. This gives the criterion for $D_3$.

The additional proper principal down-set of $M_0$ is $\{2,3\}$.
A word has its value in this set precisely when every variable occurring
in it has value in $\{2,3\}$. Its membership implication is therefore again
$c(\bq )\subseteq c(\bu )$. The down-sets $\{1\}$ and $\{3\}$ give the same
conditions as before. Since principal down-sets separate elements,
the criterion follows for $M_0$ as well. In particular,
$\Var(D_3)=\Var(M_0)$.

If $E(\bu )\cap T(\bu )\ne\varnothing$, the criterion gives $\bu \approx \bu ^2$ in $M_0$.
The same identity holds in $C_4$ because every element of $C_4$ is
multiplicatively idempotent. Hence it holds in $\mathcal J_0$.
If $\mathrm{(O)}_{\bu}$ holds, assign $3$ to $E(\bu )$ and $1$ to $T(\bu )$ in $D_3$.
Then $\bu =1$ and $\bu ^2=0$, so the identity fails.
\end{proof}

\begin{lem}\label{lem:j0-seed}
Suppose that a subvariety of $\mathcal J_0$ satisfies $\mathbf{z} \approx \mathbf{z} ^2$ and $\mathbf{z} \preceq \bu $.
Then it satisfies $\bu \approx \bu ^2$.
\end{lem}
\begin{proof}
Monotonicity of the squaring map gives $\mathbf{z} \approx \mathbf{z} ^2\preceq \bu ^2$.
Substituting $x=\mathbf{z} $ and $y=\bu $ into \eqref{eq:j0-close} yields
$\bu \approx \mathbf{z} +\bu \approx \mathbf{z} +\bu ^2\approx \bu ^2$.
\end{proof}

\begin{lem}\label{lem:j0-omit}
If $\mathcal W\subseteq\mathcal J_0$ and $M_0\notin\mathcal W$, then
$\mathcal W$ is finitely based.
\end{lem}
\begin{proof}
Within $\mathcal J_0$, we first show that $M_0\notin\mathcal W$
is equivalent to
\begin{equation}\label{eq:j0-product-square}
 xy\approx(xy)^2.
\end{equation}
Since $\Var(D_3)=\Var(M_0)$, an identity of
$\mathcal W$ fails in $D_3$. For any two different elements of $D_3$, either
the identity map or right multiplication by $1$ makes precisely one image
equal to $1$. Right multiplication is implemented by an additional
variable with prescribed value, not by adding a constant to the signature.
Representing the values $3,1,0$ by $x,y,y^2$, respectively, gives a
two-variable identity $\mathbf{a} \approx \mathbf{b} $ of $\mathcal W$ with
$\mathbf{a} (3,1)=1$ and $\mathbf{b} (3,1)\ne1$.
Every summand of $\mathbf{a} $ is $x^ky$ for some $k\ge0$, so $x\mathbf{a} \approx xy$.
At least one summand of $\mathbf{b} $ is not of this form. After multiplication
on the left by $x$, that summand either contains only $x$, or contains a
nonfinal $y$. Its last variable has already occurred earlier:
if it is $x$, use the prefixed $x$; if it is $y$, use that nonfinal $y$.
Thus $E(x\mathbf{b} )\cap T(x\mathbf{b} )\ne\varnothing$.
Lemma~\ref{lem:j0-M0} gives $x\mathbf{b} \approx(x\mathbf{b} )^2$, and $x\mathbf{a} \approx x\mathbf{b} $ proves
\eqref{eq:j0-product-square}. Conversely, \eqref{eq:j0-product-square}
fails in $M_0$ at $(x,y)=(3,1)$.

Let $N_2=S|_{\{0,1\}}$. All its products equal the greatest element $0$
of its additive order. We claim that
\[
 \mathcal J_0\cap\operatorname{Mod}\bigl(xy\approx(xy)^2\bigr)
 =\Var(C_4,N_2).
\]
Let $A$ belong to the left-hand side. The map $r(a)=a^2$ is an
idempotent endomorphism by \eqref{eq:j0-square}. Write $I=r(A)$.
Then $I\in\Var(C_4)$: every identity of $C_4$ remains valid
after replacing its variables by squares in $C_4$ and $M_0$, since
the square image $\{0,2,3\}$ of $M_0$ is a subsemiring of $C_4$.
Equation~\eqref{eq:j0-close} gives $I+A\subseteq I$, and
\eqref{eq:j0-product-square} gives $AA\subseteq I$.
The equivalence relation $\rho$ having $I$ as one class and each element
outside $I$ as a singleton is therefore a semiring congruence.

All products in $A/\rho$ equal its greatest additive element $I/\rho$.
Every semiring with this property belongs to $\Var(N_2)$.
Indeed, for every proper principal down-set of its additive semilattice,
the map taking value $1$ on the down-set and $0$ outside is a homomorphism
to $N_2$. It preserves joins because a principal down-set is join-closed;
it preserves products because the greatest element lies outside
every proper principal down-set. These maps separate different
elements. The assertion also holds for the one-element semiring.
Finally,
\[
 A\longrightarrow I\times(A/\rho),\qquad
 a\longmapsto(a^2,[a]_\rho)
\]
is an injective homomorphism: the second coordinate separates all pairs
except pairs in $I$, and $r$ restricts to the identity on $I$.
This proves one inclusion. The converse follows by checking the identity
in $C_4$ and $N_2$.
Theorem~\ref{thm:RenZeng} applies to $\Var(C_4,N_2)$:
it is contained in the join of the variety of multiplicatively idempotent
ai-semirings and the variety satisfying $xy\approx zt$.
Every subvariety of that join is finitely based, including $\mathcal W$.
\end{proof}

\subsection{Multiplicative idempotence}

Use the notation of Lemma~\ref{lem:band12}: $\mathcal L$ is generated by
the two-element ai-semiring with left-zero multiplication, $\mathcal M$
by the two-element algebra whose two operations are both join, and
$\mathcal D_{\rm lat}$ by the two-element distributive lattice. Put
\[
 \mathcal K=\mathcal D_{\rm lat}\vee\mathcal M,\qquad
 \mathcal Q=\mathcal L\vee\mathcal K,\qquad
 \mathcal B=\Var(S|_{\{0,2,4\}}),\qquad
 \mathcal S_b=\Var(S|_{\{2,3,4\}}).
\]
For $\mathcal U\subseteq\Var(C_4)$, define
\[
 \mathcal R_{\mathcal U}
 =\{A\in\mathcal J_0:r(A)\in\mathcal U\}.
\]
The variety $\mathcal R_{\mathcal U}$ is finitely based: add to a finite
equational basis of $\mathcal J_0$ the
identities obtained from a finite equational basis of $\mathcal U$ by replacing
every variable by its square. For any $\mathcal W\subseteq\mathcal J_0$, put
\begin{equation}\label{eq:j0-band-part}
 \mathcal U_{\mathcal W}
 =\mathcal W\cap\operatorname{Mod}(x^2\approx x).
\end{equation}
This is exactly the variety generated by the square images of members
of $\mathcal W$. Each image is a subalgebra of its original algebra,
and each multiplicatively idempotent member is its own square image.
If $M_0\in\mathcal W$, then $\mathcal M\subseteq\mathcal U_{\mathcal W}$.

We shall use the following conditions:
\begin{align*}
 \mathrm A(\bu ,\bq )&:\quad
   \text{there is }\bp \in\Omega\text{ with }c(\bp )\subseteq c(\bq ),\\
 \mathrm C(\bu ,\bq )&:\quad
   c(\bq )=\bigcup\{c(\bp ):\bp \in\Omega,\ c(\bp )\subseteq c(\bq )\}.
\end{align*}
Let $\operatorname{Ex}(\bq ,\bu )$ mean the following: if $b$ is a variable
at its first occurrence in $\bq $ and $D$ is the set of variables
preceding it, some $\bp \in\Omega$ has $b$ as its first variable not in $D$.
The following are the identity criteria:
\[
\begin{array}{c|l}
 \text{variety}&\text{conditions for }\bq \preceq \bu \\\hline
 \mathcal M&c(\bq )\subseteq c(\bu )\\
 \mathcal L&h(\bq )=h(\bp )\text{ for some }\bp \in\Omega\\
 \mathcal K&c(\bq )\subseteq c(\bu ),\ \mathrm A(\bu ,\bq )\\
 \mathcal Q&c(\bq )\subseteq c(\bu ),\ \mathrm A(\bu ,\bq ),\
                  h(\bq )=h(\bp )\text{ for some }\bp \in\Omega\\
 \mathcal B&\mathrm A(\bu ,\bq ),\ \operatorname{Ex}(\bq ,\bu )\\
 \mathcal S_b&\mathrm C(\bu ,\bq )\\
 \mathcal L\vee\mathcal S_b&\mathrm C(\bu ,\bq ),\
                  h(\bq )=h(\bp )\text{ for some }\bp \in\Omega.
\end{array}
\]
For $\mathcal M$ and $\mathcal L$, evaluations depend only on content
and on first variables, respectively. In the two-element distributive
lattice, $\bq \preceq \bu $ is equivalent to $\mathrm A(\bu ,\bq )$.
This gives the criteria for $\mathcal K$ and $\mathcal Q$ as well.
In $B=S|_{\{0,2,4\}}$, the element $2$ is an identity and $0,4$ are
left zeros. A word has the value of its first variable not assigned $2$,
or value $2$ if there is none. The proper principal down-sets
$\{4\},\{2,4\}$ give the displayed conditions as follows.
If $\operatorname{Ex}(\bq ,\bu )$ fails at $b$, assign $2$ to the variables
preceding its first occurrence in $\bq $, $0$ to $b$, and $4$ to the other
variables. The word $\bq $ has value $0$, whereas every summand of $\bu $
has value in $\{2,4\}$. If $\mathrm A(\bu ,\bq )$ fails, assign $2$ to
$c(\bq )$ and $4$ to the other variables; then $\bq =2$ and $\bu =4$.
Conversely, if $\bq $ has value $0$, apply $\operatorname{Ex}$ at its
first variable not assigned $2$: it supplies a summand of $\bu $
with value $0$. If $\bq $ has value $2$, all variables of $\bq $ have
value $2$, and $\mathrm A$ supplies a summand with value $2$.
If $\bq =4$, then $\bq \leq_B \bu $ since $4$ is the least element. Thus the two conditions
are sufficient as well as necessary.
In $S|_{\{2,3,4\}}$, the element $4$ is a multiplicative zero; in its
absence a word has value $2$ precisely when some variable does, and
otherwise has value $3$. Its additive order is $4<3<2$.
If $\mathrm C(\bu ,\bq )$ fails at $b\in c(\bq )$, assign $2$ to $b$,
$3$ to $c(\bq )\setminus\{b\}$, and $4$ to the other variables.
Then $\bq =2$, but every summand with content in $c(\bq )$ has value $3$
and every other summand has value $4$, so $\bu <2$.
Conversely, assume $\mathrm C(\bu ,\bq )$.
When $\bq =2$, choose a variable $b$ of $\bq $ with value $2$; the condition
supplies a summand containing $b$ whose content is contained in $c(\bq )$,
and this summand also has value $2$.
When $\bq =3$, any summand with content contained in $c(\bq )$ has value $3$;
at least one exists because $c(\bq )$ is nonempty.
When $\bq =4$, the inequality holds since $4$ is the least element.
This proves the criterion $\mathrm C(\bu ,\bq )$.
Identities of a join are precisely the common identities of its factors.

\begin{lem}\label{lem:j0-closed-band}
Suppose that $\mathcal W\subseteq\mathcal R_{\mathcal U}$,
$M_0\in\mathcal W$, and $\mathcal U_{\mathcal W}=\mathcal U$.
If $E(\bu )\cap T(\bu )\ne\varnothing$ and $\mathcal W\models \bq \preceq \bu $,
then $\mathcal R_{\mathcal U}\models \bq \preceq \bu $.
If $\mathrm{(O)}_{\bu}$ holds, every such identity satisfies
$e(\bq )\subseteq E(\bu )$ and $t(\bq )\in T(\bu )$.
In the latter case, writing $\bq =\bs t$ with $\bs $ nonempty, every summand
$\bp $ with $c(\bp )\subseteq c(\bq )$ equals $t$ or has the form $\br t$,
where $\br $ is nonempty and $c(\br )\subseteq c(\bs )$.
\end{lem}
\begin{proof}
Since $\mathcal U\subseteq\mathcal W$, the identity holds in $\mathcal U$.
Applying it in a square image gives $\bq ^2\preceq \bu ^2$ throughout
$\mathcal R_{\mathcal U}$.
If $E(\bu )\cap T(\bu )\ne\varnothing$, use Lemma~\ref{lem:j0-M0}
and $\bq \preceq \bq ^2$ to conclude $\bq \preceq \bu $.
Under $\mathrm{(O)}_{\bu}$, the two inclusion conditions follow from $M_0$.
Since $c(\bs )\subseteq E(\bu )$ and $E(\bu )\cap T(\bu )=\varnothing$, the only
possible last variable of a summand with content in $c(\bs t)$ is $t$.
\end{proof}

For $\mathrm{(O)}_{\bu}$ and $\bq =\bs t$ with $\bs $ nonempty, write
\[
 \mathrm A_+(\bu ,\bs t):\quad
 \text{there is a summand }\br t\text{ of }\bu
 \text{ with }\br \ne\varepsilon\text{ and }c(\br )\subseteq c(\bs ).
\]
Under the occurrence conditions of Lemma~\ref{lem:j0-M0},
$\mathrm A(\bu ,\bs t)$ is equivalent to $\mathrm A_+(\bu ,\bs t)$ or to $t$
being a one-letter summand of $\bu $.

\subsection{Constructing auxiliary words}

Fix the identities
\begin{align}
 \lambda &: xz\preceq z+xyz,&
 \tau &: xyz\preceq xz+xyw,\label{eq:j0-lambda-tau}\\
 \nu &: xayt\preceq xat+yaw,&
 \iota &: xyz\preceq xz+yw.\label{eq:j0-nu-iota}
\end{align}
Write $\mathcal R[\Sigma]=\mathcal R\cap\operatorname{Mod}(\Sigma)$.
All four identities hold in $M_0$ and $\mathcal B$, by the above criteria.
They therefore hold in every multiplicatively idempotent subvariety of
$\mathcal B$ used below.
Substituting $x\mapsto xa$, $z\mapsto t$, $w\mapsto aw$ into $\iota$
gives $\nu$. Substituting $a\mapsto x$, $y\mapsto xy$ into $\nu$
and deleting repeated nonfinal variables gives $\tau$.

\begin{lem}\label{lem:j0-carriers}
Suppose that $\mathrm{(O)}_{\bu}$, $\bq =\bs t$, $\bs \ne\varepsilon$, and
$c(\bs )\subseteq E(\bu )$. Each of the following sets of conditions suffices
for the indicated variety to satisfy $\bs t\preceq \bu $:
\begin{enumerate}
\item in $\mathcal R_{\mathcal K}$:
      $c(\bs )\subseteq E_t(\bu )$ and $\mathrm A_+(\bu ,\bs t)$;
\item in $\mathcal R_{\mathcal Q}$:
      the preceding conditions and $h(\bs )\in H(\bu )$;
\item in $\mathcal R_{\mathcal B}$:
      the conditions in item~1 and $\operatorname{Ex}(\bs t,\bu )$;
\item in $\mathcal R_{\mathcal L\vee\mathcal M}$:
      $c(\bs )\subseteq E_t(\bu )$ and $h(\bs )\in H(\bu )$.
\end{enumerate}
In items~1--3, adding $\lambda$ allows $\mathrm A_+$ to be replaced by $\mathrm A$.
\end{lem}
\begin{proof}
Repeated applications of \eqref{eq:j0-ST} to the summands of length
at least two ending in $t$ produce $\mathbf{w}_0 t\preceq \bu $ with $c(\mathbf{w}_0 )=E_t(\bu )$.
Their prefixes are nonempty. Such a summand exists in items~1--3 by
$\mathrm A_+(\bu ,\bs t)$; in all four items its existence also follows from
$\varnothing\ne c(\bs )\subseteq E_t(\bu )$.
The latter observation applies as well when $\mathrm A_+$ is replaced
by $\mathrm A$ in the presence of $\lambda$.
In $\mathcal R_{\mathcal K}$, square substitution in $xy\approx yx$, followed
by right multiplication by $z$ and \eqref{eq:j0-word}, gives $xyz\approx yxz$.
Hence prefix variables may be put in any order.
In $\mathcal R_{\mathcal Q}$, the same procedure applied to $xyz\approx xzy$
gives $xyzt\approx xzyt$, permitting permutation after the first prefix variable.
To put $h$ first, write the auxiliary word as $\mathbf{a} h\mathbf{r}_0 $ with $\mathbf{a} ,\mathbf{r}_0 $ nonempty;
$\mathbf{r}_0 $ includes the final $t$. Choose an original summand $h\br $.
Equation \eqref{eq:j0-HI} gives $h\mathbf{a} h\mathbf{r}_0 \preceq \mathbf{a} h\mathbf{r}_0 +h\br \preceq \bu $.
In $\mathcal R_{\mathcal B}$, apply Lemma~\ref{lem:carrier-reorder}
using $\operatorname{Ex}(\bs t,\bu )$.
For each variable in $\bs $, its specified occurrence in the corresponding
original summand is nonfinal, since it lies in $E(\bu )$ and not in $T(\bu )$.
Thus the hypotheses of that reordering lemma hold.

In items~1--3 we obtain $\bs \mathbf{z} t\preceq \bu $, where $\mathbf{z} $ consists of unnecessary
prefix variables and may be absent. If it is absent, the required identity
has been obtained.
Otherwise, the identity $xy\preceq y+xz$ of $\mathcal B$ and its subvarieties,
after square substitution and right multiplication by $t$, gives
\begin{equation}\label{eq:j0-trim}
 xyt\preceq yt+xzt.
\end{equation}
Take the summand $\br t$ given by $\mathrm A_+$ and substitute
$x=\bs $, $y=\br $, $z=\mathbf{z} $. Then $\bs \br t\preceq \br t+\bs \mathbf{z} t\preceq \bu $.
Since $c(\br )\subseteq c(\bs )$, repeated-prefix deletion gives $\bs \br t\approx \bs t$.
If only $\mathrm A$ is known and $\mathrm A_+$ fails, $t$ is an original
summand, and $\lambda$ gives $\bs t\preceq t+\bs \mathbf{z} t\preceq \bu $.

For item~4, square substitution in $x\preceq xy$ of
$\mathcal L\vee\mathcal M$, followed by right multiplication, gives
\begin{equation}\label{eq:j0-PD}
 xz\preceq xyz.
\end{equation}
After putting the first variable in place and arranging the rest,
this gives $\bs t\preceq \bs \mathbf{z} t\preceq \bu $ without condition $\mathrm A$.
\end{proof}

Under $\mathrm{(O)}_{\bu}$, define a bipartite graph $G(\bu )$ with parts
$E(\bu )$ and $T(\bu )$, with $b$ adjacent to $s$ exactly when $b\in E_s(\bu )$.
Define $G_H(\bu )$ by deleting the vertices of $E(\bu )\setminus H(\bu )$ and
their edges, retaining all vertices of $T(\bu )$. Put
\[
 O_t=\{b\in E(\bu ):b\text{ lies in the component of }t\text{ in }G(\bu )\},
\]
and
\[
 A_t=\bigcup\{E_s(\bu ):s\text{ lies in the component of }t\text{ in }G_H(\bu )\}.
\]
The set $A_t$ may thus include variables outside $H(\bu )$.
An isolated $t$ has $O_t=A_t=\varnothing$.
We retain this notation throughout the remainder of this section.

\begin{lem}\label{lem:j0-component-carriers}
In Lemma~\ref{lem:j0-carriers}, adding $\tau$ permits $E_t(\bu )$ to be
replaced by $A_t$; for $\mathcal R_{\mathcal K}$ it may instead be
replaced by $O_t$. Adding $\nu$ permits replacement by $O_t$ in all
four cases. Every other condition remains unchanged.
\end{lem}
\begin{proof}
For every nonisolated last variable $s$, equation \eqref{eq:j0-ST}
produces a word ending in $s$ whose prefix contains exactly $E_s(\bu )$.
If two such words share a variable $h\in H(\bu )$, equation
\eqref{eq:j0-HI} puts $h$ first in each; write them $h\mathbf{a} t$ and $h\mathbf{b} u$.
If $\mathbf{b} $ contributes no new variable, nothing needs to be added.
Otherwise substitute $x=h$, $y=\mathbf{b} $, $z=\mathbf{a} t$, $w=u$ into $\tau$,
obtaining $h\mathbf{b} \mathbf{a} t\preceq \bu $. The substituted word $\mathbf{a} t$ is nonempty even
when $\mathbf{a} $ is absent. Following paths in $G_H(\bu )$ gives prefix set $A_t$.
In $\mathcal R_{\mathcal K}$ any shared prefix variable can be put first,
so the same argument follows all paths in $G(\bu )$ and gives $O_t$.

For $\nu$, let $\mathbf{a} t\preceq \bu $, $\mathbf{b} u\preceq \bu $, and $a\in c(\mathbf{a} )\cap c(\mathbf{b} )$.
Write $\mathbf{a} =\mathbf{x} a\mathbf{z} $, choosing the first $a$. If $\mathbf{x} $ is nonempty, substitute
\[
 x=\mathbf{x} ,\quad a=a,\quad y=\mathbf{b} ,\quad t=\mathbf{z} t,\quad w=u
\]
into $\nu$. Its two upper-side summands are $\mathbf{a} t$ and $\mathbf{b} au$, and
$\mathbf{b} au\approx \mathbf{b} u$ by \eqref{eq:j0-word}. Its lower word is $\mathbf{x} a\mathbf{b} \mathbf{z} t$,
with prefix set $c(\mathbf{a} )\cup c(\mathbf{b} )$.
If $\mathbf{x} $ is absent, use $x=a$ and delete the extra $a$ by $x^2y\approx xy$.
Following paths in $G(\bu )$ gives $O_t$.
Adjoin the words so obtained to $\bu $; this does not change $\bu $ under
the current identities. The original summands giving $\mathrm A$,
$\mathrm A_+$, the first-variable condition, and $\operatorname{Ex}$
remain present. Apply Lemma~\ref{lem:j0-carriers}.
\end{proof}

\subsection{A uniform reduction to at most 260 variables}

Set $\mathcal U=\mathcal R_{\mathcal Q}[\iota]$.
We first prove that this variety is hereditarily finitely based.

\begin{lem}\label{lem:j0-global-sufficient}
Suppose that $\mathrm{(O)}_{\bu}$, $\bq =\bs t$, and $\bs $ is nonempty.
If
\[
 c(\bs )\subseteq E(\bu ),\qquad h(\bs )\in H(\bu ),\qquad
 \mathrm A_+(\bu ,\bs t),
\]
then $\mathcal U\models \bs t\preceq \bu $.
\end{lem}
\begin{proof}
Begin with the summand $\br t$ given by $\mathrm A_+$.
If a variable $b\in c(\bs )$ has not yet occurred in the prefix, choose an
original summand $\mathbf{b} w$ with $b\in c(\mathbf{b} )$; in particular $\mathbf{b} $ is nonempty.
In $\iota$, substitute the current nonempty prefix for $x$, $\mathbf{b} $ for $y$,
and the indicated last variables for $z,w$.
The resulting identity inserts $\mathbf{b} $ before $t$ in a word whose inequality
with $\bu $ has already been established.
Repetition gives a word ending in $t$ whose prefix contains $c(\bs )$.
Adjoin it to $\bu $ and apply item~2 of Lemma~\ref{lem:j0-carriers},
using the original summands still giving $\mathrm A_+$ and $h(\bs )\in H(\bu )$.
\end{proof}

\begin{thm}\label{thm:j0-global-hfb}
The variety $\mathcal U$ is hereditarily finitely based.
\end{thm}
\begin{proof}
The variety $\mathcal U$ is finitely based by definition and locally
finite because $\mathcal U\subseteq\mathcal J_0\subseteq\Var(S)$ and $S$
is finite. We first reduce the number of variables in an arbitrary additional
identity, whether or not it holds in $\mathcal U$.

Suppose that $\bq \preceq \bu $ satisfies $\mathrm{(O)}_{\bu +\bq }$.
Keep the last variable $t$ of $\bq $ unchanged.
If $\bq =\bs t$ has length at least two, also keep $h=h(\bs )$ unchanged and put
$D=c(\bs )$. Partition the remaining nonfinal variables into
\[
 X=(E(\bu )\cap D)\setminus\{h\},\qquad
 Y=E(\bu )\setminus D,\qquad
 Z=(D\setminus c(\bu ))\setminus\{h\}.
\]
Identify each nonempty group with one new variable.
If $h\notin E(\bu )$, it remains separate nevertheless.
The original set $E(\bu )$ is divided into at most three nonempty groups:
$\{h\}$, $X$, and $Y$, with $\{h\}$ omitted when $h\notin E(\bu )$.
For a one-letter $\bq =t$, use only $Y=E(\bu )$; neither $h$ nor $Z$ is needed.

For each $s\in T(\bu )\setminus\{t\}$ record whether $s$ is itself a
summand and the family of subsets of these three groups occurring in
prefixes of summands ending in $s$.
There are at most seven nonempty subsets of three groups.
If a summand of length at least two ends in $s$, there are at most
$2(2^7-1)$ possibilities for these data.
If none does, $s$ is a one-letter summand, giving one further possibility.
Hence there are at most
\[
 2(2^7-1)+1=255
\]
types of last variables. Identify last variables having the same type,
retaining $t$ separately. All variables representing different groups
or types are distinct. This gives a substitution consequence
$\bq _0\preceq \bu _0$ on at most
\[
 1+1+1+1+1+255=260
\]
variables: those for $h,X,Y,Z,t$, when present, and the last-variable
types. The disjointness condition is preserved.

For a reverse substitution, replace each prefix-group variable by a product
of all original members of that group.
If $X$ or $Y$ contains a member of $H(\bu )$, choose such a member first.
Keep $h$ unchanged and order $Z$ arbitrarily.
Replace each last-variable type by the sum of its original members,
keeping $t$ unchanged. Every substituted sum and product is nonempty.
The image of $\bq _0$ equals $\bq $ in $\mathcal U$: its first variable and
prefix set are unchanged, the other prefix variables may be permuted
by $xyzt\approx xzyt$, and repetitions may be deleted.

Let $\bp $ be a summand of length at least two of the resulting term $\bu '$.
If its last variable is $s$, the type data ensure that an original
summand ending in $s$ uses the same collection of prefix groups.
Its prefix is contained in $e(\bp )$, so it gives $\mathrm A_+(\bu ,\bp )$.
Also $e(\bp )\subseteq E(\bu )$.
The first group of $\bp $ contained the first variable of some original
summand. The chosen order of each group product therefore gives
$h(\bp )\in H(\bu )$. Lemma~\ref{lem:j0-global-sufficient} yields $\bp \preceq \bu $.
Each one-letter summand of $\bu '$ is already a summand of $\bu $.
Thus $\mathcal U\models \bu '\preceq \bu $, and the reverse substitution in
$\bq _0\preceq \bu _0$ proves $\bq \preceq \bu $. This argument requires only $\bu '\preceq \bu $.
This proves equivalence, over $\mathcal U$, with an identity on at most
260 variables.

Next consider an additional identity $\bu \approx \bu ^2$.
If $E(\bu )\cap T(\bu )\ne\varnothing$, it holds in $\mathcal U$ by
Lemma~\ref{lem:j0-M0}.
Otherwise identify all members of $E(\bu )$ with one variable.
Partition $T(\bu )$ into at most three classes, according as a variable
occurs only as a one-letter summand, only as the last variable of a
longer summand, or in both ways.
Identify each nonempty class with one variable. This gives
$\bu _0\approx \bu _0^2$ on at most four variables.
For a reverse substitution, replace the prefix variable by the product
of $E(\bu )$, choosing a member of $H(\bu )$ first, and replace every
last-variable class by the sum of its members.
If $E(\bu )$ is nonempty, so is $H(\bu )$.
Each resulting longer summand has prefix set $E(\bu )$, first variable in
$H(\bu )$, and an original summand ending in the same variable.
Lemma~\ref{lem:j0-global-sufficient} gives $\bu '\preceq \bu $.
Choose any summand $\bw $ of the nonempty term $\bu '$.
The reverse substitution gives
$\bw ^2\preceq \bu '^2\approx \bu '\preceq \bu $, and \eqref{eq:j0-square} gives
$\bw ^2\approx(\bw ^2)^2$. Lemma~\ref{lem:j0-seed} implies $\bu \approx \bu ^2$.
If $E(\bu )=\varnothing$, the substitutions by sums recover $\bu $ exactly,
without substituting any prefix product.

Finally let $\mathcal W\subseteq\mathcal U$ be arbitrary.
Its multiplicatively idempotent subvariety
$\mathcal W\cap\operatorname{Mod}(x^2\approx x)$ has a finite equational basis $\Gamma$
by Theorem~\ref{thm:band-literature}.
Let $\Delta$ consist of a finite equational basis of $\mathcal U$, the
square-variable instances of $\Gamma$, and a finite representative
set of all identities of $\mathcal W$ on a fixed 260-variable alphabet.
The last set is finite because the free algebra of $\mathcal U$
on that alphabet is finite. Choose one term for each of its finitely
many term classes and include every equality between representatives
that holds in $\mathcal W$.
All identities in $\Delta$ hold in $\mathcal W$.

Suppose that $\mathcal W\models \bq \preceq \bu $.
The square-variable instances of $\Gamma$ give
$\Delta\vdash \bq ^2\preceq \bu ^2$.
If $E(\bu )\cap T(\bu )\ne\varnothing$, use $\bu \approx \bu ^2$ and $\bq \preceq \bq ^2$.
If $\mathrm{(O)}_{\bu}$ holds but $\mathrm{(O)}_{\bu +\bq }$ fails, then
in $\mathcal W$,
$\bu \approx \bu +\bq \approx(\bu +\bq )^2\approx \bu ^2$.
Its four-variable consequence is represented in $\Delta$, and
the reverse substitution proved above gives $\Delta\vdash \bu \approx \bu ^2$.
Again $\bq \preceq \bu $ follows.
If $\mathrm{(O)}_{\bu +\bq }$ holds, the 260-variable reduction and its
reverse substitution give $\bq \preceq \bu $ directly.
Every term identity is equivalent to finitely many identities
saying that each summand on either side is at most the other side.
Thus $\Delta$ is a finite equational basis for $\mathcal W$.
\end{proof}

\subsection{The case $\mathcal R_{\mathcal M}$: a 15-variable bound}

In $\mathcal R_{\mathcal M}$, square substitution in $x+y\approx xy$,
followed by right multiplication by $z$, gives
\begin{equation}\label{eq:j0-degree2}
 xyz\approx xz+yz.
\end{equation}
Every term is therefore equal to a sum of one-letter and two-letter
words. Regard $xy$ as an edge between first variable $x$ and last variable $y$.
Under $\mathrm{(O)}_{\bu}$ these sets are disjoint and the graph is bipartite.
Mark a last-variable vertex when it also occurs as a one-letter summand;
retain isolated one-letter summands.
If $E(\bu )\cap T(\bu )\ne\varnothing$, then
\begin{equation}\label{eq:j0-content-square}
 \bu \approx \mathbf{k}_{c(\bu )}:=\left(\sum_{x\in c(\bu )}x\right)^2,
\end{equation}
because $\bu \approx \bu ^2$ and the two operations agree in the square image.
Put
\[
 \rho:xt\preceq xy+zy+zt,\qquad \pi_0:t\preceq xt.
\]

\begin{lem}\label{lem:j0-rect-hfb}
The variety $\mathcal R_{\mathcal M}[\rho]$ is hereditarily finitely based.
\end{lem}
\begin{proof}
Suppose that $\mathrm{(O)}_{\bu}$ holds.
Within a component of its bipartite graph, vertices in opposite parts
are joined by a path of odd length.
The identity $\rho$ adds an edge between the endpoints of any
three-edge segment. Repeated application joins the endpoints of any
such path. Thus $\bu $ equals the term consisting of the complete
bipartite graph on each component, with the original marked vertices
and isolated one-letter summands. No variables or marked vertices
are added. Otherwise use \eqref{eq:j0-content-square}.

An arbitrary identity is equivalent, using \eqref{eq:j0-degree2},
to finitely many identities $\bq \preceq \bu $ with $\ell(\bq)\leq2$.
Fix one and leave its at most two variables unchanged.
If $E(\bu )\cap T(\bu )\ne\varnothing$, identify all other variables of
$c(\bu )$ with one new variable. Replacing it by the sum of its original
members recovers \eqref{eq:j0-content-square} exactly.
At most three variables are needed, including those of $\bq $ outside $\bu $.

Suppose now $\mathrm{(O)}_{\bu}$, and take the complete bipartite components.
Distinguish nonisolated components having a marked last-variable
vertex from those having none.
Retain components containing variables of $\bq $, at most two.
Also retain a representative of each occurring type, at most two more.
Retain every isolated one-letter summand.
Map first-part vertices of each other component to a first-part
vertex in a retained component of the same type, and all its last
vertices to one last vertex there. In the marked type choose a
marked image. Fix all retained variables and all variables of $\bq $.
Each removed edge or mark maps to an existing summand, so the
substituted upper side is precisely the retained term $\bu _{\rm ret}$.
This derives $\bq \preceq \bu _{\rm ret}$ from the original identity.
Conversely, $\bu _{\rm ret}\preceq \bu $, so it implies $\bq \preceq \bu $.

Within each retained component, group the variables not in $\bq $
into first-part vertices, marked last vertices, and unmarked last
vertices, giving at most three nonempty groups.
Group the remaining isolated one-letter variables together.
Identify each group with a distinct new variable and retain the
variables of $\bq $. This gives at most
\[
 4\cdot3+1+2=15
\]
variables. Replace each group variable back by the sum of its
original members. Distributivity recovers exactly the complete
bipartite components, the marks, and the isolated summands.
The word $\bq $ remains unchanged. This proves the required equivalence.

For a subvariety, add to the finite equational basis of
$\mathcal R_{\mathcal M}[\rho]$ every identity of that subvariety of
the form $\bq \preceq \bu $ on a fixed 15-variable alphabet with summands
of length at most two. There are only $15+15^2$ such words and
finitely many nonempty sets of them.
This is therefore a finite set, and the preceding reductions
show that it implies all identities of the subvariety.
\end{proof}

\begin{lem}\label{lem:j0-M-hfb}
The variety $\mathcal R_{\mathcal M}$ is hereditarily finitely based.
\end{lem}
\begin{proof}
Subvarieties not containing $M_0$ are covered by Lemma~\ref{lem:j0-omit}.
Assume $M_0\in\mathcal W$.
If $E(\bu )\cap T(\bu )\ne\varnothing$, an identity $\bq \preceq \bu $
already follows from the content condition and \eqref{eq:j0-content-square}.
For the other identities use summands of length at most two.

Suppose an additional identity $xt\preceq \bu $ has $\mathrm{(O)}_{\bu}$ but
the edge $xt$ is absent.
The criterion of $M_0$ gives $x\in E(\bu )$ and $t\in T(\bu )$.
Keep $x,t$, identify all other first-part variables with $z$,
and all other last variables with $y$.
The upper side becomes a sum of summands of $xy+zy+zt+y+t$.
Next substitute $y\mapsto zy$ and $t\mapsto zt$.
By \eqref{eq:j0-degree2}, the left side is $xt+zt$ and every
right-hand summand is among $xy,zy,zt$. Adding missing terms gives $\rho$.
If the lower side is a one-letter word $t$ that is not a summand of $\bu $,
then $t\in T(\bu )$ by Lemma~\ref{lem:j0-M0}.
Map all first-part variables to $x$ and last variables other than
$t$ to $xt$, keeping $t$. Every summand becomes $xt$, giving $\pi_0$.

Consider first a proper subvariety of
$\mathcal R_{\mathcal M}[\pi_0]$ containing $M_0$.
Choose an identity failing in that containing variety and reduce
it to a single-summand identity as above.
The upper side must satisfy $\mathrm{(O)}_{\bu}$.
An identity with a lower word of length one cannot fail: if $t\in T(\bu )$,
either $t$
is a summand or some $xt$ is a summand, and $\pi_0$ already gives
$t\preceq \bu $. Thus an absent edge must be added and $\rho$ follows.
It follows from Lemma~\ref{lem:j0-rect-hfb} that this subvariety is finitely based.
The variety $\mathcal R_{\mathcal M}[\pi_0]$ itself is finitely based
by definition, so it is hereditarily finitely based.
Finally a proper subvariety of $\mathcal R_{\mathcal M}$ containing
$M_0$ implies $\rho$ or $\pi_0$ by the two arguments above.
Both resulting containing varieties are hereditarily finitely based.
\end{proof}

\subsection{Three consequences of an additional identity}

In the next lemma, $\mathcal W\subseteq\mathcal J_0$ contains $M_0$ and
satisfies $\bs t\preceq \bu $, where $\bs $ is nonempty and $\mathrm{(O)}_{\bu}$ holds.
Thus $c(\bs )\subseteq E(\bu )$ and $t\in T(\bu )$.

\begin{lem}\label{lem:j0-extractions}
Under these assumptions:
\begin{enumerate}
\item if $\mathrm A(\bu ,\bs t)$ holds but $\mathrm A_+(\bu ,\bs t)$ fails,
      then $\mathcal W\models\lambda$;
\item if $c(\bs )\not\subseteq E_t(\bu )$, then $\mathcal W\models\tau$;
\item if $c(\bs )\not\subseteq O_t$ and $\mathcal W$ lies in one of
      $\mathcal R_{\mathcal K}$, $\mathcal R_{\mathcal Q}$,
      $\mathcal R_{\mathcal B}$, $\mathcal R_{\mathcal L\vee\mathcal M}$,
      then $\mathcal W\models\iota$.
\end{enumerate}
\end{lem}
\begin{proof}
For item~1, the summand supplied by $\mathrm A$ must be $t$ itself.
Map $c(\bs )$ to $x$, $E(\bu )\setminus c(\bs )$ to the nonempty word $xy$,
$t$ to $z$, and every other last variable to $xyz$.
Each longer original summand ending in $t$ contains a variable
outside $c(\bs )$, so becomes $xyz$ after repeated-prefix deletion.
Other summands become $xyz$, except $t$, which becomes $z$.
The lower word is $xz$. Adding a missing right-hand term gives $\lambda$.

For item~2, choose $b\in c(\bs )\setminus E_t(\bu )$.
Map $b$ to $xy$, other members of $E(\bu )$ to $x$, $t$ to $z$,
and other last variables to $xyw$.
All prefix images start with $x$; repeated-prefix deletion shows
that summands ending in $t$ become $z$ or $xz$, other summands
become $xyw$, and the lower word becomes $xyz$.
Substitute $z\mapsto xz$ to obtain $\tau$.

For item~3, map first-part vertices in the component of $t$ in $G(\bu )$
to $x$, the others to $y$, last vertices in that component to $xz$,
and other last vertices to $yw$.
No summand connects distinct components, so the upper side maps to a
sum of summands of $\bv =xz+yw$. The resulting lower word $\bq '$ ends in $z$ and
its prefix contains $x,y$: $x$ comes from the image of $t$,
and $y$ from $c(\bs )\setminus O_t$. Hence $\bq '\preceq \bv $ holds in $\mathcal W$.
For the lower word $xyz$ and upper side $\bv +\bq '$, the prefix set for last variable
$z$ contains $\{x,y\}$, the summand $xz$ gives $\mathrm A_+$,
and $x$ is an original first variable.
Condition $\operatorname{Ex}$ is supplied by $xz$ for $x,z$
and $yw$ for $y$.
Each case of Lemma~\ref{lem:j0-carriers} gives $xyz\preceq \bv +\bq '\preceq \bv $.
This is $\iota$. An isolated original $t$ is included, since its
image still supplies the variable $x$.
\end{proof}

\subsection{The case $\mathcal R_{\mathcal K}$}

\begin{lem}\label{lem:j0-K-hfb}
The variety $\mathcal R_{\mathcal K}$ is hereditarily finitely based.
\end{lem}
\begin{proof}
Put $R=\mathcal R_{\mathcal K}$ and $H=R[\tau]$.
Every variety used below is obtained from $R$ by adding finitely many
identities and is finitely based.
Since $R[\iota]\subseteq\mathcal U$, Theorem~\ref{thm:j0-global-hfb}
applies to all subvarieties of $R[\iota]$.

Subvarieties not containing $M_0$ are covered by Lemma~\ref{lem:j0-omit}.
If $M_0\in\mathcal W$ and $\mathcal U_{\mathcal W}\subsetneq\mathcal K$,
Lemma~\ref{lem:band12} gives $\mathcal U_{\mathcal W}=\mathcal M$,
so $\mathcal W\subseteq\mathcal R_{\mathcal M}$ and
Lemma~\ref{lem:j0-M-hfb} applies.
It remains to consider proper subvarieties containing $M_0$ with
$\mathcal U_{\mathcal W}=\mathcal K$.
Choose an identity of $\mathcal W$ failing in the particular containing
variety under consideration, and reduce it to $\bq \preceq \bu $.
Lemma~\ref{lem:j0-closed-band} forces $\mathrm{(O)}_{\bu}$.
A lower word of length one is impossible: $\mathrm A(\bu ,\bq )$, forced by
$\mathcal K$, then makes it an original summand.
Thus $\bq =\bs t$ with $\bs $ nonempty; the occurrence conditions and $\mathrm A$
hold by $M_0$ and $\mathcal K$.

For $H[\lambda]$, Lemma~\ref{lem:j0-component-carriers} shows that
$c(\bs )\subseteq O_t$ and $\mathrm A$ suffice.
Failure therefore forces $c(\bs )\not\subseteq O_t$ and implies $\iota$
by Lemma~\ref{lem:j0-extractions}.
The resulting subvariety is finitely based by the first paragraph,
so $H[\lambda]$ is hereditarily finitely based.
For $H$, failure of $\mathrm A_+$ implies $\lambda$ and gives a
subvariety of $H[\lambda]$.
If $\mathrm A_+$ holds, failure forces $c(\bs )\not\subseteq O_t$
and implies $\iota$. Thus $H$ is hereditarily finitely based.
For $R[\lambda]$, the condition $c(\bs )\subseteq E_t(\bu )$ together
with $\mathrm A$ already suffices. Failure implies $\tau$,
giving a subvariety of $H[\lambda]$.
Finally, for $R$, failure of $\mathrm A_+$ implies $\lambda$;
otherwise failure forces $c(\bs )\not\subseteq E_t(\bu )$ and implies $\tau$.
The containing varieties in both cases have just been treated.
Each variety under consideration is itself finitely based, so the argument
covers all of its subvarieties.
\end{proof}

\subsection{An identity implied by failure of the $A_t$ condition}

\begin{lem}\label{lem:j0-head-extraction}
Let $R$ be $\mathcal R_{\mathcal Q}$ or
$\mathcal R_{\mathcal L\vee\mathcal M}$.
Suppose that $\mathcal W\subseteq R$, $M_0\in\mathcal W$,
$\mathcal L\subseteq\mathcal U_{\mathcal W}$, and
$\mathcal W\models \bs t\preceq \bu $, where $\bs \ne\varepsilon$,
$\mathrm{(O)}_{\bu}$ holds, and $c(\bs )\not\subseteq A_t$.
Then $\mathcal W\models\nu$.
\end{lem}
\begin{proof}
The condition from $\mathcal L$ gives $h(\bs )\in H(\bu )$:
an original summand starting in $h(\bs )$ cannot have length one,
because $\mathrm{(O)}_{\bu}$ holds and $h(\bs )\in E(\bu )$.
Choose $b\in c(\bs )\setminus A_t$.
Map members of $H(\bu )$ in the component of $t$ in $G_H(\bu )$ to $x$,
and those outside that component to $y$.

First suppose $b\notin H(\bu )$.
Keep its image as a separate variable $b$, map every other member of
$E(\bu )\setminus H(\bu )$ to $a$, last variables in the component to $xat$,
and other last variables to $yabw$.
A summand ending in the component cannot contain a member of $H(\bu )$
outside that component, nor can it contain $b$.
Its image starts in $x$ and equals $xat$ after the image of its
last variable is appended.
A longer summand ending outside starts in $y$ and contains no member
of $H(\bu )$ from the component. Permuting after its first variable
and deleting repeated prefix occurrences gives $yabw$.
The one-letter summands have these two images as well.
Thus the substituted identity implies
\[
 \bq '\preceq \bv =xat+yabw.
\]
The prefix of $\bq '$ contains $x,a,b$ and may contain $y$.
For the lower word $xabt$ and upper side $\bv +\bq '$, the prefix set for last
variable $t$ contains the needed variables, $x$ is an original
first variable, and $xat$ supplies $\mathrm A_+$.
Lemma~\ref{lem:j0-carriers} gives $xabt\preceq \bv +\bq '\preceq \bv $.
Substituting $b\mapsto y$ and deleting repeated $y$ proves $\nu$.

If $b\in H(\bu )$, it lies outside the component.
Map $E(\bu )\setminus H(\bu )$ to $a$, last variables outside to $yaw$,
and last variables inside to $xat$.
The upper side becomes a sum of summands of $xat+yaw$,
and the lower word ends in $t$ with prefix containing $x,a,y$.
Apply Lemma~\ref{lem:j0-carriers} to the enlarged upper side as before,
using $xat$ for $\mathrm A_+$ and the first-variable condition.
It gives $xayt\preceq xat+yaw$, namely $\nu$.
An isolated original $t$ is included: its image $xat$ still supplies
the summand used in the argument.
\end{proof}

\subsection{The case $\mathcal R_{\mathcal L\vee\mathcal M}$}

\begin{lem}\label{lem:j0-LM-hfb}
The variety $\mathcal R_{\mathcal L\vee\mathcal M}$ is hereditarily
finitely based.
\end{lem}
\begin{proof}
Put $R=\mathcal R_{\mathcal L\vee\mathcal M}$,
$H=R[\tau]$, and $G=R[\nu]$.
Since $R[\iota]\subseteq\mathcal U$, Theorem~\ref{thm:j0-global-hfb}
applies to all subvarieties of $R[\iota]$.
Subvarieties not containing $M_0$, or having multiplicatively idempotent
subvariety $\mathcal M$, are treated by Lemmas~\ref{lem:j0-omit}
and \ref{lem:j0-M-hfb}.
Lemma~\ref{lem:band12} shows that these cover all cases except
$M_0\in\mathcal W$ and
$\mathcal U_{\mathcal W}=\mathcal L\vee\mathcal M$.

In this remaining case, choose an identity valid in a proper subvariety
but failing in the containing variety, and reduce it to $\bq \preceq \bu $.
Lemma~\ref{lem:j0-closed-band} forces $\mathrm{(O)}_{\bu}$.
If $\bq =t$ had length one, $M_0$ would give $t\in T(\bu )$ and
$\mathcal L$ would require a summand beginning in $t$.
Disjointness then forces that summand to be $t$ itself, a contradiction.
Hence $\bq =\bs t$ with $\bs $ nonempty.
The occurrence conditions and $h(\bs )\in H(\bu )$ are already satisfied.

In $G$, the condition $c(\bs )\subseteq O_t$ suffices by item~4 of
Lemma~\ref{lem:j0-carriers} and Lemma~\ref{lem:j0-component-carriers}.
Failure implies $\iota$ by Lemma~\ref{lem:j0-extractions}.
Thus $G$ is hereditarily finitely based.
In $H$, the condition $c(\bs )\subseteq A_t$ suffices; failure implies
$\nu$ by Lemma~\ref{lem:j0-head-extraction}, giving a subvariety of $G$.
Thus $H$ is hereditarily finitely based.
In $R$, the condition $c(\bs )\subseteq E_t(\bu )$ suffices; failure implies
$\tau$, giving a subvariety of $H$.
Each variety under consideration is finitely based by definition, so this proves
the assertion for every subvariety.
\end{proof}

\subsection{The case $\mathcal R_{\mathcal Q}$}

\begin{lem}\label{lem:j0-Q-hfb}
The variety $\mathcal R_{\mathcal Q}$ is hereditarily finitely based.
\end{lem}
\begin{proof}
Put $R=\mathcal R_{\mathcal Q}$, $H=R[\tau]$, and $G=R[\nu]$.
The variety $R[\iota]=\mathcal U$ was treated in
Theorem~\ref{thm:j0-global-hfb}.
By Lemma~\ref{lem:band12}, a proper subvariety of $\mathcal Q$
containing $\mathcal M$ is contained in $\mathcal K$ or
$\mathcal L\vee\mathcal M$.
Hence Lemmas~\ref{lem:j0-K-hfb} and \ref{lem:j0-LM-hfb}
treat $\mathcal U_{\mathcal W}\subsetneq\mathcal Q$ when $M_0\in\mathcal W$.
Lemma~\ref{lem:j0-omit} treats $M_0\notin\mathcal W$.

It remains to consider $M_0\in\mathcal W$ and
$\mathcal U_{\mathcal W}=\mathcal Q$.
An additional identity failing in the containing variety must be
$\bs t\preceq \bu $ with $\bs $ nonempty and $\mathrm{(O)}_{\bu}$.
The occurrence conditions, $h(\bs )\in H(\bu )$, and $\mathrm A(\bu ,\bs t)$
are forced by $M_0,\mathcal Q$.
All six varieties below are finitely based by definition; treat
their proper subvarieties in this order:
\begin{enumerate}
\item In $G[\lambda]$, $c(\bs )\subseteq O_t$ and $\mathrm A$ suffice.
      Failure implies $\iota$, so this variety is hereditarily finitely based.
\item In $G$, failure of $\mathrm A_+$ implies $\lambda$.
      Otherwise failure forces $c(\bs )\not\subseteq O_t$ and implies $\iota$.
      It follows from item~1 or Theorem~\ref{thm:j0-global-hfb} that the
      subvariety is finitely based.
\item In $H[\lambda]$, $c(\bs )\subseteq A_t$ and $\mathrm A$ suffice.
      Failure implies $\nu$, giving a subvariety of $G[\lambda]$.
\item In $H$, failure of $\mathrm A_+$ implies $\lambda$.
      Otherwise failure forces $c(\bs )\not\subseteq A_t$ and implies $\nu$.
      The containing varieties $H[\lambda]$ and $G$ are already treated.
\item In $R[\lambda]$, $c(\bs )\subseteq E_t(\bu )$ and $\mathrm A$ suffice.
      Failure implies $\tau$, giving a subvariety of $H[\lambda]$.
\item In $R$, failure of $\mathrm A_+$ implies $\lambda$.
      Otherwise failure forces $c(\bs )\not\subseteq E_t(\bu )$ and implies $\tau$,
      giving a subvariety of $H$.
\end{enumerate}
Sufficiency follows from Lemmas~\ref{lem:j0-carriers} and
\ref{lem:j0-component-carriers}; the required consequences follow
from Lemmas~\ref{lem:j0-extractions} and \ref{lem:j0-head-extraction}.
Every proper subvariety is finitely based by one of the preceding cases,
and each variety under consideration is itself finitely based.
This completes the proof.
\end{proof}

\subsection{Two cases determined by the content condition}

\begin{lem}\label{lem:j0-S-LS}
We have
\[
 \mathcal R_{\mathcal S_b}=\Var(M_0,S|_{\{2,3,4\}}),\qquad
 \mathcal R_{\mathcal L\vee\mathcal S_b}
   =\Var(M_0,L,S|_{\{2,3,4\}}),
\]
and both varieties are hereditarily finitely based.
\end{lem}
\begin{proof}
The displayed generators belong to the respective varieties on the left.
For the reverse inclusion, consider a common identity $\bq \preceq \bu $ of the
generators on the right.
If $E(\bu )\cap T(\bu )\ne\varnothing$, the content condition from $M_0$
is implied by $\mathrm C(\bu ,\bq )$. Apply the identity in the square image,
then use $\bq \preceq \bq ^2$ and $\bu \approx \bu ^2$.
If $\mathrm{(O)}_{\bu}$ and $\bq $ is a single letter,
$\mathrm C(\bu ,\bq )$ forces it to be an original summand.

Suppose that $\mathrm{(O)}_{\bu}$ and $\bq =\bs t$ with $\bs $ nonempty.
Every summand whose content is contained in $c(\bq )$ ends in $t$.
The prefixes of those of length at least two cover exactly $c(\bs )$.
Combining them by \eqref{eq:j0-ST} gives a word with last variable $t$
and prefix set $c(\bs )$.
For $\mathcal S_b$, use $xyz\approx yxz$ to obtain the required prefix order.
For $\mathcal L\vee\mathcal S_b$, an original summand starts with $h(\bs )$.
Equation \eqref{eq:j0-HI} puts that variable first, and $xyzt\approx xzyt$
arranges the others. This proves the two equalities.
Both varieties are finitely based by the definition of $\mathcal R_{\mathcal U}$.

For a subvariety of the first, if $M_0$ is absent use Lemma~\ref{lem:j0-omit}.
If it contains $M_0$ and its multiplicatively idempotent subvariety
is $\mathcal S_b$, it contains all displayed generators and therefore
equals $\mathcal R_{\mathcal S_b}$.
Otherwise that subvariety is contained in $\mathcal K$,
and Lemma~\ref{lem:j0-K-hfb} applies.
For the second variety, equality
$\mathcal U_{\mathcal W}=\mathcal L\vee\mathcal S_b$, together with
$M_0\in\mathcal W$, likewise gives
$\mathcal W=\mathcal R_{\mathcal L\vee\mathcal S_b}$.
Every proper multiplicatively idempotent subvariety lies in
$\mathcal S_b$ or $\mathcal Q$ by Lemma~\ref{lem:band12}.
Use the first assertion just proved or Lemma~\ref{lem:j0-Q-hfb},
respectively. This covers all cases.
\end{proof}

\subsection{Three subvarieties of $\mathcal R_{\mathcal B}$}

In $\mathcal R_{\mathcal B}$, distinct prefix variables cannot in general
be freely permuted. Put $R=\mathcal R_{\mathcal B}$ and
\[
 \pi:xz\preceq z+xyw,\qquad
 K_0=R[\iota],\qquad K_1=R[\iota,\lambda],\qquad
 K_2=R[\iota,\lambda,\pi].
\]
Each is finitely based by definition.

\begin{lem}\label{lem:j0-B-global}
The varieties $K_0,K_1,K_2$ are hereditarily finitely based, and
\[
 K_2=\Var(M_0,B),\qquad B=S|_{\{0,2,4\}}.
\]
\end{lem}
\begin{proof}
We first give three sufficient conditions for $\bs t\preceq \bu $ under
$\mathrm{(O)}_{\bu}$ and $\bs \ne\varepsilon$.
Assume the occurrence conditions from $M_0$ and
$\operatorname{Ex}(\bs t,\bu )$ from $\mathcal B$.

If $\mathrm A_+(\bu ,\bs t)$ holds, start with its summand $\br t$.
For each original summand $\mathbf{b} w$ of length at least two, use $\iota$
with $x$ equal to the current prefix, $y=\mathbf{b} $, $z=t$, and the
displayed $w$ as fourth argument.
This gives a word ending in $t$ whose prefix contains $E(\bu )$.
Use Lemma~\ref{lem:carrier-reorder} to put the variables of $\bs $
first in the required order. Its nonfinal-occurrence hypotheses
follow from $\mathrm{(O)}_{\bu}$ and the original summands giving
$\operatorname{Ex}(\bs t,\bu )$. We obtain $\bs \mathbf{z} t\preceq \bu $.
The original summand from $\mathrm A_+$ and \eqref{eq:j0-trim}
remove the unnecessary $\mathbf{z} $. Thus the identity holds in $K_0$.

More generally, suppose at least one original summand of length
at least two ends in $t$. Starting there gives $\bs \mathbf{z} t\preceq \bu $ in the
same way. The criterion of $\mathcal B$ provides $\mathrm A(\bu ,\bs t)$.
If $\mathrm A_+$ holds, the previous argument applies.
Otherwise $t$ is an original one-letter summand, and $\lambda$
gives $\bs t\preceq t+\bs \mathbf{z} t\preceq \bu $. Thus the identity holds in $K_1$.

Now allow $t$ to occur only as an isolated one-letter summand.
At the first variable $h=h(\bs )$, $\operatorname{Ex}$ gives an
original summand $h\br $. It has length at least two:
$h\in E(\bu )$ and cannot belong to $T(\bu )$.
Substitute $x=h$, $y=h$, $z=t$, $w=\br $ into $\pi$ to obtain
\[
 ht\preceq t+hh\br \approx t+h\br \preceq \bu .
\]
Adjoining $ht$ makes $\mathrm A_+$ hold for $\bs t$; the first
construction now gives $\bs t\preceq \bu $.
If $E(\bu )\cap T(\bu )\ne\varnothing$, common identities of $M_0,B$
follow by applying their identity in the square image.
If $\mathrm{(O)}_{\bu}$ and the lower word has length one, the criterion of
$B$ makes it an original summand.
Thus $K_2$ satisfies all common identities of $M_0$ and $B$.
Conversely these algebras lie in $R$ and satisfy $\iota,\lambda,\pi$.
Hence $K_2=\Var(M_0,B)$.

To prove that $K_0,K_1,K_2$ are hereditarily finitely based, first apply
Lemma~\ref{lem:j0-omit} to subvarieties not containing $M_0$.
If $M_0\in\mathcal W$ and $\mathcal U_{\mathcal W}\subsetneq\mathcal B$,
then $\mathcal U_{\mathcal W}\subseteq\mathcal Q$, the greatest proper
subvariety of $\mathcal B$; use Lemma~\ref{lem:j0-Q-hfb}.
It remains to consider $M_0\in\mathcal W$ and
$\mathcal U_{\mathcal W}=\mathcal B$.
Such a subvariety of $K_2$ contains $M_0,B$ and therefore equals
$K_2$, which is finitely based. This proves the result for $K_2$.

For a proper subvariety in the remaining case inside $K_1$,
choose a failing single-summand identity $\bs t\preceq \bu $.
The identity already holds when $E(\bu )\cap T(\bu )\ne\varnothing$ or
when the lower word has length one, by the criteria for $M_0$ and $B$.
The second sufficient condition shows that no original longer
summand ends in $t$. Since $t\in T(\bu )$, it is isolated and one-letter.
Map $E(\bu )$ to $x$, $t$ to $z$, and all other last variables to $xyw$.
The upper side becomes a sum of summands of $z+xyw$ and the lower word $xz$,
proving $\pi$. Thus the subvariety lies in $K_2$.
Since $K_1$ itself is finitely based, it is hereditarily finitely based.
Finally, for a proper subvariety in the remaining case inside $K_0$,
the first sufficient condition forces failure of $\mathrm A_+$.
Condition $\mathrm A$ holds by $B$, so Lemma~\ref{lem:j0-extractions}
gives $\lambda$ and a subvariety of $K_1$.
Thus $K_0$ is hereditarily finitely based.
\end{proof}

\subsection{The $A_t$ condition without prefix commutativity}

\begin{lem}\label{lem:j0-B-head}
Suppose that $\mathcal W\subseteq R=\mathcal R_{\mathcal B}$,
$M_0\in\mathcal W$, and $\mathcal U_{\mathcal W}=\mathcal B$.
If $\mathcal W\models \bs t\preceq \bu $ with $\bs \ne\varepsilon$,
$\mathrm{(O)}_{\bu}$, and $c(\bs )\not\subseteq A_t$, then $\mathcal W\models\nu$.
\end{lem}
\begin{proof}
Let $b$ be the first variable of $\bs $ that does not belong to $A_t$.
Condition $\operatorname{Ex}(\bs t,\bu )$ supplies an original summand
$\br b\mathbf{z} $ with $c(\br )$ contained in the preceding variables of $\bs $,
and hence in $A_t$. Since $b\in E(\bu )$ and $\mathrm{(O)}_{\bu}$,
that occurrence is nonfinal and $\mathbf{z} $ is nonempty.

If $\br $ were nonempty, its first variable $h$ would lie in
$H(\bu )\cap A_t$. By definition of $A_t$, it would be adjacent in
$G_H(\bu )$ to a last vertex in the component of $t$.
It would therefore be in that component itself.
The original summand $\br b\mathbf{z} $ connects its last variable to $h$,
so that last variable would be in the same component.
As $b$ occurs in its prefix, this implies $b\in A_t$, a contradiction.
Thus $\br $ is absent and $b$ is a member of $H(\bu )$ outside the component.

Map $H(\bu )$ inside the component to $x$, $H(\bu )$ outside to $y$,
and $E(\bu )\setminus H(\bu )$ to $a$.
Map last variables inside to $xat$, and outside to $yaw$.
A summand ending inside begins in $x$ and contains no original
first variable mapped to $y$, so its image is $xat$.
A summand ending outside begins in $y$ and contains no original
first variable mapped to $x$, so its image is $yaw$.
Only two distinct prefix variables occur in either case;
repeated-prefix deletion suffices, without permuting distinct variables.
The one-letter summands have these images as well.

We obtain $\bq '\preceq \bv :=xat+yaw$.
The prefix of $\bq '$ contains $x,a,y$: $x,a$ come from the image of
the original $t$, and $y$ from $b$.
For the lower word $xayt$ and upper side $\bv +\bq '$, the condition on
prefixes ending in the same variable holds, and $xat$ gives $\mathrm A_+$.
For $\operatorname{Ex}$, use $xat$ at $x,a,t$ and $yaw$ at $y$.
Item~3 of Lemma~\ref{lem:j0-carriers} gives
$xayt\preceq \bv +\bq '\preceq \bv $, which is $\nu$.
This includes an isolated original $t$.
\end{proof}

\begin{thm}\label{thm:j0-B-hfb}
The variety $\mathcal R_{\mathcal B}$ is hereditarily finitely based.
\end{thm}
\begin{proof}
Put $R=\mathcal R_{\mathcal B}$, $H=R[\tau]$, and $G=R[\nu]$.
Lemma~\ref{lem:j0-B-global} proves the assertion for $R[\iota]$.
Subvarieties without $M_0$ are covered by Lemma~\ref{lem:j0-omit};
those with $M_0$ and $\mathcal U_{\mathcal W}\subsetneq\mathcal B$ are
covered by Lemma~\ref{lem:j0-Q-hfb}.
Consider $M_0\in\mathcal W$ and
$\mathcal U_{\mathcal W}=\mathcal B$.

An identity of a proper subvariety that fails in the containing
variety can be taken as $\bs t\preceq \bu $ with $\bs \ne\varepsilon$
and $\mathrm{(O)}_{\bu}$.
The occurrence conditions, $\operatorname{Ex}(\bs t,\bu )$, and
$\mathrm A(\bu ,\bs t)$ follow from $M_0,B$.
Lemmas~\ref{lem:j0-carriers} and \ref{lem:j0-component-carriers}
give these sufficient additional conditions:
\[
\begin{array}{c|c|c}
 \text{variety}&\text{prefix condition}&\text{summand condition}\\\hline
 H&c(\bs )\subseteq A_t&\mathrm A_+\\
 H[\lambda]&c(\bs )\subseteq A_t&\mathrm A\\
 G&c(\bs )\subseteq O_t&\mathrm A_+\\
 G[\lambda]&c(\bs )\subseteq O_t&\mathrm A.
\end{array}
\]
If the identity fails in one of these varieties, at least one of the
listed conditions fails. Every variety below is finitely based by definition.
\begin{enumerate}
\item For $G[\lambda]$, failure forces $c(\bs )\not\subseteq O_t$,
      so implies $\iota$.
\item For $G$, failure of $\mathrm A_+$ implies $\lambda$.
      Otherwise failure forces $c(\bs )\not\subseteq O_t$, so implies $\iota$.
\item For $H[\lambda]$, failure forces $c(\bs )\not\subseteq A_t$
      and implies $\nu$ by Lemma~\ref{lem:j0-B-head}.
      The subvariety lies in $G[\lambda]$.
\item For $H$, failure of $\mathrm A_+$ implies $\lambda$.
      Otherwise failure forces $c(\bs )\not\subseteq A_t$ and implies $\nu$.
      The containing varieties $H[\lambda]$ and $G$ are already treated.
\item For $R[\lambda]$, $c(\bs )\subseteq E_t(\bu )$, $\mathrm A$,
      and $\operatorname{Ex}$ suffice.
      Failure implies $\tau$ and a subvariety of $H[\lambda]$.
\item For $R$, failure of $\mathrm A_+$ implies $\lambda$.
      Otherwise failure forces $c(\bs )\not\subseteq E_t(\bu )$ and implies $\tau$.
\end{enumerate}
The consequences $\lambda,\tau,\iota$ come from
Lemma~\ref{lem:j0-extractions}. Each resulting containing variety
has been proved hereditarily finitely based earlier in this order.
Consequently every subvariety is finitely based, as required.
\end{proof}

\subsection{Completion of the proof for $\mathcal J_0$}

\begin{thm}\label{thm:j0-hfb}
The variety $\mathcal J_0=\Var(M_0,C_4)$ is hereditarily
finitely based.
\end{thm}
\begin{proof}
Let $\mathcal W\subseteq\mathcal J_0$ be arbitrary.
If $M_0\notin\mathcal W$, apply Lemma~\ref{lem:j0-omit}.
Otherwise $\mathcal U_{\mathcal W}$ contains $\mathcal M$
and is contained in $\Var(C_4)$.
If it equals $\Var(C_4)$, then $C_4\in\mathcal W$ and
\[
 \mathcal J_0=\Var(M_0,C_4)
 \subseteq\mathcal W\subseteq\mathcal J_0.
\]
Thus $\mathcal W=\mathcal J_0$, which is finitely based.
Otherwise Lemma~\ref{lem:band12} gives
\[
 \mathcal U_{\mathcal W}\subseteq\mathcal B
 \quad\text{or}\quad
 \mathcal U_{\mathcal W}\subseteq\mathcal L\vee\mathcal S_b.
\]
Consequently $\mathcal W$ is contained in
$\mathcal R_{\mathcal B}$ or $\mathcal R_{\mathcal L\vee\mathcal S_b}$.
Theorem~\ref{thm:j0-B-hfb} and Lemma~\ref{lem:j0-S-LS}
show that both varieties are hereditarily finitely based.
Hence $\mathcal W$ is finitely based in every case.
\end{proof}

The proof first establishes the 260-variable and 15-variable reductions.
It then treats square images in $\mathcal M$, $\mathcal K$,
$\mathcal L\vee\mathcal M$, and $\mathcal Q$, followed by
$\mathcal S_b$, $\mathcal L\vee\mathcal S_b$, and $\mathcal B$.
The classification of $\mathcal U_{\mathcal W}$ completes the argument.
\section{The two maximal proper subvarieties are hereditarily finitely based}
\label{sec:transfers}

Throughout this section, we consider algebras with the two binary operations
$+$ and $\cdot$ and no constant symbols. We retain the notation
\[
 \mathcal V=\Var(S),\qquad
 \mathcal A=\Var(A_4),\qquad
 \mathcal C=\Var(M_0,M_1,C_4),
\]
\[
 \mathcal D=\Var(M_1,C_4),\qquad
 \mathcal J_0=\Var(M_0,C_4),\qquad
 \mathcal E=\mathcal V\cap\operatorname{Mod}(xy\preceq x+xyz).
\]
Here $\operatorname{Mod}$ is taken relative to the ai-semiring axioms, and
$\bq\preceq \bu$ abbreviates the identity $\bu+\bq\approx \bu$. The join of two varieties
is the least variety containing them; their meet is their intersection.

For a word $\bw$, let $c(\bw)$ denote its content, $h(\bw)$ its first letter,
$t(\bw)$ its last letter, and $e(\bw)$ the set of letters occurring before its
last position. For a term $\bu$, define
\[
 E(\bu)=\bigcup_{\bw\in \bu}e(\bw),\qquad
 T(\bu)=\{t(\bw):\bw\in \bu\}.
\]
As usual, a term is identified with its nonempty finite set of word
summands. We write $\mathrm{SA}(\bu,\bq)$ when at least one of the following
explicit conditions holds:
\begin{enumerate}
 \item[\textnormal{(SA1)}] There is a summand $\mathbf{u}_1$ of $\bu$ with
 $\ell(\mathbf{u}_1)\geq2$ and $c(\mathbf{u}_1)\subseteq c(\bq)$.
 \item[\textnormal{(SA2)}] There is $a\in c(\bq)$ such that $a$ and a word
 $a\br$, $\br\ne\varepsilon$, are both summands of $\bu$.
\end{enumerate}
Here $\varepsilon$ denotes the empty word. It may occur to the left or right
of a term, but no variable is replaced by $\varepsilon$. We use
$\mathrm{SA}(\bu,\bq)$ as an abbreviation for the two alternatives above.

We shall use the following results established in the preceding sections.
\begin{enumerate}
 \item The varieties $\mathcal C$ and $\mathcal E$ are finitely based by
 Theorems~\ref{thm:C-fb} and~\ref{thm:E-fb}. Moreover,
 \[
 \mathcal D=\mathcal C\cap\operatorname{Mod}(xy\approx(xy)^2),\qquad
 \mathcal J_0=\mathcal C\cap\operatorname{Mod}(x^2\preceq x+y^2).
 \]
 The variety $\mathcal D$ is hereditarily finitely based by
 \eqref{eq:D-HFB-conclusion}, and $\mathcal J_0$ is hereditarily finitely
 based by Theorem~\ref{thm:j0-hfb}.
 \item In $\mathcal C$ the identities
 \[
 x^2y\approx xy,\qquad xyx\approx xy^2,\qquad x\preceq x^2,
 \qquad (x+y)^2\approx x^2+y^2
 \]
 hold. Lemma~\ref{lem:C-closing} has the following consequence: if
 $t(\mathbf{z})\in e(\mathbf{z})$, $\ell(\mathbf{v}_1)\geq2$, and both $\mathbf{z}\preceq \bu$ and $\mathbf{v}_1\preceq \bu$ have been
 derived, then $\mathbf{v}_1^2\preceq \bu$ follows. No inclusion between $c(\mathbf{z})$ and $c(\mathbf{v}_1)$
 is required.
 \item The criteria for identities $\bq\preceq \bu$ in
 $S,A_4,M_1$ and $\mathcal C$ have been proved. Every subvariety of
 $\mathcal V$ that does not contain the algebra $A_4$ is contained in
 $\mathcal C$.
 \item Let $T$ be the quotient semiring defined by the language in
 Lemma~\ref{lem:D-T-criterion}. It belongs to $\mathcal E\cap\mathcal C$.
 For $\ell(\bq)\geq2$,
 \[
 T\models \bq\preceq \bu
 \quad\Longleftrightarrow\quad
 M_1\models \bq\preceq \bu\ \text{and}\ \mathrm{SA}(\bu,\bq).
 \]
 It is generated by the images of $x,y,z$ and is defined by the language
 \[
 \mathscr L=\{y,z\}\ \cup\ x^+y\{x,y,z\}^*z\{x,y,z\}^*.
 \]
 The quotient definition used here is recalled below. In particular,
 different elements can be distinguished by adjoining suitable words on
 the left and right and testing membership in $\mathscr L$.
\end{enumerate}

The following distinction between the criteria for identities will be used
repeatedly. For a
letter $b$ whose preceding letters form the set $U$, an original summand
may have the form $\mathbf{v}_0b\br$, with $c(\mathbf{v}_0)\subseteq U$ and $b\notin U$.
The condition $\br\ne\varepsilon$ is stronger than $\ell(\mathbf{v}_0b\br)\geq2$.
The criteria for $S,A_4$ require the former at every nonfinal first
occurrence of $b$, unless the content of some original summand is contained
in $U$. The criterion for $M_1$ requires only the latter, but also at the
final first occurrence. We shall always state explicitly which condition
is being used.

\subsection{A quotient semiring defined by a language}

Let $X$ be a nonempty finite alphabet. Write $X^+$ for the set of nonempty
words and $X^*=X^+\cup\{\varepsilon\}$. Let $\mathscr P(X)$ be the set of
nonempty finite subsets of $X^+$. The operations
\[
 \bu+\bv=\bu\cup \bv,\qquad \bu\bv=\{\mathbf{u}_1\mathbf{v}_1:\mathbf{u}_1\in \bu,\ \mathbf{v}_1\in \bv\}
\]
make $\mathscr P(X)$ the free ai-semiring on $X$.
For a language $K\subseteq X^+$, define
\[
 \bu\equiv_K \bv
 \quad\Longleftrightarrow\quad
 \text{for every }\bs,\bt\in X^*,\quad
 \bs \bu\bt\subseteq K\ \Longleftrightarrow\ \bs \bv\bt\subseteq K.
\]
Thus the condition tested is that \emph{every} summand, after adjoining
$\bs$ and $\bt$, belongs to $K$. Empty surrounding words are permitted,
but the elements of $\mathscr P(X)$ and the words substituted for
variables remain nonempty. For related syntactic semiring constructions, see~\cite{Polak2006}.
We give direct proofs of all the properties needed below.

\begin{lem}[Separation by a quotient semiring]\label{lem:syntactic-separation}
 Let $\mathcal W$ be an ai-semiring variety, let $\bu\in\mathscr P(X)$, and put
 \[
 K=\{\bw\in X^+:\mathcal W\models \bw\preceq \bu\}.
 \]
 Then $\equiv_K$ is a congruence on $\mathscr P(X)$, and
 \[
 A_K=\mathscr P(X)/{\equiv_K}\in\mathcal W.
 \]
 If $\mathcal W\not\models \bq\preceq \bu$, then this identity fails in
 $A_K$ under the assignment sending each letter to its own class.
\end{lem}

\begin{proof}
 It is straightforward to verify that the relation is an equivalence
 relation. The inclusion
 $\bs(\bu+\br)\bt\subseteq K$ is equivalent to the conjunction of
 $\bs \bu\bt\subseteq K$ and $\bs \br\bt\subseteq K$, which proves compatibility
 with addition. If $\br=\sum_i \bt_i$ is a term, then
 \[
 \bs \bu\br\bt\subseteq K
 \quad\Longleftrightarrow\quad
 \bs \bu(\bt_i\bt)\subseteq K\text{ for every }i.
 \]
 This proves compatibility with right multiplication; compatibility with left
 multiplication follows by the same argument. Hence $\equiv_K$ is a congruence.

 Let $\mathbf{u}_1\approx \mathbf{v}_1$ be any identity of $\mathcal W$. Substitute arbitrary members
 of $\mathscr P(X)$ for its variables and denote the resulting expanded
 terms by $\mathbf{u}_1',\mathbf{v}_1'$. For every $\bs,\bt\in X^*$ we have
 $\mathcal W\models\bs \mathbf{u}_1'\bt\approx\bs \mathbf{v}_1'\bt$. Consequently,
 \[
 \begin{split}
 \bs \mathbf{u}_1'\bt\subseteq K
 &\Longleftrightarrow\mathcal W\models\bs \mathbf{u}_1'\bt\preceq \bu\\
 &\Longleftrightarrow\mathcal W\models\bs \mathbf{v}_1'\bt\preceq \bu
 \Longleftrightarrow\bs \mathbf{v}_1'\bt\subseteq K.
 \end{split}
 \]
 Thus $\mathbf{u}_1'\equiv_K\mathbf{v}_1'$. The preceding argument applies to arbitrary term substitutions. Every element of the
 quotient has a term representative, so the quotient satisfies all
 identities of $\mathcal W$. By Birkhoff's theorem, the quotient belongs to $\mathcal W$.

 Every summand of $\bu$ belongs to $K$, whereas $\bq\notin K$ when
 $\mathcal W\not\models \bq\preceq \bu$. Taking $\bs=\bt=\varepsilon$ shows
 $\bu\not\equiv_K\bu+\bq$, proving the final assertion.
\end{proof}

To verify an identity in $A_K$ it therefore suffices to check that,
after every term substitution and after adjoining arbitrary left
and right words, all summands on one side belong to $K$ if and only if all
summands on the other side belong to $K$. For $\bq\preceq \bu$, only the implication
from the summands of $\bu$ to those of $\bq$ needs to be proved.

\begin{lem}[A consequence of finite equational derivations]
\label{lem:hfb-compactness}
 Let $\mathcal T$ be finitely based, and let
 $\mathcal U_1,\ldots,\mathcal U_m$ be hereditarily finitely based varieties
 of the same signature. They need not be contained in $\mathcal T$.
 If the map
 \[
 \mathcal W\longmapsto
 (\mathcal W\cap\mathcal U_1,\ldots,\mathcal W\cap\mathcal U_m)
 \]
 is injective on the subvarieties of $\mathcal T$, then $\mathcal T$ is
 hereditarily finitely based. If injectivity is assumed only on an interval
 $[\mathcal A,\mathcal T]$, then every variety in that interval is finitely
 based.
\end{lem}

\begin{proof}
 Write $\mathcal W=\mathcal T\cap\operatorname{Mod}(\Sigma)$, and let
 $\Theta$ be a finite equational basis for $\mathcal T$. For each $i$, choose finite
 equational bases $\Gamma_i$ for $\mathcal W\cap\mathcal U_i$ and $\Delta_i$ for
 $\mathcal U_i$. Every identity in $\Gamma_i$ follows from
 $\Theta\cup\Sigma\cup\Delta_i$. An equational derivation has only
 finitely many steps. Since the union of the sets $\Gamma_i$ is finite,
 all these derivations together use only a finite subset
 $\Sigma_0\subseteq\Sigma$.

 Put $\mathcal W_0=\mathcal T\cap\operatorname{Mod}(\Sigma_0)$.
 We have $\mathcal W\subseteq\mathcal W_0$, and each
 $\mathcal W_0\cap\mathcal U_i$ satisfies $\Gamma_i$. Hence
 \[
 \mathcal W\cap\mathcal U_i
 =\mathcal W_0\cap\mathcal U_i\qquad(1\leq i\leq m).
 \]
 Injectivity implies $\mathcal W=\mathcal W_0$, which has the finite equational basis
 $\Theta\cup\Sigma_0$. For the assertion concerning an interval, if
 $\mathcal A\subseteq\mathcal W$, then also
 $\mathcal A\subseteq\mathcal W_0$; thus both varieties remain in the
 stated interval and the same argument applies.
\end{proof}

\subsection{A decomposition of every subvariety of $\mathcal C$}

The word identities imply, for every nonempty word $\bw$ and possibly empty
word $\bs$,
\[
 \bs \bw^2\approx(\bs \bw)^2;
\]
if $\bt\ne\varepsilon$, they also imply $\bs \bw^2\bt\approx\bs \bw\bt$.
Lemma~\ref{lem:C-closing} yields the identity
\begin{equation}\label{eq:C-LF}
 px^2\preceq px+py^2.
\end{equation}
Indeed, the right-hand side has a summand $py^2$ whose last letter occurs
earlier, as well as the summand $px$ of length two. The lemma gives
$(px)^2\approx px^2\preceq px+py^2$. As an identity, \eqref{eq:C-LF} permits
arbitrary term substitutions.

\begin{pro}\label{prop:C-intersections}
 For every $\mathcal W\subseteq\mathcal C$,
 \begin{equation}\label{eq:C-decomposition}
 \mathcal W=(\mathcal W\cap\mathcal D)
       \vee(\mathcal W\cap\mathcal J_0).
 \end{equation}
\end{pro}

\begin{proof}
 Choose an identity $\bq\preceq \bu$ that fails in $\mathcal W$. Let $X$ contain
 all letters occurring in it, and construct $K$ and $A_K$ as in
 Lemma~\ref{lem:syntactic-separation}.

 First suppose that $K$ contains a word $\mathbf{z}$ with $t(\mathbf{z})\in e(\mathbf{z})$.
 If $\mathbf{v}_1\in K$ and $\ell(\mathbf{v}_1)\geq2$, the identities $\mathbf{z}\preceq \bu$ and $\mathbf{v}_1\preceq \bu$
 hold in $\mathcal W$; Lemma~\ref{lem:C-closing} therefore gives
 $\mathbf{v}_1^2\in K$. We claim that $\bw\equiv_K\bw^2$ for every $\ell(\bw)\geq2$.
 If the word $\bt$ on the right is nonempty, the words
 $\bs \bw\bt$ and $\bs \bw^2\bt$ are already equal in $\mathcal C$.
 For $\bt=\varepsilon$, the preceding observation and the identity
 $(\bs \bw)^2\approx\bs \bw^2$ give
 \[
 \bs \bw\in K\quad\Longrightarrow\quad
 (\bs \bw)^2\in K\quad\Longrightarrow\quad\bs \bw^2\in K.
 \]
 The reverse implication follows from $\bs \bw\preceq\bs \bw^2$.
 Thus the claim holds for every pair of words on the left and right.
 Every product in $A_K$ is a finite sum of classes
 of words of length at least two. Since squaring preserves addition,
 $A_K\models xy\approx(xy)^2$. Thus
 $A_K\in\mathcal W\cap\mathcal D$.

 Now suppose that no word $\mathbf{z}\in K$ satisfies $t(\mathbf{z})\in e(\mathbf{z})$.
 To prove $A_K\models x^2\preceq x+y^2$, take nonempty terms
 $\mathbf{a},\mathbf{b}$ and surrounding words $\bs,\bt$, and suppose every summand of
 $\bs(\mathbf{a}+\mathbf{b}^2)\bt$ belongs to $K$. If $\bt\ne\varepsilon$, the identity
 $\mathbf{a}^2\bt\approx \mathbf{a}\bt$ proves that every summand of $\bs \mathbf{a}^2\bt$ belongs to $K$.
 If $\bt=\varepsilon$ and $\bs\ne\varepsilon$, apply \eqref{eq:C-LF}
 with $p=\bs,x=\mathbf{a},y=\mathbf{b}$ to obtain
 \[
 \mathcal W\models\bs \mathbf{a}^2\preceq\bs \mathbf{a}+\bs \mathbf{b}^2\preceq \bu.
 \]
 If both surrounding words are empty, the assumption is impossible:
 for any summand $\mathbf{b}_1$ of $\mathbf{b}$, the summand $\mathbf{b}_1^2$ of $\mathbf{b}^2$ would belong
 to $K$, although its last letter occurs earlier. These cases prove the
 required identity, and hence
 $A_K\in\mathcal W\cap\mathcal J_0$.

 In either case, $\bq\preceq \bu$ fails in a member of one of the two
 intersections on the right of \eqref{eq:C-decomposition}. Every term identity
 is equivalent to the finitely many inequalities requiring each summand
 on either side to be at most the other side. Consequently, the join on the right-hand side
 fails every identity that fails in $\mathcal W$. The reverse inclusion
 of varieties is immediate, so the two varieties are equal.
\end{proof}

The map $\mathcal W\mapsto
(\mathcal W\cap\mathcal D,\mathcal W\cap\mathcal J_0)$ is injective
by \eqref{eq:C-decomposition}. Theorem~\ref{thm:C-fb},
\eqref{eq:D-HFB-conclusion}, Theorem~\ref{thm:j0-hfb}, and
Lemma~\ref{lem:hfb-compactness} consequently give
\[
 \mathcal C\text{ is hereditarily finitely based}.
\]
Here we use the decomposition for \emph{every} subvariety; no general
result about joins of hereditarily finitely based varieties is required.

\subsection{Identities with $E(\bu)\cap T(\bu)=\varnothing$}

The condition $E(\bu)\cap T(\bu)=\varnothing$ concerns all summands of $\bu$
together. The condition $t(\bw)\notin e(\bw)$ for each individual summand
$\bw$ is weaker and is not sufficient for the following arguments.

\begin{lem}\label{lem:open-coincidence}
 If $E(\bu)\cap T(\bu)=\varnothing$, then
 \[
 \mathcal C\models \bq\preceq \bu
 \quad\Longleftrightarrow\quad
 \mathcal V\models \bq\preceq \bu.
 \]
\end{lem}

\begin{proof}
 Only the forward implication requires proof. If $\bq=x$, both criteria
 require a summand $x^k$, $k\geq1$, in $\bu$. Suppose $\ell(\bq)\geq2$.
 The $M_0$ condition in the criterion for $\mathcal C$ gives
 \[
 e(\bq)\subseteq E(\bu),\qquad t(\bq)\in T(\bu).
 \]
 In particular, $t(\bq)\notin e(\bq)$. The other conditions for $\mathcal C$
 give the required condition on contents and, at each first occurrence of
 a letter $b$, a summand $\mathbf{v}_0b\br$ such that every letter of $\mathbf{v}_0$ precedes $b$ in
 $\bq$. For a nonfinal letter $b$ of $\bq$, we have $b\in E(\bu)$, so
 $b\notin T(\bu)$. Therefore $\br\ne\varepsilon$. This is the stronger
 condition needed for $S$, and all conditions of
 Lemma~\ref{lem:S-criterion} hold.
\end{proof}

\begin{lem}\label{lem:open-meetC}
 If $\mathcal A\subseteq\mathcal W\subseteq\mathcal V$ and
 $E(\bu)\cap T(\bu)=\varnothing$, then
 \[
 \mathcal W\models \bq\preceq \bu
 \quad\Longleftrightarrow\quad
 \mathcal W\cap\mathcal C\models \bq\preceq \bu.
 \]
\end{lem}

\begin{proof}
 We first record restrictions imposed by $A_4\models \bq\preceq \bu$.
 If a letter of $\bq$ is absent from $\bu$, assigning it $0$ and assigning
 every other letter $3$ gives values $0$ and $3$ for $\bq$ and $\bu$,
 respectively, a contradiction. Hence $c(\bq)\subseteq c(\bu)$.
 Next assign $3$ to every letter of $E(\bu)$ and $1$ to every letter of
 $T(\bu)$. Every summand of $\bu$ has value $1$. Since all letters of $\bq$
 occur in $\bu$, its value is at most $1$ if and only if
 \[
 e(\bq)\subseteq E(\bu),\qquad t(\bq)\in T(\bu).
 \]
 Indeed, under this assignment a word has value $1$ if and only if all
 its nonfinal letters have value $3$ and its last letter has value $1$.

 Suppose $\bq\preceq \bu$ fails in $\mathcal W$, and form
 $A_K\in\mathcal W$ by Lemma~\ref{lem:syntactic-separation}.
 Since $A_4\in\mathcal W$, every $\bw\in K$ satisfies the two displayed
 restrictions. Consequently every term whose summands belong to
 $K$ has disjoint sets of nonfinal and final letters.
 Take any identity of $\mathcal C$, any term substitution, and any
 surrounding words. Suppose all summands on one resulting side $\bv$
 belong to $K$, and let $\mathbf{v}_1$ be a summand on the other side. We have
 $\mathcal C\models \mathbf{v}_1\preceq \bv$. Since $E(\bv)\cap T(\bv)=\varnothing$,
 Lemma~\ref{lem:open-coincidence} yields
 $\mathcal V\models \mathbf{v}_1\preceq \bv$. Therefore
 \[
 \mathcal W\models \mathbf{v}_1\preceq \bv\preceq \bu,
 \]
 and $\mathbf{v}_1\in K$. The same argument applies in the reverse direction, so
 $A_K$ satisfies every identity of $\mathcal C$.
 Thus $A_K\in\mathcal W\cap\mathcal C$, and $\bu$ and $\bu+\bq$ have distinct values in this algebra. This proves one implication. The converse follows from the inclusion
 $\mathcal W\cap\mathcal C\subseteq\mathcal W$.
\end{proof}

\subsection{Three finitely based subvarieties of $\mathcal E$}

\begin{lem}\label{lem:A4-strong-M1}
 If $\ell(\bq)\geq2$, $A_4\models \bq\preceq \bu$, and $\mathrm{SA}(\bu,\bq)$ holds,
 then $M_1\models \bq\preceq \bu$. Consequently, for $\ell(\bq)\geq2$,
 \[
 \Var(A_4,T)\models \bq\preceq \bu
 \quad\Longleftrightarrow\quad
 A_4\models \bq\preceq \bu\ \text{and}\ \mathrm{SA}(\bu,\bq).
 \]
\end{lem}

\begin{proof}
 The criterion for $A_4$ already provides a summand with content
 contained in $c(\bq)$ and all the required first-occurrence conditions.
 The condition at each nonfinal position is stronger than the corresponding
 condition for $M_1$. It remains to verify the condition for $M_1$ at a last letter $t$ that
 has not occurred earlier.
 Put $U=e(\bq)$ in this case.

 If the first alternative of $\mathrm{SA}(\bu,\bq)$ holds, choose a summand
 $\mathbf{u}_1$ with $\ell(\mathbf{u}_1)\geq2$ and $c(\mathbf{u}_1)\subseteq c(\bq)$. If $t\notin c(\mathbf{u}_1)$,
 then $c(\mathbf{u}_1)\subseteq U$. Otherwise the first letter of $\mathbf{u}_1$ outside
 $U$ is $t$, and $\mathbf{u}_1$ has the length required in the $M_1$ criterion.
 If the second alternative holds, choose the summands $a,a\br$ stated
 there. For $a\ne t$, the content of the summand $a$ is contained in $U$.
 For $a=t$, the summand $t\br$ has length at least two and starts with
 $t$. In every case the remaining $M_1$ condition is satisfied.
 The final equivalence follows from the criterion for $T$.
\end{proof}

\begin{lem}\label{lem:E-conflicting}
 If $E(\bu)\cap T(\bu)\ne\varnothing$ and $\ell(\bq)\geq2$, then
 \[
 \mathcal E\models \bq\preceq \bu
 \quad\Longleftrightarrow\quad
 A_4\models \bq\preceq \bu\ \text{and}\ \mathrm{SA}(\bu,\bq).
 \]
 For a word $\bq=x$ of length one, the condition is the existence of a
 summand $x^k$, $k\geq1$, in $\bu$.
\end{lem}

\begin{proof}
 Necessity follows from $A_4,T\in\mathcal E$. We prove sufficiency.
 If only the second alternative of $\mathrm{SA}(\bu,\bq)$ holds, the summands
 $a,a\br$ and the valid identity
 $a^2\preceq a+a\br$ give $a^2\preceq \bu$. We may add $a^2$ to $\bu$ without
 changing its value. Thus, after adding to $\bu$ a summand already at most
 $\bu$, if necessary, there is a word
 \[
 \mathbf{a}=\mathbf{a}_0a,\qquad \mathbf{a}_0\ne\varepsilon,\qquad c(\mathbf{a})\subseteq c(\bq),\qquad \mathbf{a}\preceq \bu.
 \]
 Since $\bu+\mathbf{a}\approx \bu$ has been derived, we may apply the criterion with
 right-hand side $\bu+\mathbf{a}$.
 If $c(\bq)$ has one letter, the content of $\mathbf{a}$ is already $c(\bq)$,
 and the criterion for $S$ gives $\bq\preceq \bu$.

 Suppose $|c(\bq)|\geq2$, and put $t=t(\bq)$. If $a=t$, let $\mathbf{b}$ be the
 word listing, in order of first appearance, the letters occurring
 before the last position of $\bq$. If $a\ne t$, let $\mathbf{b}$ list all letters
 of $\bq$ in their order of first appearance. In both cases $\mathbf{b}$ is nonempty
 and $c(\mathbf{a}_0\mathbf{b}a)=c(\bq)$.
 Let $\mathbf{z}$ be obtained by concatenating the proper prefixes of all
 original summands of $\bu$, and then their final letters. Since
 $E(\bu)\cap T(\bu)\ne\varnothing$, at least one of these prefixes is
 nonempty. In particular, $\mathbf{z}$ is nonempty. Consider
 \[
 \mathbf{w}_0=\mathbf{a}_0\mathbf{b}\mathbf{z}a.
 \]
 We verify the criterion for $S$ with right-hand side $\bu+\mathbf{a}$.
 The content of $\mathbf{w}_0$ contains $c(\bu)$; conversely $c(\bq)\subseteq c(\bu)$
 follows from the $A_4$ identity. Thus $c(\mathbf{w}_0)=c(\bu)$, and all original
 summands together give the required condition on contents.

 For a first occurrence in $\mathbf{a}_0$, the summand $\mathbf{a}_0a$ has that same
 occurrence before its last position. For a first occurrence in $\mathbf{b}$
 that comes from a nonfinal first occurrence in $\bq$, all letters that
 preceded it in $\bq$ have already appeared in $\mathbf{a}_0\mathbf{b}$. The corresponding
 condition from the $A_4$ criterion therefore still applies: either the
 content of a summand is contained in the set of preceding letters, or an original summand
 has the current letter as its first outside letter, at a nonfinal
 position. The only additional case is a last letter $t$ of $\bq$ that
 had not appeared earlier and for which $a\ne t$. If $t\in c(\mathbf{a})$,
 its first occurrence is already handled in $\mathbf{a}_0$. Otherwise all letters
 of $\mathbf{a}$ have occurred before this $t$, so the content of $\mathbf{a}$ is contained
 in the set of preceding letters. The $A_4$ criterion also gives the required condition on the first
 letter outside the preceding set.

 Consider now the part of $\mathbf{z}$ consisting of proper prefixes. At the first occurrence of each
 letter, its corresponding original summand has all
 earlier letters already present and has a nonempty suffix after
 this letter. After these prefixes, every letter of
 $E(\bu)$ has appeared. Choose $b\in E(\bu)\cap T(\bu)$ and an original
 summand ending in $b$. Its content is contained in $E(\bu)$.
 Hence every letter in the content of this summand has appeared before any new letter in the
 remaining final-letter part of $\mathbf{z}$. Each such new letter also has an
 original summand ending in it, all of whose earlier letters lie in
 $E(\bu)$. These observations verify both first-occurrence requirements.
 The final letter $a$ has already occurred. Thus $S\models \mathbf{w}_0\preceq \bu+\mathbf{a}\approx \bu$.

 In the defining identity $xy\preceq x+xyz$ of $\mathcal E$, substitute
 $x=\mathbf{a}_0,y=\mathbf{b},z=\mathbf{z}$ and multiply on the right by $a$. All substituted
 words are nonempty, and we obtain
 \[
 \mathbf{a}_0\mathbf{b}a\preceq \mathbf{a}_0a+\mathbf{a}_0\mathbf{b}\mathbf{z}a\preceq \bu.
 \]
 We may now add $\mathbf{a}_0\mathbf{b}a$ to $\bu$. Its content is exactly $c(\bq)$,
 and the original first-occurrence conditions from $A_4$ remain valid.
 The criterion for $S$ consequently gives $\bq\preceq \bu$.
 The one-letter assertion follows from the criterion for $A_4$ for
 necessity and from $x\preceq x^2$ for sufficiency.
\end{proof}

We now introduce the identities
\[
 I:\ xyz\preceq xz+yt,\qquad
 J:\ xy\preceq y+xyz,\qquad
 L_0:\ yx\preceq x+z^2+yzx,
\]
and the varieties
\[
 \mathcal H=\Var(A_4,T),\qquad
 \mathcal H_1=\mathcal H\cap\operatorname{Mod}(J),\qquad
 \mathcal H_2=\mathcal H_1\cap\operatorname{Mod}(L_0).
\]

\begin{pro}\label{prop:H-hosts}
 We have
 \[
 \mathcal H=\mathcal E\cap\operatorname{Mod}(I),\qquad
 \mathcal H_1=\Var(A_4,M_1).
 \]
 In particular, $\mathcal H,\mathcal H_1,\mathcal H_2$ are finitely
 based and contain $\mathcal A$.
\end{pro}

\begin{proof}
 The identity $I$ holds in $A_4$ and $T$. Indeed, the summand $xz$
 has length two and its content is contained in $c(xyz)$. For the
 first occurrences of $x,y,z$ in the word, use respectively the
 summands $xz,yt,xz$; the first two occurrences are nonfinal in those
 summands. Thus both identity criteria hold.

 Conversely, consider an identity $\bq\preceq \bu$ of $\Var(A_4,T)$.
 If $\bq$ has length one, the assertion follows from the common criterion
 requiring a summand that is a power of that letter. Suppose $\ell(\bq)\geq2$. If $E(\bu)\cap T(\bu)\ne\varnothing$, apply
 Lemma~\ref{lem:E-conflicting}. Suppose instead that
 $E(\bu)\cap T(\bu)=\varnothing$. The $A_4$ restrictions in the proof of
 Lemma~\ref{lem:open-meetC} show that $\bq=\mathbf{u}_0t$ with $\mathbf{u}_0\ne\varepsilon$,
 $c(\mathbf{u}_0)\subseteq E(\bu)$, and $t\in T(\bu)$. The second alternative of
 $\mathrm{SA}(\bu,\bq)$ cannot occur: a letter that is both a one-letter
 summand and the first letter of a longer summand would belong to
 $E(\bu)\cap T(\bu)$. Hence there is a summand $\br{}t$, with
 $\br\ne\varepsilon$ and $c(\br)\subseteq c(\mathbf{u}_0)$.

 Consider the letters of $\mathbf{u}_0$ in their order of first occurrence.
 If a letter $b$ does not yet occur in the prefix of the auxiliary word
 ending in $t$, choose an original summand $\mathbf{v}_0b\mathbf{z}$ given by the $A_4$
 criterion. Its prefix letters lie among the preceding letters of
 $\bq$. Moreover, $\mathbf{z}\ne\varepsilon$, since $b\in E(\bu)$ and
 $E(\bu)\cap T(\bu)=\varnothing$. From $I$ we obtain
 \[
 \br\mathbf{v}_0bt\preceq \br{}t+\mathbf{v}_0b\mathbf{z}\preceq \bu.
 \]
 This replaces the auxiliary word by another word at most $\bu$, still
 having its content contained in $c(\bq)$, and containing $b$. Repeating this
 argument gives a word with content exactly $c(\bq)$. Adding it
 to $\bu$ gives the required condition on contents in the criterion for $S$;
 the original first-occurrence conditions are retained. This proves
 $\bq\preceq \bu$ in $\mathcal E\cap\operatorname{Mod}(I)$, and hence the first
 equality of varieties.

 The identity $J$ holds in $A_4,M_1$, and these algebras belong to
 $\mathcal H$. Consider a common identity $\bq\preceq \bu$ of these two
 algebras with $\ell(\bq)\geq2$. If $\mathrm{SA}(\bu,\bq)$ holds, the criterion for
 $T$ already gives the identity in $\mathcal H$. Otherwise every
 summand with content contained in $c(\bq)$ is a one-letter summand,
 and none of these letters begins an original summand of length at
 least two. Such one-letter summands exist by the criterion for $M_1$.
 Choose the earliest of these letters, say $a$, in the first-occurrence
 order of $\bq$. No original summand has its content contained in the set of
 letters preceding $a$: such a summand would have to be an earlier
 one-letter summand of the same kind. The $M_1$ criterion therefore
 supplies an original summand $\mathbf{v}_0a\mathbf{z}$ of length at least two, whose first
 letter outside that preceding set is $a$.
 We have $\mathbf{v}_0\ne\varepsilon$, since otherwise $a$ begins a longer
 summand; also $\mathbf{z}\ne\varepsilon$, since otherwise $\mathbf{v}_0a$ itself is a
 longer summand with content contained in $c(\bq)$. Applying $J$ gives
 \[
 \mathbf{v}_0a\preceq a+\mathbf{v}_0a\mathbf{z}\preceq \bu.
 \]
 Add $\mathbf{v}_0a$ to $\bu$. The first alternative of $\mathrm{SA}$ now holds,
 and the already valid identities in $A_4,M_1$ remain valid. The
 criterion for $\mathcal H$ gives the desired conclusion. The assertion for words of length one again follows from the criterion
 requiring a summand that is a power of that letter. This proves the
 second equality.

 Finally, $I,J,L_0$ all hold in $A_4$. For $L_0$, the summand $x$
 has content contained in $c(yx)$, $yzx$ provides the nonfinal first
 occurrence of $y$, and the summand $x$ provides the last-letter
 occurrence of $x$. Thus the three varieties contain $\mathcal A$.
 Their finite equational bases are obtained by adding successively $I,J,L_0$
 to a finite equational basis for $\mathcal E$.
\end{proof}

\subsection{Proper subvarieties above $\mathcal A$}

\begin{lem}\label{lem:H-first-cover}
 If $\mathcal A\subseteq\mathcal W\subsetneq\mathcal H$, then
 $\mathcal W\subseteq\mathcal H_1$.
\end{lem}

\begin{proof}
 Since $\mathcal H=\mathcal A\vee\Var(T)$, we have
 $T\notin\mathcal W$. Choose an identity of $\mathcal W$ and an
 assignment in $T$ under which it fails. Represent the assigned elements
 by nonempty terms on $x,y,z$. The definition of the congruence
 defining $T$ supplies left and right words that distinguish the two
 resulting term classes by membership of all summands in
 $\mathscr L$. Substituting the representatives and adjoining these
 words, and interchanging the two sides if necessary, gives an identity
 $\bu'\approx \bv'$ of $\mathcal W$ such that
 \[
 \bu'\subseteq\mathscr L,\qquad \bv'\nsubseteq\mathscr L.
 \]

 Put $\br=y+z+xyz$. By the normal form for words, the words of $\mathscr L$
 reduce to
 \[
 y,\quad z,\quad xyz,\quad xyz^2.
 \]
 The first three are summands of $\br$. For the last, the criterion for
 $S$ holds because the content of the summand $y$ is contained in the set of letters
 preceding the first $z$. Thus $\mathcal V\models \bu'\preceq \br$, and
 $\mathcal W\models \bv'\preceq \br$. Since $A_4\in\mathcal W$, every summand
 of $\bv'$ satisfies $\bw\preceq \br$ in $A_4$.

 We now determine all such words. The only words of length one
 are $y,z$. A word of length at least two cannot start with $y$ or $z$, since the only
 summands of $\br$ starting with these letters have length one. It must
 therefore start with $x$, and must contain $y$ or $z$ to have any
 summand of $\br$ whose content is contained in its own.
 If the order of first occurrences is $xz$, the only allowed normal
 form is $xz$. The alternative $xz^2$ would require the first $z$ to
 occur at a nonfinal position in an appropriate original summand,
 whereas no summand lies wholly in $\{x\}$ and only the one-letter
 summand $z$ has $z$ as its first letter outside $\{x\}$.
 For the same reason the order $xzy$ is impossible. For the orders $xy$
 and $xyz$, both normal forms are possible, with a new last letter or with a
 repeated last letter. Hence the complete list is
 \[
 y,\quad z,\quad xy,\quad xy^2,\quad xz,\quad xyz,\quad xyz^2.
 \]
 Among these, a word not in $\mathscr L$ has normal form
 $xy,xy^2$, or $xz$. Thus $\mathcal W$ satisfies at least one of
 $xy\preceq \br$, $xy^2\preceq \br$, and $xz\preceq \br$.

 From $xy^2\preceq \br$ and $xy\preceq xy^2$ we obtain $xy\preceq \br$.
 Substituting $z\mapsto xyz$ in $xy\preceq y+z+xyz$ and deleting repeated
 nonfinal occurrences gives $xy\preceq y+xyz$, namely $J$.
 If instead $xz\preceq \br$, substitute $y\mapsto xyz$ and $z\mapsto y$.
 This gives
 \[
 xy\preceq y+xyz+xyzy.
 \]
 The identity $xyzy\preceq xyz+y$, valid in $S$, again yields $J$.
 Therefore $\mathcal W\subseteq\mathcal H_1$.
\end{proof}

\begin{lem}\label{lem:H-second-cover}
 If $\mathcal A\subseteq\mathcal W\subsetneq\mathcal H_1$, then
 $\mathcal W\subseteq\mathcal H_2$.
\end{lem}

\begin{proof}
 Since $\mathcal H_1=\mathcal A\vee\Var(M_1)$, the algebra
 $M_1$ does not belong to $\mathcal W$. In the natural order of $M_1$,
 the set $I_1=\{1,4\}$ is the principal down-set of $1$. It is closed
 under addition, and a finite sum belongs to $I_1$ exactly when every
 summand belongs to $I_1$. The following table distinguishes all four elements, with $1$ indicating
 membership in $I_1$ and $0$ indicating nonmembership:
 \[
 \begin{array}{c|ccc}
 u&u\in I_1&u2\in I_1&u4\in I_1\\\hline
 0&0&0&0\\
 1&1&0&0\\
 2&0&0&1\\
 4&1&1&1
 \end{array}
 \]
 Under $(x,y,z)=(1,2,4)$, the elements $0,1,2,4$ are represented by
 $x^2,x,y,z$, respectively. Choose an identity of $\mathcal W$
 that fails in $M_1$. Use these representatives and, as given by the
 table, either no additional factor or a right factor $y$ or $z$.
 Reversing the two sides if needed, we obtain an identity $\bu'\approx \bv'$ of
 $\mathcal W$ for which every summand of $\bu'$ has value in $I_1$,
 while at least one summand of $\bv'$ has value outside $I_1$.

 A word has value in $I_1$ under this assignment precisely when it is
 the one-letter word $x$, or its first letter different from $y$ is
 $z$. Indeed, any initial string of $y$'s has value $2$. If the next letter is $z$, the value is $4$ and remains $4$ after any
 further multiplication. If the next letter is $x$ and the word has length
 at least two, its value is $0$. A word consisting only of $y$'s
 has value $2$.

 Put $\br=x+z^2+yzx$. Every word just described satisfies $\bw\preceq \br$ in
 $\mathcal H_1$. This is immediate for $\bw=x$. Every other such word
 contains $z$, so the content of $z^2$ is contained in $c(\bw)$. The summand $z^2$
 supplies all required conditions at a first $z$; the summand $yzx$
 does the same at a first $y$. Before any first $x$ in such a word,
 the letter $z$ has already appeared, so the content of the summand $z^2$
 is contained in the set of preceding letters. These facts verify both the $A_4$
 and $M_1$ criteria. It follows that
 $\mathcal W\models \bv'\preceq \br$.

 It remains to determine words that have value outside $I_1$ but satisfy
 $\bw\preceq \br$ in $A_4$. A word of length at least two beginning with $x$
 fails the first-position condition, because the only summand of $\br$
 beginning with $x$ is $x$ itself. A pure power of $y$ has no original
 summand whose content is contained in its own. Words starting with
 $z$, or whose first letter different from $y$ is $z$, already have
 value in $I_1$. The only remaining orders of first occurrence are $yx$
 and $yxz$. If the first $x$ is nonfinal, the set of preceding letters is
 $\{y\}$, which contains the content of no summand of $\br$. Moreover, the only
 summand whose first letter outside $\{y\}$ is $x$ is the one-letter
 summand $x$. The $A_4$ criterion therefore fails. The only remaining
 normal form is $yx$.
 At least one summand of $\bv'$ consequently gives
 \[
 yx\preceq x+z^2+yzx,
 \]
 which is $L_0$. Hence $\mathcal W\subseteq\mathcal H_2$.
\end{proof}

\subsection{Completion of the proof for $\mathcal E$}

\begin{lem}\label{lem:H2-closed-agreement}
 If $E(\bu)\cap T(\bu)\ne\varnothing$, then
 \[
 A_4\models \bq\preceq \bu
 \quad\Longleftrightarrow\quad
 \mathcal H_2\models \bq\preceq \bu.
 \]
\end{lem}

\begin{proof}
 Only the forward implication requires proof. For a word of length one,
 the assertion follows from the criterion requiring a summand that is a
 power of that letter. Suppose $\ell(\bq)\geq2$. If
 $M_1\models \bq\preceq \bu$, the identity holds in
 $\mathcal H_1=\Var(A_4,M_1)$ and hence in $\mathcal H_2$.
 Otherwise the criteria show that failure in $M_1$ can occur only at
 a last letter $t$ that has not appeared earlier. Using the
 normal form for words, write $\bq=\mathbf{u}_0t$ with $\mathbf{u}_0\ne\varepsilon$ and
 $t\notin c(\mathbf{u}_0)$. There is no summand whose content is contained in
 $c(\mathbf{u}_0)$, and there is no summand of length at least two whose first
 letter outside $c(\mathbf{u}_0)$ is $t$.
 Therefore no summand of length at least two can have its content
 contained in $c(\bq)$: if it did not contain $t$, this would contradict the
 first condition; if it did, it would contradict the second.
 The only possible summand of length one with content contained in $c(\bq)$
 is $t$, and the $A_4$ criterion
 guarantees that it is indeed a summand of $\bu$.

 Let $\br$ be the concatenation of all proper prefixes of original
 summands, followed by all their final letters, and put $\mathbf{z}=\br^2$.
 We verify the two identities
 \[
 \mathcal V\models \mathbf{z}\preceq \bu,\qquad
 \mathcal V\models \mathbf{u}_0\mathbf{z}t\preceq \bu.
 \]
 Both words have content $c(\bu)$, so the condition on contents holds.
 In the proper-prefix part of $\br$, a first occurrence of
 a letter is also a nonfinal occurrence in its corresponding original
 summand, all of whose earlier letters have already appeared.
 After this part, every letter of $E(\bu)$ has appeared. Since
 $E(\bu)\cap T(\bu)\ne\varnothing$, there is an original summand ending
 in an element of $E(\bu)$, and its content is contained in
 $E(\bu)$. The content of this summand is contained in the set of preceding letters at every new
 occurrence in the remaining final-letter part. Each such occurrence
 also has an original summand ending in that letter. This verifies
 the criterion for $\mathbf{z}$. For $\mathbf{u}_0\mathbf{z}t$, the first occurrences in $\mathbf{u}_0$ satisfy
 the $A_4$ nonfinal-position conditions from the given identity.
 The same argument applies to the subsequent part $\mathbf{z}$, in which further
 letters may already have occurred. The final letter $t$ has already
 occurred. Thus the criterion for $S$ proves both identities.

 As $\mathbf{z}=\br^2$, the word identities give $\mathbf{z}^2\approx \mathbf{z}$.
 Substitute $x=t,y=\mathbf{u}_0,z=\mathbf{z}$ in $L_0$ to obtain
 \[
 \mathbf{u}_0t\preceq t+\mathbf{z}^2+\mathbf{u}_0\mathbf{z}t\preceq \bu.
 \]
 All substituted factors are nonempty. This proves the required
 identity in $\mathcal H_2$. The reverse implication follows from
 $A_4\in\mathcal H_2$.
\end{proof}

\begin{pro}\label{prop:E-transfer-HFB}
 If $\mathcal C$ is hereditarily finitely based, then so is $\mathcal E$.
\end{pro}

\begin{proof}
 We first prove that $\mathcal H_2$ is hereditarily finitely based.
 For $\mathcal A\subseteq\mathcal W\subseteq\mathcal H_2$, identities
 $\bq\preceq \bu$ with $E(\bu)\cap T(\bu)\ne\varnothing$ are determined by
 Lemma~\ref{lem:H2-closed-agreement}: they are precisely the corresponding
 identities of $A_4$. Those with $E(\bu)\cap T(\bu)=\varnothing$ are
 determined by $\mathcal W\cap\mathcal C$, by
 Lemma~\ref{lem:open-meetC}. Since arbitrary term identities are
 equivalent to finitely many inequalities with a word on the lower side, the map
 \[
 \mathcal W\longmapsto\mathcal W\cap\mathcal C
 \]
 is injective on $[\mathcal A,\mathcal H_2]$. The assertion concerning intervals in
 Lemma~\ref{lem:hfb-compactness} applies: $\mathcal H_2$ is finitely
 based and $\mathcal C$ is hereditarily finitely based. Thus every
 variety in this interval is finitely based. Any remaining subvariety
 of $\mathcal H_2$ does not contain $A_4$, and is therefore contained
 in $\mathcal C$ and finitely based as well.

 Every proper subvariety of $\mathcal H_1$ either contains $A_4$, in
 which case Lemma~\ref{lem:H-second-cover} gives its inclusion in $\mathcal H_2$,
 or does not contain $A_4$, in which case it is contained in
 $\mathcal C$. Together with the finite equational basis of $\mathcal H_1$
 itself, this proves that $\mathcal H_1$ is hereditarily finitely
 based. The same argument using Lemma~\ref{lem:H-first-cover} proves
 that $\mathcal H$ is hereditarily finitely based.

 We now show that, for every $\mathcal W\subseteq\mathcal E$,
 \begin{equation}\label{eq:E-decomposition}
 \mathcal W=(\mathcal W\cap\mathcal H)
        \vee(\mathcal W\cap\mathcal C).
 \end{equation}
 Choose an identity $\bq\preceq \bu$ that fails in $\mathcal W$, and construct
 $A_K\in\mathcal W$ as above. Define
 \[
 E(K)=\bigcup_{\bw\in K}e(\bw),\qquad T(K)=\{t(\bw):\bw\in K\}.
 \]
 If $E(K)\cap T(K)=\varnothing$, then every term with summands
 in $K$ has disjoint nonfinal and final letter sets. The second half
 of the proof of Lemma~\ref{lem:open-meetC}, using
 Lemma~\ref{lem:open-coincidence}, shows that every identity of
 $\mathcal C$ is respected by $\equiv_K$ under arbitrary term
 substitutions. Hence $A_K\in\mathcal C$. This argument does
 not require $A_4\in\mathcal W$: the required disjointness now follows
 directly from the assumption on $K$.

 If $E(K)\cap T(K)\ne\varnothing$, choose two words in $K$ (possibly
 the same word) such that a letter in this intersection occurs nonfinally in one word
 and finally in the other. Their sum
 is a term $\bv_0$ with summands in $K$ and
 $E(\bv_0)\cap T(\bv_0)\ne\varnothing$.
 Take any identity of $\mathcal H$, any term substitution, and
 any left and right surrounding words. Suppose all summands on one
 resulting side $\bv$ belong to $K$, and let $\mathbf{v}_1$ be a summand on the
 other side. Then
 \[
 \mathcal H\models \mathbf{v}_1\preceq \bv+\bv_0.
 \]
 The right-hand term has intersecting nonfinal and final letter
 sets. Lemmas~\ref{lem:A4-strong-M1} and~\ref{lem:E-conflicting} show
 that $\mathcal H$ and $\mathcal E$ have exactly the same inequalities
 with this right-hand side, including the one-letter case. Hence
 \[
 \mathcal E\models \mathbf{v}_1\preceq \bv+\bv_0,
 \qquad
 \mathcal W\models \mathbf{v}_1\preceq \bv+\bv_0\preceq \bu.
 \]
 Thus $\mathbf{v}_1\in K$. The argument applies in both directions of the
 identity, so $A_K$ satisfies all identities of $\mathcal H$.
 Therefore $A_K\in\mathcal H$.

 In either case, the original identity fails in an algebra belonging to one of
 the two intersections on the right of \eqref{eq:E-decomposition}. As in the
 proof of \eqref{eq:C-decomposition}, this proves equality of the varieties.
 The resulting map to the two intersections is injective. Since
 $\mathcal E$ is finitely based and $\mathcal C,\mathcal H$ are
 hereditarily finitely based, Lemma~\ref{lem:hfb-compactness} proves
 that $\mathcal E$ is hereditarily finitely based.
\end{proof}
\section{Failure of strong nonfinite basability}\label{sec:not-strong}

We construct a five-element finitely based ai-semiring $H_5$ such that
$S\in\Var(H_5)$.

\subsection{A finitely based ai-semiring of order five}

Let $H_5=\{0,1,2,3,4\}$ with the operations given in Table~\ref{tab:H5}.

\begin{table}[ht]
\caption{The Cayley tables of $H_5$}\label{tab:H5}
\begin{tabular}{c|ccccc}
+&0&1&2&3&4\\\hline
0&0&0&0&0&0\\
1&0&1&0&1&1\\
2&0&0&2&2&2\\
3&0&1&2&3&3\\
4&0&1&2&3&4
\end{tabular}
\qquad
\begin{tabular}{c|ccccc}
$\cdot$&0&1&2&3&4\\\hline
0&0&0&0&0&4\\
1&0&0&0&0&4\\
2&0&1&2&3&4\\
3&3&3&3&3&4\\
4&4&4&4&4&4
\end{tabular}
\end{table}

Although $H_5$ and $S$ use the same five labels, their operations differ.

\begin{lem}\label{snfb:lem:Haxioms}
The algebra $H_5$ is an ai-semiring.
\end{lem}
\begin{proof}
The addition is the join operation of the partial order
\[
 4<3<1<0,\qquad 4<3<2<0.
\]
The only incomparable pair is $\{1,2\}$ and its join is $0$. Thus addition is associative,
commutative and idempotent.

For multiplication, put $L_s(t)=st$. For $i\in\{0,3\}$, let $f_i$
fix $4$ and send every other element to $i$. The multiplication table gives
\[
 L_0=L_1=f_0,\qquad L_2=\operatorname{id}_{H_5},\qquad
 L_3=f_3,\qquad L_4(t)=4.
\]
Direct composition gives $L_{st}=L_s\circ L_t$. Indeed,
$f_i\circ f_j=f_i$ for $i,j\in\{0,3\}$, every displayed map fixes
$4$, and composition with the constant map $L_4$ is constant $4$ on
either side. The cases involving $L_2$ are immediate. These are exactly
the cases in the multiplication table, including $L_1=L_0$.
Therefore $(st)u=s(tu)$.

Each left translation preserves joins: this is immediate for the identity
and constant maps, and for $f_0,f_3$ it follows because a join is $4$
exactly when both its arguments are $4$. This proves left distributivity.

For right distributivity, $R_4:s\mapsto s4$ is constant $4$. If
$t\in\{0,1,2,3\}$, the map $R_t:s\mapsto st$ has values
\[
 (R_t(4),R_t(3),R_t(1),R_t(2),R_t(0))=(4,3,0,t,0).
\]
It is order preserving on both chains displayed above, since $3\leq_{H_5}t$
for every such $t$. For the only incomparable pair,
\[
 R_t(1+2)=R_t(0)=0=0+t=R_t(1)+R_t(2).
\]
Thus every right translation preserves joins, as required.
\end{proof}

We use the notation for words and sums of words introduced in
Section~\ref{sec:nfb}.

\begin{lem}\label{snfb:lem:evaluation}
The product $s_1\cdots s_m$ in $H_5$, where $s_i\in H_5$, is
determined as follows.
\begin{enumerate}[(i)]
\item If some $s_i=4$, the value is $4$.
\item If every $s_i=2$, the value is $2$.
\item Otherwise, assume that no $s_i=4$ and let $j$ be the first position
with $s_j\ne2$. If $s_j=0$ or $s_j=3$, the value is $s_j$. If $s_j=1$,
the value is $1$ when $j=m$, and is $0$ when $j<m$.
\end{enumerate}
\end{lem}
\begin{proof}
The element $4$ is a two-sided multiplicative zero. The element $2$ is
a left identity. On $H_5\setminus\{4\}$, both $0$ and $3$ are left zero
elements, while $1$ multiplied by any element of this set gives $0$.
These facts, read directly from Table~\ref{tab:H5}, give all three assertions.
\end{proof}

\subsection{An identity criterion for $H_5$}

Let $\bu$ be a term and let $\bq$ be a word. Define
\[
 D_{\bq}(\bu)=\{\bp \in \bu:c(\bp)\subseteq c(\bq)\}.
\]
When this set is nonempty, we also regard it as a term.
For a variable $b\in c(\bq)$, let $U_b$ be the set of variables
occurring before the \emph{first} occurrence of $b$ in $\bq$.
In a decomposition $\bp=\br b\bs$, the contexts $\br,\bs$ may be empty. Empty
contexts mean omitted factors; variables in identities are never
replaced by the empty word.

\begin{lem}\label{snfb:lem:criterion}
The inequality $\bq\preceq\bu$ holds in $H_5$ if and only if the
following two conditions hold for every $b\in c(\bq)$.
\begin{enumerate}[(1)]
\item\label{snfb:cond:ordinary} There is a word $\bp=\br b\bs \in D_{\bq}(\bu)$ with
$c(\br)\subseteq U_b$.
\item\label{snfb:cond:nonfinal} If the first occurrence of $b$ is not the
last position of $\bq$, then either there is such a word $\bp=\br b\bs \in D_{\bq}(\bu)$
with $\bs$ nonempty, or there is a word $\bt \in D_{\bq}(\bu)$ with
$c(\bt)\subseteq U_b$.
\end{enumerate}
In particular, when the inequality holds,
\begin{equation}\label{snfb:eq:coverage}
 D_{\bq}(\bu)\ne\varnothing,\qquad
 \bigcup_{\bp \in D_{\bq}(\bu)}c(\bp)=c(\bq).
\end{equation}
\end{lem}
\begin{proof}
First suppose the inequality holds. Assign $4$ to all variables outside
$c(\bq)$; every summand outside $D_{\bq}(\bu)$ then has value $4$, the least
element of the additive order, and does not affect a sum with other
summands.

If condition~\ref{snfb:cond:ordinary} fails at $b$, assign $2$ to $U_b$,
$0$ to $b$, and $3$ to the other variables of $\bq$. The value of $\bq$
is $0$. A word in $D_{\bq}(\bu)$ whose content is contained in $U_b$ has value
$2$. Any other word in $D_{\bq}(\bu)$ has first variable outside $U_b$
different from $b$, since the required decomposition does not exist;
its value is therefore $3$. The right side has value $4$, $3$ or $2$,
none of which is at least $0$. This is a contradiction.

Next suppose condition~\ref{snfb:cond:nonfinal} fails at a first occurrence
of $b$ that is not final. Assign $2$ to $U_b$, $1$ to $b$, and $3$ to
the other variables of $\bq$. Now $\bq$ has value $0$. Every word in $D_{\bq}(\bu)$
has a variable outside $U_b$, since there is no word with content
contained in $U_b$. If its first such variable differs from $b$, its
value is $3$. If that first variable is $b$, it must be the last
position of the word, since the first alternative of condition~\ref{snfb:cond:nonfinal} fails; the word has
value $1$. The sum on the right is at most $1$, again a contradiction.
This proves necessity. Condition~\ref{snfb:cond:ordinary} also implies
\eqref{snfb:eq:coverage}.

Conversely, assume both conditions and take an arbitrary assignment
into $H_5$. It suffices to use the summands in $D_{\bq}(\bu)$, since adjoining
summands can only increase a sum. If $\bq$ contains a factor with value
$4$, its value is the least element and the inequality is automatic.
If all variables of $\bq$ have value $2$, every word in the nonempty set
$D_{\bq}(\bu)$ also has value $2$.

In the remaining case, let $b$ be the variable at the first position
of $\bq$ whose value differs from $2$. All variables in $U_b$ have value
$2$. If $b$ has value $0$ or $3$, a summand $\br b\bs$ satisfying
condition~\ref{snfb:cond:ordinary} has the same value as $\bq$, by
Lemma~\ref{snfb:lem:evaluation}.

Suppose $b$ has value $1$. If this is the last position of $\bq$, then
$\bq$ has value $1$, and a summand satisfying condition~\ref{snfb:cond:ordinary}
has value $1$ or $0$. Otherwise $\bq$ has value $0$. A summand from the
first alternative of condition~\ref{snfb:cond:nonfinal} also has value $0$.
If instead that
condition provides $\bt$ with $c(\bt)\subseteq U_b$, then $\bt$ has
value $2$, while a summand satisfying condition~\ref{snfb:cond:ordinary} has value $1$ or $0$.
Their sum is $0$, because $1+2=0$. This proves sufficiency.
\end{proof}

\subsection{A finite equational basis for $H_5$}

The following are twelve individual identities, written using the
inequality abbreviation. Every letter in them is an ordinary variable.
\begin{align}
 x^2y&\approx xy,\label{snfb:eq:W1}\\
 xyx&\approx xy^2,\label{snfb:eq:W2}\\
 x&\preceq x^2,\label{snfb:eq:L}\\
 xyz&\preceq xz+yz,\label{snfb:eq:J}\\
 xyzt&\preceq xzt+yt,\label{snfb:eq:I}\\
 (x^2y)^2&\preceq x^2+y,\label{snfb:eq:E}\\
 (xvuyt)^2&\preceq xvut^2+uy+v,\label{snfb:eq:G}\\
 (xvyt)^2&\preceq xvt^2+y+v,\label{snfb:eq:G0}\\
 (xabzcy)^2&\preceq xy+az+bzc,\label{snfb:eq:A}\\
 (xbzcy)^2&\preceq xy+z+bzc,\label{snfb:eq:A0}\\
 (xazcy)^2&\preceq xy+az+zc,\label{snfb:eq:A1}\\
 (xzcy)^2&\preceq xy+z+zc.\label{snfb:eq:A00}
\end{align}
Let $\Sigma$ consist of these twelve identities and the six ai-semiring
axioms.

\begin{thm}\label{snfb:thm:Hbasis}
The finite set $\Sigma$ is an equational basis for $H_5$.
\end{thm}

\subsubsection{Validity of the identities}

\begin{lem}\label{snfb:lem:validity}
Every identity in $\Sigma$ holds in $H_5$.
\end{lem}
\begin{proof}
The semiring axioms hold by Lemma~\ref{snfb:lem:Haxioms}.
For \eqref{snfb:eq:W1}, the only element not fixed by squaring is $1$,
whose square is $0$, and $L_1=L_0$. For \eqref{snfb:eq:W2}, if either
variable has value $4$, both sides have value $4$. Otherwise:
when $x=0$ or $x=1$, both sides are $0$; when $x=3$, both are $3$;
and when $x=2$, the equality reduces to $y2=y^2$, which holds for
$y=0,1,2,3$. Identity~\eqref{snfb:eq:L} follows from $1\leq_{H_5}0$ and the
fact that the other four elements are idempotent.

For the other nine identities we apply Lemma~\ref{snfb:lem:criterion}.
Every summand on the right uses only variables of the lower word, and
the union of their contents is the content of that word. The table below
lists the summands required by the criterion. A pair $b:\bp$ in the middle
column means that $\bp$ satisfies condition~\ref{snfb:cond:ordinary} for $b$
and that the corresponding occurrence of $b$ in $\bp$ is nonfinal.
An entry $b:\bp;\bt$ in the last column gives a summand $\bp$ satisfying
condition~\ref{snfb:cond:ordinary}, together with a summand $\bt$ whose
content is contained in the set of variables preceding the first
occurrence of $b$ in the lower word. The word ``final'' means that
this first occurrence is the last position of the lower word, so
condition~\ref{snfb:cond:nonfinal} imposes no further requirement.

\begin{center}\small
\renewcommand{\arraystretch}{1.22}
\begin{tabular}{@{}p{.12\textwidth}p{.39\textwidth}p{.40\textwidth}@{}}
\hline
Identity & Nonfinal occurrences & Other occurrences\\
\hline
\eqref{snfb:eq:J} & $x:xz,\ y:yz$ & $z:xz$ (final)\\
\eqref{snfb:eq:I} & $x:xzt,\ y:yt,\ z:xzt$ & $t:yt$ (final)\\
\eqref{snfb:eq:E} & $x:x^2$ & $y:y;x^2$\\
\eqref{snfb:eq:G} & $x,v,u,t:xvut^2$ & $y:uy;v$\\
\eqref{snfb:eq:G0} & $x,v,t:xvt^2$ & $y:y;v$\\
\eqref{snfb:eq:A} & $x:xy,\ a:az,\ b,z:bzc$ & $c:bzc;az,\ y:xy;az$\\
\eqref{snfb:eq:A0} & $x:xy,\ b,z:bzc$ & $c:bzc;z,\ y:xy;z$\\
\eqref{snfb:eq:A1} & $x:xy,\ a:az,\ z:zc$ & $c:zc;az,\ y:xy;az$\\
\eqref{snfb:eq:A00} & $x:xy,\ z:zc$ & $c:zc;z,\ y:xy;z$\\
\hline
\end{tabular}
\end{center}
For example, in \eqref{snfb:eq:A} the first-occurrence order of the lower word is
$x,a,b,z,c,y$. For $c$, the summand $bzc$ satisfies
condition~\ref{snfb:cond:ordinary}, and the content of $az$ is contained
in $\{x,a,b,z\}$. For $y$, the summand $xy$ satisfies
condition~\ref{snfb:cond:ordinary}, while the same summand $az$ satisfies
the second alternative of condition~\ref{snfb:cond:nonfinal}. The other
rows are read in exactly the specified first-occurrence orders; thus
both conditions of the criterion have been checked for every variable
in every lower word. All nine identities follow.
\end{proof}

\subsubsection{Word normal forms}

For a nonempty word $\bw$, let $\init(\bw)$ list its distinct variables in
order of first occurrence. Thus $\init(\bw)$ is a word in which each
variable occurs once. We distinguish whether the last letter of $\bw$
has already appeared earlier in $\bw$.

The identities \eqref{snfb:eq:W1} and \eqref{snfb:eq:W2} are the two word
identities used in Lemma~\ref{lem:wordnormal}. Its proof therefore gives,
using only these two identities,
\begin{align*}
 \bw &\approx\init(\bw) &&\text{if the last letter occurs once},\\
 \bw &\approx\init(\bw)^2 &&\text{if the last letter has occurred earlier}.
\end{align*}
Consequently, words with the same first-occurrence order and the same
last-letter repetition condition are equivalent. In the second case,
$\bw$ and $\bw^2$ have the same first-occurrence order and a repeated last
letter, so $\bw^2\approx \bw$ also follows from these two identities.
As in Lemma~\ref{lem:wordnormal}, substitutions for variables are always
nonempty; only the surrounding contexts may be absent.

\subsubsection{Completeness when the last letter is repeated}

Suppose that $H_5\models\bq\preceq\bu$, where $\bu$ is a term.
By Lemma~\ref{snfb:lem:criterion},
$H_5\models\bq\preceq D_{\bq}(\bu)$, and
$D_{\bq}(\bu)\preceq\bu$ follows from the additive axioms.
It therefore suffices to derive the inequality with upper side
$D_{\bq}(\bu)$. Replacing $\bu$ by this term, we may assume that
every summand of $\bu$ has content contained in $c(\bq)$.
In particular,
\begin{equation}\label{snfb:eq:relevant}
 c(\bp)\subseteq c(\bq)\ (\bp \in \bu),\qquad
 \bigcup_{\bp \in \bu}c(\bp)=c(\bq).
\end{equation}
All deductions in the completeness argument below are from $\Sigma$.

Assume first that the last letter of $\bq$ appears earlier. Let $h$ be
its first letter. We construct a word below $\bu$ with the same
first-occurrence order and a repeated last letter.

\begin{lem}\label{snfb:lem:start}
There is a derivable inequality $\bw \preceq \bu$ such that $\bw$ starts with
$h$, has a repeated last letter, and $c(\bw)\subseteq c(\bq)$.
\end{lem}
\begin{proof}
Let $N$ be the set of variables occurring in a nonfinal position of
at least one word of $\bu$. Some word of $\bu$ ends in a letter of $N$.
For suppose otherwise. Then every word of $\bu$ ends outside $N$.
If $c(\bq)\subseteq N$, this contradicts \eqref{snfb:eq:relevant}.
Otherwise choose the first variable $b$ in the first-occurrence order
of $\bq$ that is outside $N$. All preceding variables lie in $N$.
Since $\bq$ has a repeated last letter, the first occurrence of $b$
is not final. Since $b\notin N$, no summand has the required
nonfinal occurrence of $b$.
The criterion therefore requires a word of $\bu$ whose content is
contained in the preceding variables, hence in $N$. This is impossible
because every word of $\bu$ ends outside $N$.

We may therefore choose words
\[
 \mathbf{a} z\in \bu,\qquad \mathbf{b} z\mathbf{c} \in \bu,\qquad \mathbf{c} \ne\varepsilon,
\]
where $\mathbf{a},\mathbf{b}$ may be empty contexts. At the first letter $h$ of $\bq$,
the preceding set is empty and cannot contain the content of a
nonempty word. The strengthened condition of the criterion supplies
a word $h\mathbf{y} \in \bu$ with $\mathbf{y}$ nonempty.

If $\mathbf{a},\mathbf{b}$ are nonempty, substitute these words in \eqref{snfb:eq:A}, with $x=h$.
This gives $(h\mathbf{a} \mathbf{b} z\mathbf{c} \mathbf{y})^2\preceq h\mathbf{y}+\mathbf{a} z+\mathbf{b} z\mathbf{c} \preceq \bu$. If $\mathbf{a}$ is empty,
use \eqref{snfb:eq:A0}; if $\mathbf{b}$ is empty, use \eqref{snfb:eq:A1}; if both are empty, use \eqref{snfb:eq:A00}.
In every case the left side is a square starting with $h$ and using
only variables of $\bq$. This is the required $\bw$.
\end{proof}

\begin{lem}\label{snfb:lem:extend}
The word $\bw$ in Lemma~\ref{snfb:lem:start} can be chosen with
$c(\bw)=c(\bq)$.
\end{lem}
\begin{proof}
For each $\bs \in \bu$, apply \eqref{snfb:eq:E} with $x=\bw$ and $y=\bs$:
\[
 (\bw^2\bs)^2\preceq \bw^2+\bs \approx \bw+\bs \preceq \bu.
\]
Here $\bw^2\approx \bw$ follows from Lemma~\ref{lem:wordnormal}.
The new word is still a square starting with $h$, and its content is
$c(\bw)\cup c(\bs)$. Iterate over the finite set $\bu$ and use
\eqref{snfb:eq:relevant}.
\end{proof}

We now put the first occurrences in the order prescribed by $\bq$.
Suppose that the current word has full content $c(\bq)$, has a
repeated last letter, and begins in first-occurrence order with a
nonempty correct prefix $\mathbf{a}$. Let $\mathbf{b}$ list the remaining variables in
their current first-occurrence order. If $\mathbf{b}$ is empty, the task is
complete. Otherwise $\mathbf{b}$ is a nonempty word, and the current word is
equivalent to $\mathbf{a} \mathbf{b}^2$ by Lemma~\ref{lem:wordnormal}. Let $d$ be the next
variable in the required order. If it is already the first letter
of $\mathbf{b}$, simply extend the correct prefix. In the other case proceed
as follows.

If there is a summand $\mathbf{r}  d \mathbf{s} \in \bu$ with the required nonfinal occurrence of $d$, where
$c(\mathbf{r})\subseteq c(\mathbf{a})$ and $\mathbf{s}$ nonempty, substitute
\[
 x=\mathbf{a},\quad y=\mathbf{r} d,\quad z=\mathbf{b},\quad t=\mathbf{s}
\]
in \eqref{snfb:eq:I}. All four substitutions are nonempty, even if the context $\mathbf{r}$
is empty. We obtain
\begin{equation}\label{snfb:eq:insert-nonfinal}
 \mathbf{a} \mathbf{r} d\mathbf{b} \mathbf{s} \preceq \mathbf{a} \mathbf{b} \mathbf{s}+\mathbf{r} d\mathbf{s} \preceq \bu.
\end{equation}
The word $\mathbf{a} \mathbf{b} \mathbf{s}$ has the same first-occurrence order as $\mathbf{a} \mathbf{b}$ and a
repeated last letter, because $\mathbf{a} \mathbf{b}$ already contains every variable
of $\bq$. It is therefore equivalent to the current word $\mathbf{a} \mathbf{b}^2$.
The word on the left of \eqref{snfb:eq:insert-nonfinal} has full content,
a repeated last letter, and correct first-occurrence prefix $\mathbf{a} d$.

If no such summand exists, the two conditions of the
criterion give
\[
 \mathbf{r} d\in \bu,\quad \bv \in \bu,\quad
 c(\mathbf{r}),c(\bv)\subseteq c(\mathbf{a}).
\]
The word $\bv$ is nonempty, while $\mathbf{r}$ may be empty. If $\mathbf{r}$ is nonempty,
substitute $x=\mathbf{a}$, $u=\mathbf{r}$, $y=d$, $t=\mathbf{b}$, and the displayed word $\bv$
for $v$ in \eqref{snfb:eq:G}. This gives
\[
 (\mathbf{a} \bv \mathbf{r} d\mathbf{b})^2\preceq \mathbf{a} \bv \mathbf{r} \mathbf{b}^2+\mathbf{r} d+\bv \preceq \bu.
\]
The first word on the right is equivalent to $\mathbf{a} \mathbf{b}^2$, because all
letters of $\bv,\mathbf{r}$ have already occurred in $\mathbf{a}$. The left side has
first-occurrence prefix $\mathbf{a} d$. If $\mathbf{r}$ is empty, \eqref{snfb:eq:G0} gives instead
\[
 (\mathbf{a} \bv d\mathbf{b})^2\preceq \mathbf{a} \bv \mathbf{b}^2+d+\bv \preceq \bu,
\]
with the same conclusion. Iteration increases the length of the correct
prefix until all variables are in the prescribed order. The resulting
word is equivalent to $\bq$ by Lemma~\ref{lem:wordnormal}. Thus $\bq \preceq \bu$
has been derived whenever the last letter of $\bq$ is repeated.

\subsubsection{Completeness when the last letter occurs once}

Now suppose that the last letter $z$ occurs only once in $\bq$.
If $H_5\models \bq^2\preceq \bu$, the preceding argument derives that
inequality, and \eqref{snfb:eq:L}, with $x=\bq$, gives $\bq \preceq \bq^2\preceq \bu$.

It remains to consider
\begin{equation}\label{snfb:eq:unique-case}
 H_5\models \bq \preceq \bu,\qquad H_5\not\models \bq^2\preceq \bu.
\end{equation}
The words $\bq,\bq^2$ have the same content and first-occurrence order.
The only extra requirement in Lemma~\ref{snfb:lem:criterion} for $\bq^2$
is the strengthened condition at $z$. Its failure has two consequences:
no word of $\bu$ has content contained in $c(\bq)\setminus\{z\}$,
and no word of $\bu$ has a nonfinal first occurrence of $z$ preceded
only by letters different from $z$. Hence every word of $\bu$ contains
$z$ exactly once, at its last position. We can write
\begin{equation}\label{snfb:eq:commonlast}
 \bu=\{\mathbf{u}_i z:i\in I\},\qquad z\notin c(\mathbf{u}_i),
\end{equation}
where a context $\mathbf{u}_i$ may be empty.

If $\bq=z$, then \eqref{snfb:eq:relevant} and \eqref{snfb:eq:commonlast} force
$\bu=\{z\}$ and the inequality is immediate. Otherwise let $h\ne z$
be the first letter of $\bq$. A summand satisfying
condition~\ref{snfb:cond:ordinary} for $h$ starts with $h$ and ends
with $z$. Start with this word
below $\bu$. If the current word is $\mathbf{r} z$, where $\mathbf{r}$ is nonempty and
does not contain $z$, and $\mathbf{u}_i$ is nonempty, \eqref{snfb:eq:J} yields
\[
 \mathbf{r} \mathbf{u}_i z\preceq \mathbf{r} z+\mathbf{u}_i z\preceq \bu.
\]
If $\mathbf{u}_i$ is empty, its summand $z$ introduces no new variable and can
be skipped. Iterating over $I$ gives a word below $\bu$ which starts
with $h$, has full content $c(\bq)$, and contains $z$ exactly
once, in its last position.

Use Lemma~\ref{lem:wordnormal} to remove repetitions before that last
position. Suppose the resulting first-occurrence order is $\mathbf{a} \mathbf{b} z$,
where $\mathbf{a}$ is a nonempty correct prefix and $\mathbf{b}$ lists the still
unplaced variables other than $z$. If $\mathbf{b}$ is empty we have obtained
a word equivalent to $\bq$. Otherwise let $d$ be the next variable in the required order;
as before, no deduction is needed if it already follows $\mathbf{a}$.
A summand satisfying condition~\ref{snfb:cond:ordinary} for $d$ has the form
\[
 \mathbf{r} d\mathbf{s} z\in \bu,\qquad c(\mathbf{r})\subseteq c(\mathbf{a}),
\]
where $\mathbf{r},\mathbf{s}$ may be empty and neither contains $z$. Apply \eqref{snfb:eq:I} with
\[
 x=\mathbf{a},\quad y=\mathbf{r} d,\quad z=\mathbf{b},\quad t=\mathbf{s} z.
\]
The word $\mathbf{s} z$ is nonempty.
It gives
\begin{equation}\label{snfb:eq:insert-final}
 \mathbf{a} \mathbf{r} d\mathbf{b} \mathbf{s} z\preceq \mathbf{a} \mathbf{b} \mathbf{s} z+\mathbf{r} d\mathbf{s} z\preceq \bu.
\end{equation}
The word $\mathbf{a} \mathbf{b} \mathbf{s} z$ is equivalent to the current word $\mathbf{a} \mathbf{b} z$: all letters
of $\mathbf{s}$ already occur in $\mathbf{a} \mathbf{b}$, while $z$ occurs just once at the end.
The left side has full content, still has a unique final $z$, and
has correct first-occurrence prefix $\mathbf{a} d$. Iteration gives the required
order, and Lemma~\ref{lem:wordnormal} then yields $\bq \preceq \bu$.

The substitution $t=\mathbf{s} z$ in \eqref{snfb:eq:insert-final} preserves the
unique occurrence of $z$ at the end of the lower word.

\begin{proof}[Completion of the proof of Theorem~\ref{snfb:thm:Hbasis}]
Lemma~\ref{snfb:lem:validity} proves that all identities in $\Sigma$ hold
in $H_5$. The two cases above derive every valid inequality of the form
$\bq \preceq \bu$, with $\bq$ a word and $\bu$ a nonempty sum of words, for
arbitrary finite sets of variables.

Finally, expand any identity of $H_5$ into sums of words using
distributivity. An equality
$\sum_i \bp_i\approx\sum_j \bq_j$ is equivalent, in an ai-semiring,
to all inequalities $\bp_i\preceq\sum_j \bq_j$ and
$\bq_j\preceq\sum_i \bp_i$. Each has just been shown derivable from
$\Sigma$, and their sums give the original equality. Thus $\Sigma$
is a finite equational basis for $H_5$.
\end{proof}

\subsection{A homomorphic image of a subsemiring of $H_5^2$}

In the direct square $H_5^2$, with coordinatewise operations, put
\begin{align*}
 b_0&=(0,0),&b_1&=(0,2),&b_2&=(0,4),&b_3&=(1,0),\\
 b_4&=(1,2),&b_5&=(2,0),&b_6&=(2,2),&b_7&=(3,4),
\end{align*}
and let $B_8=\{b_0,b_1,\ldots,b_7\}$. The operation tables induced
from $H_5^2$ are given in Table~\ref{tab:B8}. An entry $j$ denotes $b_j$.

\begin{table}[ht]
\caption{The Cayley tables of $B_8$}\label{tab:B8}
\setlength{\tabcolsep}{3.3pt}
\begin{tabular}{c|cccccccc}
+&0&1&2&3&4&5&6&7\\\hline
0&0&0&0&0&0&0&0&0\\
1&0&1&1&0&1&0&1&1\\
2&0&1&2&0&1&0&1&2\\
3&0&0&0&3&3&0&0&3\\
4&0&1&1&3&4&0&1&4\\
5&0&0&0&0&0&5&5&5\\
6&0&1&1&0&1&5&6&6\\
7&0&1&2&3&4&5&6&7
\end{tabular}
\qquad
\begin{tabular}{c|cccccccc}
$\cdot$&0&1&2&3&4&5&6&7\\\hline
0&0&0&2&0&0&0&0&2\\
1&0&1&2&0&1&0&1&2\\
2&2&2&2&2&2&2&2&2\\
3&0&0&2&0&0&0&0&2\\
4&0&1&2&0&1&0&1&2\\
5&0&0&2&3&3&5&5&7\\
6&0&1&2&3&4&5&6&7\\
7&7&7&7&7&7&7&7&7
\end{tabular}
\end{table}

For example, $b_5b_4=(2,0)(1,2)=(1,0)=b_3$, and
$b_4+b_5=(1,2)+(2,0)=(0,0)=b_0$. Every displayed entry is in
$B_8$, proving that $B_8$ is a subsemiring of $H_5^2$.

Define $\pi:B_8\to S$ by
\begin{equation}\label{snfb:eq:pi}
\begin{array}{c|cccccccc}
t&b_0&b_1&b_2&b_3&b_4&b_5&b_6&b_7\\\hline
\pi(t)&0&0&0&0&1&2&3&4
\end{array}.
\end{equation}

\begin{lem}\label{snfb:lem:quotient}
The map $\pi$ is a surjective homomorphism. In particular,
$S\cong B_8/\ker\pi$.
\end{lem}
\begin{proof}
Surjectivity is immediate from \eqref{snfb:eq:pi}. We verify that the
five fibres
\[
 F_0=\{b_0,b_1,b_2,b_3\},\quad F_1=\{b_4\},\quad
 F_2=\{b_5\},\quad F_3=\{b_6\},\quad F_4=\{b_7\}
\]
give exactly the operations of $S$ displayed in Section~\ref{sec:nfb}.

For addition, Table~\ref{tab:B8} shows that
$F_0+F_i\subseteq F_0$ for every $i$. Among the singleton fibres,
the diagonal sums stay in the same fibre, and the off-diagonal sums are
\[
 F_1+F_2,F_1+F_3\subseteq F_0,\qquad
 F_1+F_4\subseteq F_1,\qquad
 F_2+F_3,F_2+F_4\subseteq F_2,\qquad
 F_3+F_4\subseteq F_3.
\]
Together with commutativity these are all cases, and they reproduce
the addition table of $S$.

For multiplication, the first four rows of Table~\ref{tab:B8} give
$F_0F_i\subseteq F_0$ for all $i$, and the row of $b_4$ gives
$F_1F_i\subseteq F_0$ for all $i$. The remaining three rows give
\begin{align*}
 &(F_2F_0,F_2F_1,F_2F_2,F_2F_3,F_2F_4)
       \subseteq(F_0,F_0,F_2,F_2,F_4),\\
 &(F_3F_0,F_3F_1,F_3F_2,F_3F_3,F_3F_4)
       \subseteq(F_0,F_1,F_2,F_3,F_4),\\
 &F_4F_i\subseteq F_4\qquad(0\leq i\leq4),
\end{align*}
where inclusion between the two tuples means coordinatewise inclusion.
These are precisely the five rows of the multiplication table of $S$.
Thus $\pi(s+t)=\pi(s)+\pi(t)$ and $\pi(st)=\pi(s)\pi(t)$ for
every $s,t\in B_8$. The kernel is consequently a congruence and the
displayed quotient is isomorphic to $S$.
\end{proof}

\begin{pro}\label{snfb:prop:not-strong}
The ai-semiring $S$ belongs to $\Var(H_5)$ and is not strongly
nonfinitely based.
\end{pro}
\begin{proof}
The algebra $H_5$ is a five-element ai-semiring by
Lemma~\ref{snfb:lem:Haxioms}, and is finitely based by
Theorem~\ref{snfb:thm:Hbasis}. The subset $B_8$ is an eight-element
subsemiring of $H_5^2$, and Lemma~\ref{snfb:lem:quotient} exhibits $S$ as
its homomorphic image. Identities are preserved under direct products, subalgebras and
homomorphic images. Therefore
\[
 S\in\mathsf{HS}(H_5^2)\subseteq\Var(H_5),
 \qquad\text{and hence}\qquad \Var(S)\subseteq\Var(H_5).
\]
Since $H_5$ is finite and finitely based, the definition of strong
nonfinite basability implies that $S$ is not strongly nonfinitely based.
\end{proof}

\section{The least order of a finite generator}\label{sec:minimality}

In this section we prove that no ai-semiring of order at most four
can generate a limit variety. All algebras considered here have the
binary signature $\{+,\cdot\}$, without constant symbols. The order in
this statement is the number of elements of the generating semiring;
it is not the number of members of the generated variety.

We use the complete classification of semirings of order four obtained
in the preprint of Wang, Wang, Yin and Wu~\cite{Wang2026}.
Their Theorem~1.1 and Proposition~12.2 show that, up to isomorphism,
there are $866$ ai-semirings of order four, of which $821$ are finitely
based and $45$ are nonfinitely based. The latter representatives are
listed in their Theorem~11.5. The inclusions of varieties needed below
are proved explicitly from their operation tables.

\begin{lem}\label{lem:min-S7}
The variety generated by the three-element ai-semiring $S_7$ has a
nonfinitely based proper subvariety.
\end{lem}
\begin{proof}
Write $S_7=\{z,e,a\}$. Its addition is given by $u+u=u$ and
$u+v=z$ for $u\ne v$. For multiplication, $z$ is a zero, $e$ is an
identity element, and $a^2=z$. These names do not introduce constant
symbols into the signature.

Let $K=S_c(abc)$. More explicitly, the nonzero elements of $K$ are
\[
 a,\ b,\ c,\ ab,\ ac,\ bc,\ abc,
\]
regarded as the nonempty subsets of $\{a,b,c\}$. Addition is flat: two distinct elements have sum $0$.
The product of two nonzero elements is their union if their subsets
are disjoint, and is $0$ otherwise; $0$ is a multiplicative zero.
By \cite[Theorem~3.6]{RenJackson2023}, $\Var(K)$ is a limit variety.
In particular, it is nonfinitely based.

For completeness, we give the finite construction showing that
$K\in\Var(S_7)$. In the direct power $S_7^3$, put
\[
 D=S_7^3\setminus\{(e,e,e)\},\qquad
 I=\{(s_1,s_2,s_3)\in D:s_j=z\text{ for some }j\}.
\]
In $S_7$, an equality $u+v=e$ forces $u=v=e$, and so does an
equality $uv=e$. Hence the omitted tuple cannot be the sum or product
of two tuples of $D$. Thus $D$ is a subsemiring of $S_7^3$.
Since $z$ absorbs both addition and multiplication, we have
\[
 I+D\subseteq I,\qquad ID\subseteq I,\qquad DI\subseteq I.
\]
Consequently, the relation which identifies all elements of $I$ and
leaves every element of $D\setminus I$ in a singleton class is a
congruence on $D$.

Every tuple in $D\setminus I$ has its coordinates in $\{e,a\}$ and
has at least one coordinate equal to $a$. Identify it with the nonempty
set of positions at which $a$ occurs. Two distinct such tuples have
sum in $I$, because they differ at a coordinate where their sum is
$e+a=z$. Their product is outside $I$ exactly when their sets of
$a$-coordinates are disjoint, and then corresponds to the union of
these sets. Thus the quotient of $D$ by the displayed congruence is
isomorphic to $K$. It follows that $\Var(K)\subseteq\Var(S_7)$.

Every element of $K$ has square $0$, so $K$ satisfies $x^2\approx y^2$.
This identity fails in $S_7$, since $e^2=e\ne z=a^2$. Therefore
\[
 \Var(K)\subsetneq\Var(S_7),
\]
and the proper subvariety on the left is nonfinitely based.
\end{proof}

\begin{lem}\label{lem:min-four}
If $A$ is a nonfinitely based ai-semiring of order four, then
$S_7\in\Var(A)$.
\end{lem}
\begin{proof}
For this proof, let $C_i$ denote the algebra called $A_i$ in the local
catalogue of \cite[Appendix~B]{Wang2026}. Its underlying set is
$\{0,1,2,3\}$, with the operations printed there. These indices are
not the indices $S_{(4,j)}$ used in the four-element classification
series, and they are unrelated to the algebra $A_4$ used earlier in
this paper.

The $45$ nonfinitely based representatives in
\cite[Theorem~11.5]{Wang2026} are precisely the $41$ representatives
listed in Table~\ref{tab:min-direct} and the four representatives
\[
 C_{1442},\quad C_{1886},\quad C_{1892},\quad C_{1950}.
\]
We give a subsemiring and a surjective homomorphism onto $S_7$ for
each of them, allowing a direct square for these four exceptions.

\begin{table}[ht]
\caption{Subsemirings with a homomorphism onto $S_7$}
\label{tab:min-direct}
\begin{tabular}{p{.50\textwidth}cc}
Indices $i$ & Ordered domain & Images\\
\hline
$450,452,456,460,469,471,474,$
$1354,1359,1364,1368,1385,1387,$
$1390,1394,1398,1406,1409$
& $(0,1,2)$ & $(z,e,a)$\\[3pt]
$1407,1436,1440,1466,1475,1503,$
$1505,1511,1856,1858,1885,1887,$
$1890,1930,1937,1938,1953,1958$
& $(0,1,3)$ & $(z,e,a)$\\[3pt]
$1411,1469$ & $(0,1,2,3)$ & $(z,e,z,a)$\\[3pt]
$1476$ & $(0,1,2,3)$ & $(z,e,a,z)$\\[3pt]
$1769,1770$ & $(0,1,2)$ & $(e,z,a)$
\end{tabular}
\end{table}

Each row of Table~\ref{tab:min-direct} applies to every index in its
first column. The domain column lists the elements of a subsemiring
of $C_i$, and the last column lists their images in the same order.
Here is the verification from the operation tables. In the first,
second and fifth rows, after the displayed relabelling, distinct
elements have sum $z$, $e$ is a multiplicative identity, $z$ is a
multiplicative zero, and $a^2=z$. Thus these domains are isomorphic
to $S_7$. In the third row, the classes are
\[
 \{0,2\},\quad\{1\},\quad\{3\},
\]
with respective images $z,e,a$; in the fourth row they are
\[
 \{0,3\},\quad\{1\},\quad\{2\},
\]
with the same respective images. In each case the sum of distinct
classes is the $z$-class, the $z$-class absorbs multiplication,
the $e$-class is an identity for multiplication on the classes, and
the square of the $a$-class is the $z$-class. These assertions follow
entry by entry from the displayed operations of $C_{1411}$,
$C_{1469}$ and $C_{1476}$ in \cite[Appendix~B]{Wang2026}.
They show both that the indicated relation is a congruence and that
the quotient is $S_7$. This proves $S_7\in\Var(C_i)$ for all
$41$ indices in Table~\ref{tab:min-direct}.

For each of the four remaining indices, take the five tuples listed
in Table~\ref{tab:min-squares}, in the given order, and denote them by
$d_0,d_1,d_2,d_3,d_4$. All operations on tuples are coordinatewise.

\begin{table}[ht]
\caption{Five-element subsemirings of the exceptional squares}
\label{tab:min-squares}
\begin{tabular}{c|ccccc}
$i$ & $d_0$ & $d_1$ & $d_2$ & $d_3$ & $d_4$\\
\hline
$1442$ & $(1,2)$ & $(3,0)$ & $(3,2)$ & $(0,0)$ & $(0,2)$\\
$1886$ & $(0,1)$ & $(2,3)$ & $(0,3)$ & $(2,0)$ & $(0,0)$\\
$1892$ & $(0,1)$ & $(2,3)$ & $(0,3)$ & $(2,0)$ & $(0,0)$\\
$1950$ & $(1,2)$ & $(3,1)$ & $(3,2)$ & $(0,1)$ & $(0,2)$
\end{tabular}
\end{table}

In all four cases the two operations on these tuples are exactly those
in Table~\ref{tab:min-common}, where an entry $j$ denotes $d_j$.
In particular, the tables prove closure under both operations.
For example, for $i=1442$ the coordinate calculations give
\[
 (1,2)+(3,0)=(3,2),\qquad (3,0)(3,0)=(0,0),
\]
which are the entries $d_0+d_1=d_2$ and $d_1^2=d_3$.
The other entries are obtained in the same way from the four
operation tables in \cite[Appendix~B]{Wang2026}.

\begin{table}[ht]
\caption{The common operations on the five tuples}
\label{tab:min-common}
\begin{tabular}{c|ccccc}
$+$ & $0$ & $1$ & $2$ & $3$ & $4$\\
\hline
$0$ & $0$ & $2$ & $2$ & $4$ & $4$\\
$1$ & $2$ & $1$ & $2$ & $3$ & $4$\\
$2$ & $2$ & $2$ & $2$ & $4$ & $4$\\
$3$ & $4$ & $3$ & $4$ & $3$ & $4$\\
$4$ & $4$ & $4$ & $4$ & $4$ & $4$
\end{tabular}\qquad
\begin{tabular}{c|ccccc}
$\cdot$ & $0$ & $1$ & $2$ & $3$ & $4$\\
\hline
$0$ & $0$ & $1$ & $2$ & $3$ & $4$\\
$1$ & $1$ & $3$ & $3$ & $3$ & $3$\\
$2$ & $2$ & $3$ & $4$ & $3$ & $4$\\
$3$ & $3$ & $3$ & $3$ & $3$ & $3$\\
$4$ & $4$ & $3$ & $4$ & $3$ & $4$
\end{tabular}
\end{table}

In these common tables, $J=\{d_2,d_3,d_4\}$ is closed under addition
with every element and under multiplication on both sides by every
element. Collapsing $J$ therefore defines a congruence. The quotient
has three elements: $[d_0]$, $[d_1]$ and $J$. The class $[d_0]$ is a
multiplicative identity, $J$ is a multiplicative zero, and
$[d_1]^2=J$. The sum of any two distinct classes is $J$. Hence the map
\[
 d_0\longmapsto e,\qquad d_1\longmapsto a,\qquad
 d_2,d_3,d_4\longmapsto z
\]
is a surjective homomorphism onto $S_7$. This proves
$S_7\in\Var(C_i)$ for the four exceptional indices as well.
The classification now proves the lemma.
\end{proof}

\begin{thm}\label{thm:minimum-order}
No ai-semiring of order at most four generates a limit variety.
\end{thm}
\begin{proof}
Let $A$ be an ai-semiring of order at most four. If $A$ is finitely
based, then $\Var(A)$ is finitely based and is not a limit variety.
Suppose that $A$ is nonfinitely based.

The classification through order three, recalled in
\cite[Corollary~5.1]{Jackson2022}, shows that if $|A|\leq3$, then
$A$ is isomorphic to $S_7$. Lemma~\ref{lem:min-S7} excludes this case.
If $|A|=4$, Lemma~\ref{lem:min-four} gives
$\Var(S_7)\subseteq\Var(A)$. With the notation of
Lemma~\ref{lem:min-S7}, we have
\[
 \Var(K)\subsetneq\Var(S_7)\subseteq\Var(A).
\]
Thus $\Var(K)$ is a nonfinitely based proper subvariety of
$\Var(A)$, so $\Var(A)$ is not a limit variety. This completes the
proof.
\end{proof}

\begin{remark}
The restriction to the binary signature is essential. In the signature
$\{+,\cdot,1\}$ with a distinguished multiplicative identity,
$S_7$ does generate a limit variety; see
\cite[Section~4.3]{RenJackson2023}. Theorem~\ref{thm:minimum-order}
and the five-element example in this paper concern the signature
$\{+,\cdot\}$.
\end{remark}

\section{Proof of the main theorem}
\label{sec:conclusion}

\begin{proof}[Proof of Theorem~\ref{thm:main}]
Lemma~\ref{lem:S-ai} verifies that the two tables define a five-element
ai-semiring $S$. Theorem~\ref{thm:S-nfb} proves that
$\mathcal V=\Var(S)$ is nonfinitely based.

By Theorem~\ref{thm:two-coatoms}, $\mathcal V$ has exactly two maximal
proper subvarieties,
\[
 \mathcal C=\Var(M_0,M_1,C_4),\qquad
 \mathcal E=\mathcal V\cap\Mod(xy\preceq x+xyz),
\]
and every proper subvariety is contained in one of them.
Theorems~\ref{thm:C-fb} and~\ref{thm:E-fb} give finite bases for them.
Equation~\eqref{eq:D-HFB-conclusion} and Theorem~\ref{thm:j0-hfb} show that
$\mathcal D=\Var(M_1,C_4)$ and $\mathcal J_0=\Var(M_0,C_4)$
are hereditarily finitely based. Proposition~\ref{prop:C-intersections}
gives, for every $\mathcal W\subseteq\mathcal C$,
\[
 \mathcal W=(\mathcal W\cap\mathcal D)\vee
             (\mathcal W\cap\mathcal J_0).
\]
Lemma~\ref{lem:hfb-compactness} therefore proves that $\mathcal C$ is
hereditarily finitely based. Proposition~\ref{prop:E-transfer-HFB}
proves the same conclusion for $\mathcal E$.

Every proper subvariety of $\mathcal V$ is consequently finitely based.
Together with nonfinite basability of $\mathcal V$, this proves that
$\mathcal V$ is a limit variety.

By Theorem~\ref{snfb:thm:Hbasis}, the ai-semiring $H_5$ of order five
is finitely based. Lemma~\ref{snfb:lem:quotient} gives an eight-element
subsemiring $B_8$ of $H_5^2$ and a surjective homomorphism
$B_8\to S$. Hence $S\in\Var(H_5)$, and $S$ is not strongly
nonfinitely based.

By Theorem~\ref{thm:minimum-order}, an ai-semiring generating a limit
variety has at least five elements. Since the algebra $S$ has five
elements and generates a limit variety, we obtain
\[
 \min\{\,|A|: A\text{ is a finite ai-semiring and }\Var(A)
                    \text{ is a limit variety}\,\}=5.
\]
This completes the proof.
\end{proof}

\begin{remark}
A Cross variety is finitely based, generated by a finite algebra, and
has only finitely many subvarieties. A limit variety is nonfinitely
based, so $\Var(S)$ is not a Cross variety.
Moreover, $\Var(S)\subsetneq\Var(H_5)$, since the former is
nonfinitely based and the latter is finitely based. Thus the finite basis
property of the proper subvarieties of $\Var(S)$ is compatible with
its containment in a variety generated by a finite, finitely based
ai-semiring.
\end{remark}

\end{document}